\documentclass[11pt,a4paper]{article}

\usepackage{amsmath,amssymb}
\newcommand{\D}{\Delta}
\newcommand{\NN}{\frac{N}{p}}

\newcommand{\e}{\epsilon}

\newcommand{\va}{\varphi}
\newcommand{\n}{\nabla}
\newcommand{\w}{\omega}
\newcommand{\N}{\frac{N}{2}}
\newcommand{\NNN}{\frac{N}{p_{1}}}

\newcommand{\p}{\partial}

\newcommand{\R}{\mathbb{R}}
\newcommand{\T}{\mathbb{T}}
\newcommand{\h}{\hookrightarrow}

\newtheorem{definition}{Definition}
\newtheorem{theorem}{Theorem}

\newtheorem{proposition}{Proposition}

\newtheorem{remarka}{Remark}

\newtheorem{lemme}{Lemma}

\title{Regularity of weak solutions of the compressible barotropic Navier-Stokes equations}
\author{Boris Haspot \footnote{Karls Ruprecht Universit\"at Heidelberg, Institut for Applied Mathematics, Im Neuenheimer Feld 294,
D-69120 Heildelberg, Germany. Tel. 49(0)6221-54-6112} \footnote{Basque Center for Applied Mathematics (BCAM), Bizkaia Technology Park, Building 500. 48160 Derio Basque Country (Spain). 
} \footnote{Ceremade UMR CNRS 7534
Universit\'e de Paris Dauphine,
Place du MarŽchal DeLattre De Tassigny
75775 PARIS CEDEX 16 , haspot@ceremade.dauphine.fr }}
\date{}
\begin{document}
\maketitle
\begin{abstract}
Regularity and uniqueness of weak solutions of the compressible barotropic Navier-Stokes equations with constant viscosity coefficients is proven for small time in dimension $N=2,3$ under periodic boundary conditions. In this paper, the initial density is not required to have a positive lower bound and the pressure law is assumed to satisfy a condition that reduces to  $P(\rho)=a\rho^{\gamma}$ with $\gamma>1$ (in dimension three, additional conditions of size will be ask on $\gamma$). The second part of
the paper is devoted to  blow-up criteria for slightly subcritical initial data for the scaling of the equations when the viscosity coefficients  $(\mu,\lambda)$ are assumed constant
provided that their ratio is large enough (in particular $0<\lambda<\frac{5}{4}\mu$). More precisely we prove that under the condition $\rho$  belongs to  $L^{\infty}((0,T)\times\T^{N})$ then we can extend the unique solution beyond $T>0$.
Finally, we prove  that weak solutions in the torus $\mathbb{T}^{N}$ turn out to be smooth as long as the density remains bounded in $L^{\infty}(0,T,L^{(N+1+\e)\gamma}(\mathbb{T}^{N}))$ with $\e>0$ arbitrary small. This result may be considered as a Prodi-Serrin theorem (see \cite{prodi} and \cite{serrin}) for compressible Navier-Stokes system.
\end{abstract}
\section{Introduction}
The Navier-Stokes are the basic model describing the evolution of a viscous compressible gas.
Let us first recall that the periodic compressible barotropic Navier-Stokes equations in the torus
$\T^{N}$ ($N\geq2$) read as follows:
\begin{equation}
\begin{cases}
\begin{aligned}
&\p_{t}\rho+{\rm div}(\rho u)=0,\\
&\p_{t}(\rho u)+{\rm div}(\rho u\otimes u)-{\rm
div}(2\mu D(u))-\n(\lambda {\rm div}u)
+\n P(\rho)=\rho g.
\end{aligned}
\end{cases}
\label{1}
\end{equation}
The unknowns $\rho\geq0$ and  $u\in\mathbb{R}^{N}$  correspond to the density of the gas  and its velocity
field, respectively. The last equation of (\ref{1}) involves the pressure $P$ which is assumed to be a given increasing function of  $\rho$ and $D(u) = \frac{1}{2}[\n u + ^{t}\n u]$  is the strain tensor.
Recall that in the barotropic case, we have $P(\rho)=P_{\gamma,a}=a\rho^{\gamma}$ for some positive constant $a$ and some $\gamma\geq1$. The viscosity coefficients are assumed to satisfy $\mu>0$, and the physical case where $N\lambda+2\mu>0$, and the external force $g$ belongs to $L^{2}((0,T)\times\T^{N})^{N}$ for all $T>0$. Finally we complement the above system with the following initial conditions
\begin{equation}
\begin{cases}
\begin{aligned}
&\rho_{|t=0}=\rho_{0}\geq 0,\\
&\rho u_{|t=0}=m_{0}.
\end{aligned}
\end{cases}
\label{2}
\end{equation}
As emphasized in many papers related to compressible fluid dynamics \cite{10,11,16,19,20}, vacuum is a major difficulty when trying to prove global existence and strong regularity results. A second main difficulty corresponds to get estimate in $L^{\infty}$ norm for the density in order to control the non-linear term as the pressure but also for preserving the regularity of the velocity by assuming that the density stays in some multiplier spaces (for more details on this point we refer to \cite{H3}).\\
As a matter of fact, starting from bounded initial densities that have positive lower bounds, local existence of smooth solutions can be proved by classical means, since lower and upper bounds on the density persist for small enough time. To obtain global existence of strong solution, one of the main difficulty corresponds to the loss of control on the vacuum and on the $L^{\infty}$ norms of the density. An other point is related to the difficulty to estimate for long time the Lipschitz norm of the velocity $u$, that is why many criterion of blow-up assumed that the Lipschitz norm of the velocity $\|\n u\|_{L^{1}_{T}(L^{\infty})}$ is bounded (in particular in this case it allows  to control easily the $L^{\infty}$ norm of the velocity).
In this paper, we want to follow an radically different approach inasmuch as we want to define blow-up conditions depending only on the Lebesgue integrability of the density.  More precisely we want to prove that the norm $L^{\infty}(L^{q}(\T^{N}))$ on the pressure $P(\rho)$ controls the breakdown of strong solutions of the Navier-Stokes equations in dimension $N=2,3$ for large enough $q>1$. In other words, if a solution of the Navier-Stokes equations is initially suitably smooth and loses its regularity at some later time, then the norm $L^{\infty}(L^{q}(\T^{N}))$ of the pressure grows without bounds as the critical time approaches.\\
Before stating our main results, I would like to recall some important results concerning the theory of the existence of global weak solutions and the existence of strong solution (it means unique solutions)  for barotropic Navier-Stokes system. Indeed the results obtained in this article combine the different tools used for proving the existence of global weak solution and of strong solution. 
\subsection{Existence of global weak solutions}
Before the work by P-L Lions in \cite{13}, very little was known about solutions of the compressible barotropic Navier-Stokes equation at least when $N\geq 2$. There he proved a global existence theorem and weak stability results for $P_{\gamma,a}$ pressure laws under some conditions on $\gamma$ and with the following assumptions on the initial data:
\begin{equation}
\begin{cases}
\begin{aligned}
&\rho_{0}\in L^{1}(\T^{N})\cap L^{\gamma}(\T^{N}), \rho_{0}\geq0,\\
&\frac{m_{0}^{2}}{\rho_{0}}\in L^{1}(\T^{N}),
\end{aligned}
\end{cases}
\label{3}
\end{equation}
where we agree that $\frac{m_{0}^{2}}{\rho_{0}}=0$ on $\{x\in\T^{N}\;\mbox{such that}\;\rho_{0}(x)=0\}$.
Here is a up-to-date statement of Lions' result later improved by E. Feireisl et al in \cite{5F1,5F2,5F3}.
\begin{theorem}
\label{theorem1}
We assume (\ref{3}) and $\gamma>1$ if $N=2$, $\gamma>\frac{3}{2}$ if $N=3$. Then there exists a solution $(\rho,u)\in L^{\infty}(0,\infty;L^{\gamma}(\T^{N}))\times L^{2}(0,\infty;H^{1}(\T^{N}))^{N}$ satisfying in addition $\rho\in C([0,\infty),L^{p}(\T^{N}))$ if $1\leq p<\gamma$, $\rho|u|^{2}\in L^{\infty}(0,\infty;L^{1}(\T^{N}))$, $\rho\in L^{q}_{loc}([0,\infty);L^{q}(\T^{N}))$ for $1\leq q\leq \gamma-1+\frac{2}{N}\gamma$. Moreover, when $f=0$, for almost all $t\geq 0$, we have
\begin{equation}
\begin{aligned}
&\int_{\T^{N}}(\frac{1}{2}\rho|u|^{2}+\frac{a}{\gamma-1}\rho^{\gamma})(t,x)dx+\int^{t}_{0}ds
\int_{\T^{N}}(\mu|\n u|^{2}+(\lambda+\mu)({\rm div}u)^{2})dxds\\
&\hspace{7cm}\leq \int_{\T^{N}}(\frac{1}{2}\rho_{0}|u_{0}|^{2}+\frac{a}{\gamma-1}\rho_{0}^{\gamma})(x)dx.
\end{aligned}
\label{4}
\end{equation}
\end{theorem}
Lions  proved similar results for more general pressure laws $P(\rho)$ such that
\begin{equation}
\begin{cases}
\begin{aligned}
&\int^{1}_{0}\frac{P(s)}{s^{2}}ds<+\infty,\\
&\lim\inf\rightarrow_{s\rightarrow+\infty}\frac{P(s)}{s^{\gamma}}>0,
\end{aligned}
\end{cases}
\label{e4}
\end{equation}
for some $\gamma$ satisfying the above condition of theorem \ref{theorem1}.
\medbreak
Let us stress that  the main difficulty for proving Lions'  theorem consists
in exhibiting strong compactness properties of the density $\rho$ in
$L^{p}_{loc}(\R^{+}\times\R^{N})$ spaces required to pass to the
limit in the pressure
term $P(\rho)=a\rho^{\gamma}$.
As a matter of fact, in his pioneering work, Lions made
the additional assumption that  $\gamma-1+\frac{2\gamma}{N}\geq2$
 in dimension $N=2,3.$
 However, later on,  Feireisl and his collaborators in \cite{5F1,5F2,5F3} generalized the result to
any $\gamma>\frac{N}{2}$ for $N\geq 2$ by establishing that we can obtain renormalized
solution \emph{without} imposing that $\rho\in
L^{2}_{loc}(\R^{+}\times\R^{N})$ (see also \cite{NS}). This improvement
was based on  the concept of oscillation defect measure evaluating the
loss of compactness.
\subsection{Existence of unique solution}
The problem of existence of global solutions for Navier-Stokes equations was addressed in one dimension for
smooth enough data by Kazhikov and Shelukin in \cite{5K1}, and for
discontinuous ones, but still with densities away from zero, by Serre
in  \cite{16} and Hoff in \cite{5H1}.  Those results have been
generalized to higher dimension by Hoff  in \cite{5H2,5H3}.
The existence and uniqueness of local classical solutions for (\ref{1})
with smooth initial data such that the density $\rho_{0}$ is bounded
and bounded away from zero has been stated by Nash in \cite{5Na}. Let us emphasize that no stability condition was required there. On the other hand, for small smooth perturbations of a stable equilibrium with constant positive density, global well-posedness
has been proved in \cite{5MN}. Refined functional analysis has been used
for the last decades, ranging from Sobolev, Besov, Lorentz and
Triebel spaces to describe the regularity and long time behavior of
solutions to the compressible model \cite{5So}, \cite{20}, \cite{5H4}, \cite{5HZ}, \cite{5K1}.
\\
Guided in the approach by numerous works dedicated to the incompressible Navier-Stokes equation (see e.g \cite{Meyer}) we want here to recall the fundamental notion of \textit{critical} regularity.\\
By critical, we mean that we want to solve the system (\ref{1}) in functional spaces with norm
 invariant by the changes of scales which leaves (\ref{1}) invariant.
In the case of barotropic fluids, it is easy to see that the transformations:
\begin{equation}
(\rho(t,x),u(t,x))\longrightarrow (\rho(l^{2}t,lx),lu(l^{2}t,lx)),\;\;\;l\in\R,
\label{bbb1}
\end{equation}
have that property, provided that the pressure term has been changed accordingly.\\
The use of critical functional frameworks led to several new well-posedness results for compressible
fluids (see \cite{DL,DW, H1, H2}). In addition to have a norm invariant by (\ref{bbb1}),
appropriate functional spaces for solving (\ref{1}) must provide a control on the $L^{\infty}$
norm of the density (in order to avoid vacuum and loss of parabolicity) but also on the Lipschitz norm of the velocity (it means $\|\n u\|_{L^{1}_{T}(L^{\infty})}$) in order to be able to estimate the density via the  mass equation. For that reason,
the study is restricted to the case where the initial data $(\rho_{0},u_{0})$ and external force $f$
are in homogeneous Besov spaces such that, for some positive constant:
\begin{equation}
\rho_{0}\in L^{\infty}\cap B^{\N}_{p,\infty},\;u_{0}\in B^{\frac{N}{p_{1}}-1}_{p_{1},1}\;\;\mbox{and}\;\;f\in L^{1}_{loc}(\R^{+},\in B^{\frac{N}{p_{1}}-1}_{p_{1},1})
\label{critique}
\end{equation}
with $(p,p_{1})\in [1,+\infty[$ suitably chosen.
\subsubsection*{Local existence of strong solution with large initial data}
The first result concerning the existence of strong solutions in spaces invariant for the scaling of the equations is due to R. Danchin in
\cite{DL}. More precisely he obtains strong solutions for
initial data in $B^{\frac{N}{2}}_{2,1}\times (B^{\frac{N}{2}-1}_{2,1})^{N}$. Here compared with
the result on incompressible Navier-Stokes, he needs to control the vacuum hence the norm  of the density is in $L^\infty$ in order to take advantage of the parabolicity of the momentum equation.
That is why he is working with a third index $r=1$ for the previous Besov space because $B^{\N}_{2,1}$ is embedded in $L^{\infty}$. In \cite{DW}, R. Danchin generalizes the previous results by working with more general Besov space of the type $B^{\NN}_{p,1}\times (B^{\NN-1}_{p,1})^{N}$ with some restrictions on the choice of $p$. Indeed comparing with the general setting in (\ref{critique}), we can observe that R. Danchin needs to have $p=p_{1}$ leading to the limitation $p<2N$ for the existence and $p\leq N$ for the uniqueness due to some limitation concerning the paraproduct laws when he treats some non-linear terms. The fact that  $p=p_{1}$ is a consequence of the  strong coupling between the density and the velocity equations. To be more precise, in \cite{DW} the pressure term is considered as a remainder for the parabolic operator in the momentum equation of (\ref{1}).\\
In \cite{H3}, we address the question of local well-posedness in the critical functional framework under the assumption that the initial density and the initial velocity belong to critical Besov spaces with different index of integrability.  In this article we improve the results of \cite{DW} inasmuch as we have no restriction on the size of $p$. In particular, we can observe that when $p$ goes to infinity, we are close to get strong solution with initial data  $(q_{0},u_{0})$ in
$B^{0}_{\infty,1}\times B^{1}_{N,1}$ (at less when $P(\rho)=K\rho$ and weak solution else when the pressure are more general). It means that this theorem bridges the gap of the result of D. Hoff where the initial density is in $L^{\infty}$ but where the initial velocity is more regular (it means not critical) and the results of R. Danchin in \cite{DW}. \\
For proving this result, we adapt the spirit of the results of \cite{AP} and \cite{H} which treat the case of the density-dependent incompressible Navier-Stokes equations  (at the difference that in these works the velocity and the density are naturally decoupled).
The key of  \cite{H3} is to introduce a new unknown  $v_{1}$ \textit{the effective velocity} to avoid the coupling between the density and the velocity, we analyze then by a new way the pressure term. This idea  originates from the works by  D. Hoff,  P-L Lions and D. Serre (see \cite{5H4, 13, 16}) where the so-called  \textit{effective pressure} has been first introduced.  We refer also to the fundamental result of  V. A. Vaigant and A. V Kazhikhov on the existence of global strong solution in two dimension which use crucially this notion of \textit{ effective pressure} (in this paper the viscosity coefficients are variable and have a specific form) .\\
In \cite{H3} the divergence of this \textit{ effective velocity} is exactly the \textit{effective pressure}. More precisely we write the gradient of the pressure as a Laplacian of some vector-field $v$, and we introduce this term in the linear part of the momentum equation ( in other words, $v={\cal G}P(\rho)$ where
${\cal G}P(\rho)$ stands for some pseudo-differential operator of order $-1$). We then introduce the effective velocity $v_{1}=u-v$.  By this way, we have canceled out the coupling between $v_{1}$ and the density. We next verify easily that we have a Lipschitz control on the gradient of $u$ (it is crucial to estimate the density via the mass equation). 
\subsubsection*{Global existence of strong solution with small initial data}
In this subsection we would like to recall some important results of existence of global strong solutions with small initial data. More precisely we would like to emphasize the importance of the notion of \textit{effective pressure} (see \cite{13}, \cite{16}, \cite{5H4}) or of \textit{effective velocity} (see \cite{H4}) in the resolution of this problem. Indeed this last notion introduced in \cite{H3} will be the key of our proofs of blow-up, in some sense this unknown of \textit{effective velocity} is more regular than the velocity (we will give a precise definition of the  \textit{effective velocity} later) and we will explain why.\\
The global existence of strong solutions for initial data with high regularity order and close to a stable equilibrium has been proved by Matsumura and Nishida in \cite{5MN} for three-dimensional polytropic ideal fluids and no outer force with initial data such that $(\rho_{0}-1,u_{0})\in H^{3}(\R^{3})\times H^{3}(\R^{3})$. More recently D. Hoff in \cite{5H4,5H2,10} stated the existence of global weak solutions with small initial data including discontinuous initial data (namely $\rho_{0}-1$ is small in $L^{2}(\R^{N})\cap L^{\infty}(\R^{N})$ and $u_{0}$ is small in $L^{4}(\R^{2})$ if $N=2$ and small
in $L^{8}(\R^{3})$ if $N=3$). One of the major interest of the results of D. Hoff is to exhibit some smoothing effects on the incompressible part of the velocity $u$ and on \textit{the effective pressure} $F=(2\mu+\lambda){\rm div}u-P(\rho)+P(\bar{\rho})$ (see also the work by D. Serre in \cite{16}). This also plays a crucial role in the proof of P.-L. Lions for the existence of global weak solution (see \cite{13}). However if the results of D. Hoff are critical in the sense of the scaling for the density, it is not the case for the initial velocity. In \cite{5H2}, D. Hoff shows a very interesting theorem of weak-strong uniqueness when $P(\rho)=K\rho$ with $K>0$. To speak roughly under the conditions that two solutions $(\rho,u)$, $(\rho_{1},u_{1})$ satisfy a control $L^{\infty}$ on the density and a control Lipschitz on the velocity, with additional property of regularity on the strong solution $(\rho_{1},u_{1})$ then we obtain $(\rho,u)=(\rho_{1},u_{1})$.
D. Hoff uses this result to show that the solutions of \cite{5H2} are unique.\\
Finally R. Danchin in \cite{DG} shows for the first time a result of existence of global strong solution close to a stable equilibrium in critical space for the scaling of the system. More precisely the initial data are chosen as follows $(\rho_{0}-1,u_{0})\in (B^{\N}_{2,1}\cap B^{\N-1}_{2,1})\times B^{\N-1}_{2,1}$. The main difficulty is to get estimates on the linearized system given that the velocity and the density are coupled via the pressure. What is crucial in this work is the smoothing effect on the velocity and a $L^{1}$ decay on $\rho-1$ (this plays a key role to control the pressure term). In this work, R. Danchin uses some tricky energy inequalities on the system in Fourier variable. This explains in particular why the result is obtained in Besov space with a Lebesgue index $p=2$. Recently Q. Chen et al in \cite{CMZ} and F. Charve and R. Danchin in \cite{CD} improve the previous result by working in more general Besov space, for that they study the linear part of the system by including the convection terms in the linearized system. The main idea is then to "paralinearize" the convection terms, which avoids any troubles concerning the coupling of the linear system between low and high frequencies.\\
In \cite{H4} we make a connection between the article of D. Hoff \cite{5H4,5H2}
and those of  F. Charve and R. Danchin and Q. Chen et al in \cite{CD} and \cite{CMZ}. In fact we extend the results \cite{CD} and \cite{CMZ} to the case where the Lebesgue index of Besov spaces are not the same for the density and the velocity. To do that, as in \cite{H3} we introduce the unknown of \textit{effective velocity} in high frequencies so as to "kill" the relation of coupling between the velocity and the pressure. This effective velocity enables us to get as in R. Danchin in \cite{DG} a $L^{1}$ decay on $\rho-1$ in the high frequency regime. In low frequencies, the first order terms predominate, so that (\ref{1}) has to be treated by means of hyperbolic energy methods (roughly $\rho-1$ and the potential part of the velocity verify a wave equation). This implies that we can treat the low regime only in space constructed on $L^{2}(\R^{N})$ as it is classical that hyperbolic systems are ill-posed in general $L^{p}(\R^{N})$ spaces. So as in \cite{CMZ} and \cite{CD}, the system has to be handled differently in low and high frequencies. In particular in \cite{H4}, we are able to deal with very critical spaces of initial data, in particular $(\rho_{0}-1,u_{0})\in \widetilde{B}^{\N-1,0}_{2,\infty,1}\times \widetilde{B}^{\N-1,0}_{2,N,1}$ (see \cite{H4} for the definition of the hybrid Besov spaces).
\subsubsection*{Existence of strong solution with vacuum}
On the other hand there have been few existence results on the strong solutions for the general case of nonnegative initial densities. The first result was proved by R. Salvi and I. Straskraba. They showed in \cite{CCK17} that
if $\Omega$ is a bounded domain, $P=P(\cdot)\in C^{2}[0,\infty)$, $\rho_{0}\in H^{2}$, $u_{0}\in H^{1}_{0}\cap H^{2}$ and the compatibility condition:
\begin{equation}
Lu_{0}+\n P(\rho_{0})=\rho_{0}^{\frac{1}{2}}g,\;\;\mbox{for some}\;g\in L^{2},
\label{1.3}
\end{equation}
is satisfied, then there exists a unique local strong solution $(\rho,u)$ to the initial boundary value problem (\ref{1}). H. J. Choe and H. Kim proved in \cite{CCK3} a similar existence result when $\Omega$ is either a bounded domain or the whole space, $P(\rho)=a\rho^{\gamma}$ ($a>0$, $\gamma>1$), $\rho_{0}\in L^{1}\cap H^{1}\cap W^{1,6}$, $u_{0}\in D^{1}_{0}\cap D^{2}$ and the condition (\ref{1.3}) is satisfied.\\
B. Desjardins in \cite{CCK5} proved the local existence of a weak solution solution $(\rho,u)$ with a bounded nonnegative density to the periodic boundary value problem (\ref{1}) as long as $\sup_{0\leq t\leq T^{*}}(\|\rho(t)\|_{L^{\infty}(\T^{3})}+\|\n u(t)\|_{L^{2}(\T^{3})})<+\infty$. In dimension $N=2$, the regularizing effects proved in  \cite{CCK5} hold as long as $\sup_{0\leq t\leq T^{*}}(\|\rho(t)\|_{L^{\infty}(\T^{2})})<+\infty$. We would like to draw the attention on this last result because in the sequel we will generalize this assertion to the case of the dimension $N=3$ in the context of the continuation of strong solution.
\subsection{Notations and main results}
In a first time this paper will be devoted to improve the works \cite{CCK5} and \cite{CCK3} by choosing a larger class of initial velocity data. A crucial point will be to explain how it is possible to obtain strong solutions and not only weak solutions as in \cite{CCK5} and \cite{CCK3}. In the sequel we will give new results of blow-up which will be our main results.\\
In the sequel we will note $\frac{d}{dt}=\p_{t}+u\cdot\n$ and $\dot{f}=\frac{d}{dt}f$. We will define also the unknown $\omega={\rm curl}u$ as the rotational of the velocity.
The viscosity coefficients are constant and are assumed to satisfy:
\begin{equation}
\mu>0,\;\;0<\lambda<\frac{5}{4}\mu.
\label{condviscosity}
\end{equation}
It follows that there is a $p>6$, which will be fixed throughout, such that:
\begin{equation}
\frac{\mu}{\lambda}>\frac{(p-2)^{2}}{4(p-1)}.
\label{condviscosity1}
\end{equation}
In the sequel we will assume that $g\in E^{1}_{T}$ with:
\begin{equation}
\begin{aligned}
&\|g\|_{E^{1}_{T}}=\|g\|_{L^{\infty}_{T}(L^{2}(\T^{N}))}+\|g\|_{L^{2}_{T}(L^{2}(\T^{N}))}+\|g\|_{L^{1}_{T}(L^{N+\e}(\T^{N}))}\\
&\hspace{2cm}+\int^{T}_{0}f^{7}(s)\|\n g\|_{L^{4}(\T^{N})}^{2}ds+\int^{T}_{0}\int_{\T^{N}}f^{5}(s)|\p_{s}g|^{2}dsdx,
\end{aligned}
\label{condig}
\end{equation}
where $f(s)=\min(1,s)$ and $\e>0$.\\
Let us now state a first theorem of weak-strong solutions.
\begin{theorem}
\label{theo1}
Let $N=2,3$. Assume that $\mu$ and $\lambda$ verify (\ref{condviscosity}) and (\ref{condviscosity1}). We assume that $\rho_{0}\in L^{\infty}(\T^{N})$, $\rho_{0}^{\frac{1}{p}}u_{0}\in L^{p}(\T^{N})$ and $\rho_{0}^{\frac{1}{2}}u_{0}\in L^{2}(\T^{N})$. Moreover $g$ is in $E^{1}_{T}$ for any $T>0$ and $g\in L^{\infty}(L^{p}(\T^{N}))$. Finally we assume that $\gamma>\frac{p-1}{p-6}$ if $N=3$ and $\gamma>1$ if $N=2$. Here $p$ verifies:
\begin{equation}
 \begin{cases}
  \begin{aligned}
&p>2,\;\;\;\mbox{if}\;\;N=2,\\
&p>6,\;\;\;\mbox{if}\;\;N=3.
  \end{aligned}
\end{cases}
\label{1c24}
\end{equation}
\begin{itemize}
\item  There exists $T_{0}\in (0,+\infty]$ and a weak solution $(\rho,u)$ to the system (\ref{1}) in $[0,T_{0}]$ such that for all $T<T_{0}$ (with $f(t)=\min(t,1)$, $\omega={\rm curl}u$):
\begin{equation}
 \begin{aligned}
&\sup_{0<t\leq T_{0}}\int_{\T^{N}}[\frac{1}{2}\rho(t,x)|u(t,x)|^{2}+|\Pi(\rho(t,x))|+f(t)|\n u(t,x)|^{2}]dx\\
&+\sup_{0<t\leq T_{0}}\int_{\T^{N}}[\frac{1}{2}f(t)^{N}(\rho|\dot{u}(t,x)|^{2}+|\n \omega(t,x)|^{2})]dx\\
&+\int^{T_{0}}_{0}\int_{\T^{N}}[|\n u|^{2}+f(s)(\rho|\dot{u}|^{2}+|\n\omega|^{2})+f^{N}(s)|\n \dot{u}|^{2}]dxdt\leq C C_{0}.
 \end{aligned}
\label{1.21}
\end{equation}
\begin{equation}
\begin{aligned}
&\sup_{0<t\leq T_{0}}\int_{\T^{N}}\rho(t,x)\,|u(t,x)|^{p}dx\leq C_{0,g}
\end{aligned}
 \label{1.211}
\end{equation}
with:
$$\Pi(z)=\big(\int^{s}_{0}\frac{P(z)}{z^{2}}dz\big).$$
Furthermore $C_{0,g}$ depends only on the initial data $\rho_{0}$,  $u_{0}$ and on $g$.
\item In addition if we assume that $u_{0}\in H^{\frac{N}{2}-1+\e}(\T^{N})$ with $\e>0$ and $\frac{1}{\rho_{0}}\in L^{\infty}(\T^{N})$, we obtain the following estimates:
\begin{equation}
 \begin{aligned}
&\sup_{0<t\leq T_{0}}\int_{\T^{N}}[\frac{1}{2}\rho(t,x)|u(t,x)|^{2}+\Pi(\rho)(t,x)+t^{2-\N-\e}|\n u(t,x)|^{2}]dx\\
&+\sup_{0<t\leq T_{0}}\int_{\T^{N}}[\frac{1}{2}t^{\sigma}(\rho|\dot{u}(t,x)|^{2}+|\n \omega(t,x)|^{2})dx]\\
&+\int^{T_{0}}_{0}\int_{\T^{N}}[|\n u|^{2}+t^{2-\N-\e}|\dot{u}|^{2}+t^{\sigma}|\n \dot{u}|^{2}]dxdt\leq C C^{'}_{0,g},
 \end{aligned}
\label{b1.22}
\end{equation}
with  $C^{'}_{0,g}$ depends on the initial data and on $g$ with:
\begin{equation}
 \begin{cases}
  \begin{aligned}
&\sigma=2-\e,\;\;\;\mbox{if}\;\;N=2,\\
&\sigma=\frac{3}{2}-\e,\;\;\;\mbox{if}\;\;N=3,
  \end{aligned}
\end{cases}
\label{c24}
\end{equation}
and:
\begin{equation}
\sup_{0\leq t\leq T_{0}}\|u(t,\cdot)\|_{H^{\N-1+\e}(\T^{N})}\leq CC^{1}_{0,g},\;\;\mbox{and}\;\;
\frac{1}{\rho}\in L^{\infty}_{T_{0}}(L^{\infty}(\T^{N})),
 \label{b1.23}
\end{equation}
where $C^{1}_{0,g}$ depends on the initial data and on $g$ and:
\begin{equation}
\n u\in L^{1}_{T_{0}}(BMO(\T^{N})).
\label{c1.23}
\end{equation}
\item The regularity properties (\ref{1.21}),  (\ref{1.211}), (\ref{b1.22}), (\ref{b1.23}) and (\ref{c1.23}) hold as long as:
\begin{equation}
\sup_{t\in[0,T_{0}]}\|\rho\|_{L^{\infty}_{t}(L^{\infty}(\T^{N}))}<+\infty.
\label{11.23}
\end{equation}
\end{itemize}
\end{theorem}
\begin{remarka}
The crucial point of this theorem and his main interest is that estimates (\ref{1.21}),  (\ref{1.211}), (\ref{b1.22}), (\ref{b1.23}) and (\ref{c1.23}) can hold as long as the density $\rho$ belongs to $L^{\infty}_{T_{0}}(\T^{N})$ in dimension $2$ and $3$. It generalizes in particular the result of B. Desjardins in \cite{CCK5} where this result was obtained only in the case of the dimension $N=2$.
 The main ingredient to obtain our theorem compared with \cite{CCK5} is to benefit of the gain of integrability that we can obtain for compressible Navier-Stokes system when the pressure has enough integrability.\\
In the next theorem, we are going to construct strong solutions by adding a slight hypothesis of regularity on the initial data $\rho_{0}$. Furthermore we shall extend this strong solution under an hypothesis of control on $\rho$ in $L^{\infty}_{T_{0}}(\T^{N})$.\\
 In passing, we observe that in this theorem we improve also \cite{CCK5} inasmuch as we can choose more general initial data. Indeed in \cite{CCK5}, the initial data $u_{0}$ is assumed to belong to $H^{1}(\T^{N})$, here we do not ask any regularity on the velocity but only integrability.
 In fact this is possible because at the difference of \cite{CCK5}, we obtain energy inequalities by multiplying by $f(t)$ which allows us to reduce the regularity assumptions on the initial data. This idea was in particular introduce by Kato in \cite{Kato} for incompressible Navier-Stokes system and develop in the case of barotropic  Navier-Stokes system by D. Hoff in \cite{5H2}.
\end{remarka}
\begin{remarka}
However compared with \cite{CCK5} and \cite{H3} , our hypothesis on the viscosity coefficients and on the choice of $\gamma$ are more restrictive. Indeed in \cite{CCK5} B. Desjardins needs only to assume that $\gamma>3$ in dimension $N=3$. In fact in our case we have to add these conditions because at the difference of  \cite{CCK5} we will obtain only a control on $\sqrt{f(t)}P(\rho)$ in $L^{\infty}(L^{2})$ and not on  $P(\rho)$ in $L^{\infty}(L^{2})$ . That is why we need of additional integrability condition on the density to control in $L^{\infty}$ norm the term coming from $(\D)^{-1}{\rm div}(\rho u)$ (we refer to the proof for more details).\\
We want to mention that this condition on $\gamma$ is purely technic and do not play any role in the criterion of blow-up in theorem \ref{corollaire1} and \ref{theo3}. Indeed this condition is useful only to obtain in the theorem \ref{theo3} the fact that $\rho$ is in $L^{\infty}$. If we assume in a blow-up criterion that the density $\rho$ is $L^{\infty}$, then this technical assumption on $\gamma$ can be avoided.\\
Finally we have to impose certain conditions on the viscosity coefficients in order to obtain gain of integrability on the velocity.
\label{gammarem}
\end{remarka}
\begin{remarka}
In this theorem, if  we assume as in \cite{H3}  that $\frac{1}{\rho_{0}}\in L^{\infty}(\T^{N})$ then we obtain  a control on the gradient of the velocity $\n u$ in $L^{1}(BMO(\T^{N}))$. We know that for incompressible Navier-Stokes equations this hypothesis is enough to get uniqueness. In this sense we can consider our result as a theorem of \textit{strong-weak} solutions (indeed it is not enough to prove the uniqueness but \textit{almost}). It means that this result improves the results of \cite{H3} inasmuch as  we do not need any other assumption on $\rho_{0}$ than $\rho_{0}\in L^{\infty}(\T^{N})$.  We are then absolutely critical for the scaling of the equations on the density. For the initial velocity $u_{0}$ we need to be a little bit subcritical as $u_{0}\in H^{\N-1+\e}(\T^{N})$. The only thing is that in dimension $N=3$ we need extra assumption of the type $u_{0}\in L^{p}(\T^{N})$.
 \end{remarka}
In the following theorem we obtain strong solutions if we assume more regularity on the initial density and that $\rho_{0}$ is bounded away from the vacuum.  This will supply a  Lipschitz control  on the velocity $u$. By this way, we can show that the results of \cite{H3} are very critical as it seems necessary to add extra regularity to get a control of $\n u$ in $L^{1}_{T}(L^{\infty}(\T^{N}))$.
\begin{theorem}
\label{corollaire1}
Let $\e>0$. Under the hypothesis of theorem \ref{theo1} with in addition $u_{0}\in H^{\N-1+\e}(\T^{N})$ and $\inf\rho_{0}>0$, $\rho_{0}\in B^{\e}_{\infty,\infty}(\T^{N})$ if $P(\rho)=K\rho$ and $\rho_{0}\in B^{1}_{N,1}(\T^{N})\cap B^{\e}_{\infty,\infty}(\T^{N})$ for $N=3$ if $P$ is a general pressure law,
the solutions of theorem \ref{theo1} are  unique and verify locally in time (\ref{1.21}), (\ref{b1.22}), (\ref{b1.23}) and:
$$\n u\in L^{1}_{T_{0}}(L^{\infty}(\T^{N})).$$
Moreover if $\rho$ is in $L^{\infty}((0,T_{0})\times\T^{N})$ then the solutions can be extended beyond $T_{0}$.
\end{theorem}
\begin{remarka}
In fact by using exactly the same arguments than the proof of theorem \ref{corollaire1}, it would be easily possible to prove that this theorem can be adapted to the strong solution constructed in \cite{H3} (indeed the proof uses essentially energy inequalities) . It would suffices to choose initial data as in theorem \ref{corollaire1} except that we would not need to assume extra conditions on $\gamma$ when $N=3$ (see remark \ref{gammarem} for more explanations). It means that we could obtain theorem  \ref{corollaire1} with the condition $\gamma\geq 1$.
\label{gamma}
\end{remarka}
In the following theorem,  we want to improve the above blow-up criterion. More precisely we prove that it  suffices  to control the norm $L^{\infty}(L^{(N+1+\e)\gamma}(\T^{N}))$ of the density  with $\e>0$ when $P(\rho)=a\rho^{\gamma}$  to obtain global strong solutions. We refer to
\begin{theorem}
\label{theo3}
Let $\lambda=0$, $\gamma$ as in theorem \ref{theo1} and $g$ as in theorem \ref{theo1} and $g\in L^{\infty}(L^{\infty}(\T^{N}))$. Let $P(\rho)=a\rho^{\gamma}$ with $a>0$ and $\gamma\geq 1$. Assume that $(\rho_{0},u_{0})\in (L^{\gamma}(\T^{N})\cap L^{\infty}(\T^{N})\cap B^{1+\e}_{N,\infty}(\T^{N}))\times (L^{2}(\T^{N})\cap L^{\infty}(\T^{N})\cap H^{\N-1+\e}(\T^{N}))$ with $\e>0$ and that $\rho_{0}$ is bounded away from zero.\\
Let $(\rho,u)$ be a strong solution as in theorem \ref{corollaire1} of system (\ref{1}) on $[0,T)$ with the previous initial data which satisfies additionally the following conditions:
\begin{equation}
\begin{aligned}
&\rho\in L^{\gamma+1}_{T}(L^{(N+1+\e)\gamma}(\T^{N}))\;\;\mbox{and}\;\;\rho\in L^{\infty}_{T}( L^{9+\e}(\T^{N})\cap L^{3\gamma+\frac{3}{2}}(\T^{N}))\;\;\mbox{if}\;\;N=3,\\
&\rho\in L^{\gamma+1}_{T}(L^{(N+1+\e)\gamma}(\T^{N}))\;\;\mbox{and}\;\;\rho\in L^{\infty}_{T}(L^{2\gamma+1}(\T^{N}))\;\;\mbox{if}\;\;N=2 ,
\end{aligned}
\label{egcrucial}
\end{equation}
with $\e>0$. Then the solution  $(\rho,u)$ can be extended beyond $T$ and we have:
$$\n u\in L^{1}_{T}(L^{\infty}(\T^{N})).$$
\end{theorem}
\begin{remarka}
\begin{itemize}
\item This result has to be seen as a Prodi-Serrin theorem (see \cite{prodi} and \cite{serrin}) for compressible Navier-Stokes system. The main difference compared with incompressible Navier-Stokes system is that "the good variable" is the pressure and not the velocity. In some way, it is the integrability of the pressure which gives the regularity of the solutions. This result is the first one up to my knowledge which requires only condition of integrability on the density to get global strong solutions.
\item More precisely, as long than (\ref{egcrucial}) is verified, the regularity properties ( (\ref{1.21}),  (\ref{1.211}), (\ref{b1.22}), (\ref{b1.23}), (\ref{c1.23})  and $\n u\in L^{1}_{T}(L^{\infty}(\T^{N}))$ hold. It means in particular as $\n u\in L^{1}_{loc}(L^{\infty})$, that we can extend the solutions.
\item As explained in remark \ref{gamma}, we could extend the result for $\gamma\geq 1$.
\end{itemize}
\end{remarka}
\begin{remarka}
We believe that the assumption $\lambda=0$ may be weakened as the fact that $u_{0}$ belongs to $L^{\infty}(\T^{N})$. In return we would need stronger condition of integrability for the density.
\end{remarka}
\begin{remarka}
We believe that our method can be adapted to the euclidian space $\mathbb{R}^{N}$. This is the object of our future work.
\end{remarka}
Our paper is structured as follows. In section \ref{section2}, we give a few notation and briefly introduce the basic Fourier analysis
techniques needed to prove our result. In section \ref{section3} and \ref{section4}, we prove  a priori estimates on the density and the velocity,  more precisely we prove in particular that the \textit{effective velocity} is Lipschitz if we are able to control the norm $L^{\infty}$ of the density. In section \ref{section5} we prove the theorem \ref{theo1} by mollifying the initial data and by proving that the estimates of section \ref{section3} and \ref{section4} are uniform and independent of the mollifying process on the initial data. In section \ref{section6}, we prove additional regularity on the initial density which allows to get a Lipschitz control on the velocity when $\rho$ remains $L^{\infty}$. In the section \ref{section51}, we prove theorem \ref{theo1} and \ref{corollaire1} by constructing in particular approximate solutions. In section \ref{section7} we will give the proof of the blow-up theorem \ref{theo3}. We will conclude in section \ref{section8} by some comments and open problems.  An inescapable commutator estimate is postponed in the appendix in the section \ref{section9}.
\section{Littlewood-Paley theory and Besov spaces}
\label{section2}
Throughout the paper, $C$ stands for a constant whose exact meaning depends on the context. The notation $A\lesssim B$  means
that $A\leq CB$.
For all Banach space $X$, we denote by $C([0,T],X)$ the set of continuous functions on $[0,T]$ with values in $X$.
For $p\in[1,+\infty]$, the notation $L^{p}(0,T,X)$ or $L^{p}_{T}(X)$ stands for the set of measurable functions on $(0,T)$
with values in $X$ such that $t\rightarrow\|f(t)\|_{X}$ belongs to $L^{p}(0,T)$.
Littlewood-Paley decomposition  corresponds to a dyadic
decomposition  of the space in Fourier variables.
Let $\varphi\in C^{\infty}(\mathbb{R}^{N})$,
supported in the shell
${\cal{C}}=\{\xi\in\R^{N}/\frac{3}{4}\leq|\xi|\leq\frac{8}{3}\}$ and  $\chi\in C^{\infty}(\mathbb{R}^{N})$
supported in the ball $B(0,\frac{4}{3})$. $\varphi$ and  $\chi$ are valued in $[0,1]$.
We set ${\cal Q}^{N}=(0,2\pi)^{N}$ and $\widetilde{\mathbb{Z}}^{N}=(\mathbb{Z}/1)^{N}$ the dual lattice associated to $\mathbb{T}^{N}$. We decompose now $u\in{\cal S}^{'}(\T^{N})$ into Fourier series:
$$u(x)=\sum_{\beta\in \widetilde{\mathbb{Z}}^{N}}\hat{u}_{\beta}e^{i\beta\cdot x}\;\;\;\mbox{with}\;\;
\hat{u}_{\beta}=\frac{1}{|\T^{N}|}\int_{\T^{N}}e^{-i\beta\cdot y}u(y)dy.$$
Denoting;
$$h_{q}(x)=\sum_{\beta\in \widetilde{\mathbb{Z}}^{N}}\va(2^{-q}\beta)e^{i\beta\cdot x},$$
one can now define the periodic dyadic blocks as:
$$\D_{q}u(x)=\sum_{\beta\in \widetilde{\mathbb{Z}}^{N}}\va(2^{-q}\beta)\hat{u}_{\beta}e^{i\beta\cdot x}=\frac{1}{|\T^{N}|}\int_{\T^{N}}h_{q}(y)u(x-y)dy,\;\;\mbox{for all}\;\;q\in\mathbb{Z}$$
and the low frequency cutt-off:
$$S_{q}u(x)=\hat{u}_{0}+\sum_{p\leq q-1}\D_{p}u(x)=\sum_{\beta\in \widetilde{\mathbb{Z}}^{N}}\chi(2^{-q}\beta)\hat{u}_{\beta}e^{i\beta\cdot x}.$$
It is obvious that:
$$u=\hat{u}_{0}+\sum_{k\in\mathbb{Z}}\D_{k}u.$$
This decomposition is called non-homogeneous Littlewood-Paley
decomposition.
\\
Furthermore we have the following proposition where $\widetilde{{\cal C}}=B(0,\frac{2}{3})+{\cal{C}}$
\begin{proposition}
\label{d210}
\begin{equation}
|k-k^{'}|\geq 2\implies {\rm supp}\va(2^{-k}\cdot)\cap {\rm supp}\va(2^{-k^{'}}\cdot)=\emptyset,
\label{d28a}
\end{equation}
\begin{equation}
k\geq 1\implies {\rm supp}\chi\cap {\rm supp}\va(2^{-k}\cdot)=\emptyset,
\label{d29a}
\end{equation}
\begin{equation}
|k-k^{'}|\geq 5\implies 2^{k^{'}}\widetilde{{\cal C}}\cap 2^{k}{\cal{C}}=\emptyset.
\label{d210a}
\end{equation}
\end{proposition}
\subsection{Non-Homogeneous Besov spaces and first properties}
\begin{definition}
For
$s\in\R,\,\,p\in[1,+\infty],\,\,q\in[1,+\infty],\,\,\mbox{and}\,\,u\in{\cal{S}}^{'}(\T^{N})$
we set:
$$\|u\|_{B^{s}_{p,q}}=(\sum_{l\in\mathbb{Z}}(2^{ls}\|\D_{l}u\|_{L^{p}})^{q})^{\frac{1}{q}}.$$
The Besov space $B^{s}_{p,q}$ is the set of temperate distribution $u$ such that $\|u\|_{B^{s}_{p,q}}<+\infty$.
\end{definition}
\begin{proposition}
\label{derivation,interpolation}
The following properties holds:
\begin{enumerate}
\item there exists a constant universal $C$
such that:\\
$C^{-1}\|u\|_{B^{s}_{p,r}}\leq\|\n u\|_{B^{s-1}_{p,r}}\leq
C\|u\|_{B^{s}_{p,r}}.$
\item If
$p_{1}<p_{2}$ and $r_{1}\leq r_{2}$ then $B^{s}_{p_{1},r_{1}}\hookrightarrow
B^{s-N(1/p_{1}-1/p_{2})}_{p_{2},r_{2}}$.
\item $B^{s^{'}}_{p,r_{1}}\hookrightarrow B^{s}_{p,r}$ if $s^{'}> s$ or if $s=s^{'}$ and $r_{1}\leq r$.
\end{enumerate}
\label{interpolation}
\end{proposition}
Let now recall a few product laws in Besov spaces coming directly from the paradifferential calculus of J-M. Bony
(see \cite{5BJM}) and rewrite on a generalized form in \cite{AP} by H. Abidi and M. Paicu.
\begin{proposition}
\label{produit1}
We have the following laws of product:
\begin{itemize}
\item For all $s\in\R$, $(p,r)\in[1,+\infty]^{2}$ we have:
\begin{equation}
\|uv\|_{B^{s}_{p,r}}\leq
C(\|u\|_{L^{\infty}}\|v\|_{B^{s}_{p,r}}+\|v\|_{L^{\infty}}\|u\|_{B^{s}_{p,r}})\,.
\label{2.2}
\end{equation}
\item Let $(p,p_{1},p_{2},r,\lambda_{1},\lambda_{2})\in[1,+\infty]^{2}$ such that:
$$\frac{1}{p}\leq\frac{1}{p_{1}}+\frac{1}{p_{2}},\;p_{1}\leq\lambda_{2},\;p_{2}\leq\lambda_{1},\;\frac{1}{p}\leq\frac{1}{p_{1}}+\frac{1}{\lambda_{1}}\;\;\mbox{and}\;\;\frac{1}{p}\leq\frac{1}{p_{2}}+\frac{1}{\lambda_{2}}.$$
Assume that:
$$s_{1}+s_{2}+N\inf(0,1-\frac{1}{p_{1}}-\frac{1}{p_{2}})>0,\;s_{1}+\frac{N}{\lambda_{2}}<\frac{N}{p_{1}}\;\;\mbox{and}\;\;s_{2}+\frac{N}{\lambda_{1}}<\frac{N}{p_{2}},$$
then we have the following inequalities:
\begin{equation}
\|uv\|_{B^{s_{1}+s_{2}-N(\frac{1}{p_{1}}+\frac{1}{p_{2}}-\frac{1}{p})}_{p,r}}\lesssim\|u\|_{B^{s_{1}}_{p_{1},r}}
\|v\|_{B^{s_{2}}_{p_{2},\infty}}.
\label{2.3}
\end{equation}
When $s_{1}+\frac{N}{\lambda_{2}}=\frac{N}{p_{1}}$ (resp $s_{2}+\frac{N}{\lambda_{1}}=\frac{N}{p_{2}}$) we obtain a similar result as (\ref{2.3}) by replacing
$\|u\|_{B^{s_{1}}_{p_{1},r}}\|v\|_{B^{s_{2}}_{p_{2},\infty}}$ (resp $\|v\|_{B^{s_{2}}_{p_{2},\infty}}$) by
$\|u\|_{B^{s_{1}}_{p_{1},1}}\|v\|_{B^{s_{2}}_{p_{2},r}}$ (resp $\|v\|_{B^{s_{2}}_{p_{2},\infty}\cap L^{\infty}}$).\\
When $s_{1}+\frac{N}{\lambda_{2}}=\frac{N}{p_{1}}$ and $s_{2}+\frac{N}{\lambda_{1}}=\frac{N}{p_{2}}$ we replace by  $\|u\|_{B^{s_{1}}_{p_{1},1}}\|v\|_{B^{s_{2}}_{p_{2},1}}$ .
\end{itemize}
\end{proposition}
We now want to give some similar results for the critical case for the paraproduct laws, it means when  $s_{1}+s_{2}=0$.
\begin{proposition}
\label{produit2}
Assume that $s_{1}+s_{2}=0$, $s_{1}\in(\frac{N}{\lambda_{1}}-\frac{N}{p_{2}},\frac{N}{p_{1}}-\frac{N}{\lambda_{2}}]$ and
$\frac{1}{p_{1}}+\frac{1}{p_{2}}\leq 1$ then:
\begin{equation}
\|uv\|_{B^{-N(\frac{1}{p_{1}}+\frac{1}{p_{2}}-\frac{1}{p})}_{p,\infty}}\lesssim\|u\|_{B^{s_{1}}_{p_{1},1}}
\|v\|_{B^{s_{2}}_{p_{2},\infty}}.
\label{2.4}
\end{equation}
If $|s|<\NN$ for $p\geq2$ and $-\frac{N}{p^{'}}<s<\NN$ else, we have:
\begin{equation}
\|uv\|_{B^{s}_{p,r}}\leq C\|u\|_{B^{s}_{p,r}}\|v\|_{B^{\NN}_{p,\infty}\cap L^{\infty}}.
\label{2.5}
\end{equation}
\end{proposition}
\begin{remarka}
In the sequel $p$ will be either $p_{1}$ or $p_{2}$ and in this case $\frac{1}{\lambda}=\frac{1}{p_{1}}-\frac{1}{p_{2}}$
if $p_{1}\leq p_{2}$, resp $\frac{1}{\lambda}=\frac{1}{p_{2}}-\frac{1}{p_{1}}$
if $p_{2}\leq p_{1}$.
\end{remarka}
\begin{proposition}
\label{produit3}
Let $r\in [1,+\infty]$, $1\leq p\leq p_{1}\leq +\infty$ and $s$ such that:
\begin{itemize}
\item $s\in(-\frac{N}{p_{1}},\frac{N}{p_{1}})$ if $\frac{1}{p}+\frac{1}{p_{1}}\leq 1$,
\item $s\in(-\frac{N}{p_{1}}+N(\frac{1}{p}+\frac{1}{p_{1}}-1),\frac{N}{p_{1}})$ if $\frac{1}{p}+\frac{1}{p_{1}}> 1$,
\end{itemize}
then we have if $u\in B^{s}_{p,r}$ and $v\in B^{\frac{N}{p_{1}}}_{p_{1},\infty}\cap L^{\infty}$:
$$\|uv\|_{B^{s}_{p,r}}\leq C\|u\|_{B^{s}_{p,r}}\|v\|_{B^{\frac{N}{p_{1}}}_{p_{1},\infty}\cap L^{\infty}}.$$
\end{proposition}
The study of non stationary PDE's requires space of type $L^{\rho}(0,T,X)$ for appropriate Banach spaces $X$. In our case, we
expect $X$ to be a Besov space, so that it is natural to localize the equation through Littlewood-Paley decomposition. But, in doing so, we obtain
bounds in spaces which are not type $L^{\rho}(0,T,X)$ (except if $r=p$).
We are now going to
define the spaces of Chemin-Lerner (see \cite{5CL}) in which we will work, which are
a refinement of the spaces
$L_{T}^{\rho}(B^{s}_{p,r})$.
\begin{definition}
Let $\rho\in[1,+\infty]$, $T\in[1,+\infty]$ and $s_{1}\in\R$. We set:
$$\|u\|_{\widetilde{L}^{\rho}_{T}(B^{s_{1}}_{p,r})}=
\big(\sum_{l\in\mathbb{Z}}2^{lrs_{1}}\|\D_{l}u(t)\|_{L^{\rho}_{T}(L^{p})}^{r}\big)^{\frac{1}{r}}\,.$$
We then define the space $\widetilde{L}^{\rho}_{T}(B^{s_{1}}_{p,r})$ as the set of temperate distribution $u$ over
$(0,T)\times\T^{N}$ such that 
$\|u\|_{\widetilde{L}^{\rho}_{T}(B^{s_{1}}_{p,r})}<+\infty$.
\end{definition}
We set $\widetilde{C}_{T}(\widetilde{B}^{s_{1}}_{p,r})=\widetilde{L}^{\infty}_{T}(\widetilde{B}^{s_{1}}_{p,r})\cap
{\cal C}([0,T],B^{s_{1}}_{p,r})$.
Let us emphasize that, according to Minkowski inequality, we have:
$$\|u\|_{\widetilde{L}^{\rho}_{T}(B^{s_{1}}_{p,r})}\leq\|u\|_{L^{\rho}_{T}(B^{s_{1}}_{p,r})}\;\;\mbox{if}\;\;r\geq\rho
,\;\;\;\|u\|_{\widetilde{L}^{\rho}_{T}(B^{s_{1}}_{p,r})}\geq\|u\|_{L^{\rho}_{T}(B^{s_{1}}_{p,r})}\;\;\mbox{if}\;\;r\leq\rho
.$$
\begin{remarka}
It is easy to generalize propositions \ref{produit1}, \ref{produit2} and corollary \ref{produit3}
to $\widetilde{L}^{\rho}_{T}(B^{s_{1}}_{p,r})$ spaces. The indices $s_{1}$, $p$, $r$
behave just as in the stationary case whereas the time exponent $\rho$ behaves according to H\"older inequality.
\end{remarka}
In the sequel we will need of composition estimates in $\widetilde{L}^{\rho}_{T}(B^{s}_{p,r})$ spaces.
\begin{proposition}
\label{composition}
Let $s>0$, $(p,r)\in[1,+\infty]$ and $u\in \widetilde{L}^{\rho}_{T}(B^{s}_{p,r})\cap L^{\infty}_{T}(L^{\infty})$.
\begin{enumerate}
 \item Let $F\in W_{loc}^{[s]+2,\infty}(\T^{N})$ such that $F(0)=0$. Then $F(u)\in \widetilde{L}^{\rho}_{T}(B^{s}_{p,r})$. More precisely there exists a function $C$ depending only on $s$, $p$, $r$, $N$ and $F$ such that:
$$\|F(u)\|_{\widetilde{L}^{\rho}_{T}(B^{s}_{p,r})}\leq C(\|u\|_{L^{\infty}_{T}(L^{\infty})})\|u\|_{\widetilde{L}^{\rho}_{T}(B^{s}_{p,r})}.$$
\item Let $F\in W_{loc}^{[s]+3,\infty}(\T^{N})$ such that $F(0)=0$. Then $F(u)-F^{'}(0)u\in \widetilde{L}^{\rho}_{T}(B^{s}_{p,r})$. More precisely there exists a function $C$ depending only on $s$, $p$, $r$, $N$ and $F$ such that:
$$\|F(u)-F^{'}(0)u\|_{\widetilde{L}^{\rho}_{T}(B^{s}_{p,r})}\leq C(\|u\|_{L^{\infty}_{T}(L^{\infty})}\|u\|^{2}_{\widetilde{L}^{\rho}_{T}(B^{s}_{p,r})}.$$
\end{enumerate}
\end{proposition}
Here we recall a result of interpolation which explains the link
of the space $B^{s}_{p,1}$ with the space $B^{s}_{p,\infty}$, see
\cite{DFourier}.
\begin{proposition}
\label{interpolationlog}
There exists a constant $C$ such that for all $s\in\R$, $\e>0$ and
$1\leq p<+\infty$,
$$\|u\|_{\widetilde{L}_{T}^{\rho}(B^{s}_{p,1})}\leq C\frac{1+\e}{\e}\|u\|_{\widetilde{L}_{T}^{\rho}(B^{s}_{p,\infty})}
\biggl(1+\log\frac{\|u\|_{\widetilde{L}_{T}^{\rho}(B^{s+\e}_{p,\infty})}}
{\|u\|_{\widetilde{L}_{T}^{\rho}(B^{s}_{p,\infty})}}\biggl).$$ \label{5Yudov}
\end{proposition}
Now we give some result on the behavior of the Besov spaces via some pseudodifferential operator (see \cite{DFourier}).
\begin{definition}
Let $m\in\R$. A smooth function function $f:\T^{N}\rightarrow\R$ is said to be a ${\cal S}^{m}$ multiplier if for all muti-index $\alpha$, there exists a constant $C_{\alpha}$ such that:
$$\forall\xi\in\T^{N},\;\;|\p^{\alpha}f(\xi)|\leq C_{\alpha}(1+|\xi|)^{m-|\alpha|}.$$
\label{smoothf}
\end{definition}
\begin{proposition}
Let $m\in\R$ and $f$ be a ${\cal S}^{m}$ multiplier. Then for all $s\in\R$ and $1\leq p,r\leq+\infty$ the operator $f(D)$ is continuous from $B^{s}_{p,r}$ to $B^{s-m}_{p,r}$.
\label{singuliere}
\end{proposition}
We now focus on the mass equation associated to (\ref{1})
\begin{equation}
\begin{cases}
 \begin{aligned}
&\p_{t}\rho+u\cdot\n \rho+\frac{\rho}{2\mu+\lambda} \big(h(\rho)-\int_{\T^{N}}h(\rho)(t,x)dx\big) =-\frac{1}{2\mu+\lambda}\rho{\rm div}v_{1},\\
&q_{/ t=0}=q_{0}.
 \end{aligned}
\end{cases}
\label{21}
\end{equation}
where $h\in C^{\infty}$, $h(0)=0$ and $h^{'}\in W^{s,\infty}(\R,\R)$.  Here $v_{1}$ belongs in $\widetilde{L}^{1}(B^{\NN+\e}_{p_{1},1})$ with $\e>0$ and $p_{1}\in[1,+\infty]$.
\begin{proposition}
Let 
$1\leq p\leq p_{1}\leq+\infty$ and $1\leq r\leq+\infty$ with $p^{'}=\frac{p}{p-1}$. Let assume that:
\begin{equation}
-N\min(\frac{1}{p_{1}},\frac{1}{p^{'}})<\sigma.
 \label{3.21}
\end{equation}
Assume that $\rho_{0}\in B^{\sigma}_{p,r}$, $\n u\in \widetilde{L}^{1}_{T}(B^{\frac{N}{p_{1}}}_{p_{1},\infty}\cap L^{\infty})$, ${\rm div} v_{1}\in L^{1}(0,T;L^{\infty})$, ${\rm div}v_{1}\in \widetilde{L}^{1}_{T}(B^{\sigma}_{p,r})$  and that $\rho\in\widetilde{L}^{\infty}_{T}(B^{\sigma}_{p,r})\cap L^{\infty}$
satisfies (\ref{21}).
There exists a constant $C$ depending only on $N$ such that for all
$t\in[0,T]$ and $m\in\mathbb{Z}$, we have:
\begin{equation}
\|q\|_{\widetilde{L}^{\infty}_{t}(B^{\sigma}_{p,r})}\leq e^{CV(t)}\big(\|\rho_{0}\|_{B^{\sigma}_{p,r}}+\int^{t}_{0}C \|\rho(\tau,\cdot)\|_{L^{\infty}}\|{\rm div}v_{1}(\tau,\cdot)\|_{B^{\sigma}_{p,r}}d\tau\big),
\label{22}
\end{equation}
with:
$$
\begin{cases}
\begin{aligned}
V(t)&=\int^{t}_{0}\big(\|\n u\|_{B^{\frac{N}{p_{1}}}_{p_{1},\infty}\cap L^{\infty}}+\|{\rm div}v_{1}\|_{B^{\frac{N}{p_{1}}}_{p_{1},\infty}\cap L^{\infty}}+\|\rho\|^{\alpha+1}_{L^{\infty}}+1\big)d\tau\:\:\:\mbox{if}\;\;
\sigma<\frac{N}{p_{1}}+1,\\
&=\int^{t}_{0}\big(\|\n u\|_{B^{\frac{N}{p_{1}}}_{p_{1},1}}+\|{\rm div}v_{1}\|_{B^{\frac{N}{p_{1}}}_{p_{1},1}}+\|\rho\|^{\alpha+1}_{L^{\infty}}+1\big)d\tau\:\:\:\mbox{if}\;\;
\sigma=\frac{N}{p_{1}}+1\;\;\mbox{and}\;\;r=1.
\end{aligned}
\end{cases}
$$
with $\alpha$ the smallest integer such that $\alpha\geq s$.
\label{transport}
 \end{proposition}
{\bf Proof:}
 Applying $\D_{l}$ to (\ref{21}) yields:
$$\p_{t}\D_{l}\rho+u\cdot\n\D_{l}\rho=R_{l}-\D_{l}(\rho{\rm div}v_{1})-\D_{l}(\rho l(\rho))\;\;\;\mbox{with}\;\;R_{l}=[u\cdot\n,\D_{l}]\rho,$$
with $l(\rho)=h(\rho)-\int_{\T^{N}}h(\rho)(t,x)dx$. Multiplying by $\D_{l}\rho|\D_{l}\rho|^{p-2}$ and performing a time integration, we easily get:
$$
\begin{aligned}
&\|\D_{l}\rho(t)\|_{L^{p}}d\lesssim\|\D_{l}\rho_{0}\|_{L^{p}}+\int^{t}_{0}\big(\|R_{l}\|_{L^{p}}
+\|{\rm div}u\|_{L^{\infty}}\|\D_{l}\rho\|_{L^{p}}\\
&\hspace{7cm}+\|\D_{l}(\rho{\rm div}v_{1})\|_{L^{p}}+\|\D_{l}(\rho\, l(\rho))\|_{L^{p}}\big)d\tau.
\end{aligned}
$$
By paraproduct (see proposition \ref{produit3}), there exists a constant $C$ and a positive sequence $(c_{l})\in l^{r}$ such that:
$$\|\D_{l}(\rho{\rm div}v_{1})\|_{L^{p}}\leq C c_{l}2^{-l\sigma}\|\rho\|_{B^{\sigma}_{p,r}\cap L^{\infty}}\|{\rm div}v_{1}\|_{B^{\frac{N}{p_{1}}}_{p_{1},\infty}\cap L^{\infty}}.$$
Similarly by proposition \ref{composition} we have:
$$\|\D_{l}(\rho l(\rho))\|_{L^{p}}\leq C c_{l}2^{-q\sigma}\|\rho\|_{B^{\sigma}_{p,r}}(1+\|\rho\|^{\alpha+1}_{L^{\infty}}).$$
Next the term $\|R_{l}\|_{L^{p}}$ may be bounded according to the inequality of the lemma \ref{aKlemme3} in the appendix:
\begin{equation}
\begin{aligned}
\|\big(2^{l\sigma}\|R_{l}\|_{L^{p}}\big)_{l}\|_{l^{r}}&\leq C\|\n u \|_{B^{\frac{N}{p_{1}}}_{p_{1},\infty}\cap L^{\infty}}\|\rho\|_{B^{\sigma}_{p,r}}.
\end{aligned}
\label{abc13}
\end{equation}
We end up with multiplying the previous inequality by $2^{l\sigma}$ and summing up on $\mathbb{Z}$:
$$\|\rho(t)\|_{B^{\sigma}_{p,r}}\leq \|\rho_{0}\|_{B^{\sigma}_{p,r}}
+\int^{t}_{0}C V^{'}(\tau)\|\rho(\tau,\cdot)\|_{B^{\sigma}_{p,r}}d\tau+\int^{t}_{0}C \|\rho(\tau,\cdot)\|_{L^{\infty}}\|{\rm div}v_{1}(\tau,\cdot)\|_{B^{\sigma}_{p,r}}d\tau.$$
Gr\"onwall lemma yields inequality (\ref{22}).
\hfill{$\Box$}
\section{A priori bounds on the density}
\label{section3}
In this section we make a formal analysis of the partial differential equations of (\ref{1}) and begin by classical energy estimates. Multiplying the equation of conservation of momentum by $u$, we obtain
 \begin{equation}
\begin{aligned}
&\int_{\T^{N}}\big(\frac{1}{2}\rho
|u|^{2}(t,x)+\Pi(\rho)(t,x)\big)dx+\int_{0}^{t}\int_{\T^{N}}(\mu
D(u):D(u)(s,x)\\
&\hspace{2cm}+(\lambda+\mu)|{\rm div} u|^{2}(s,x))dsdx
\leq\int_{\T^{N}}\big(\frac{|m_{0}|^{2}}{2\rho}(x)+\Pi(\rho_{0})(x)\big)dx,
\end{aligned}
\label{17}
\end{equation}
where $\Pi$ is defined by
\begin{equation}
 \Pi(s)=s\biggl(\int^{s}_{0}\frac{P(z)}{z^{2}}dz\biggl),
\end{equation}
It follows classically that we have the following bounds
\begin{equation}
 \begin{cases}
  \begin{aligned}
&\rho\in L^{\infty}(0,\infty;L^{\gamma}),\\
&\sqrt{\rho}u\in L^{\infty}(0,\infty;L^{2}),\\
&\n u\in L^{2}(0,\infty;L^{2})^{N^{2}}.
  \end{aligned}
\end{cases}
\label{19}
\end{equation}
\subsection{$L^{\infty}$ bound on $\rho$}
In this section we are interesting in obtaining some $L^{\infty}$ bound on the density $\rho$. Indeed in the sequel it will be crucial for controlling the nonlinear terms as the pressure but alsp for some reasons related to the notion of multiplier space. Roughly to obtain some properties of regularity on the velocity $u$, we will need in a certain way that $\rho\D u$ keep the same regularity than $\D u$. It means that  the density has to belong to some multiplier space and in particular the density $\rho$ must be in $L^{\infty}$.\\
Let us mention that one of the main ingredient for this  is a partial differential equation close to a transport equation derived from (\ref{1}) involving $\log\rho$. In particular the notion of \textit{effective pressure} will appear. We recall in particular that it was one of the main ingredient used by  P-L Lions in \cite{13} to prove global existence of weak solutions of (\ref{1}).  It is also one of the key to the proof of B. Desjardins in \cite{CCK5}.\\
In this section for simplicity, we assume that $g=0$ (the general case is an easy extension). Applying now formally the operator $(\D)^{-1}{\rm div}$  on the equation of conservation of momentum, we obtain then:
\begin{equation}
(2\mu+\lambda){\rm div}u-P(\rho)+\int_{\T^{N}}P(\rho)dx=\p_{t}\D^{-1}{\rm
div}(\rho\,u)+R_{i}R_{j}(\rho\,u_{i}u_{j}),
\label{20}
\end{equation}
where $\D^{-1}$ denotes the inverse Laplacian with zero mean value on $\T^{N}$
and $R_{i}$ the usual Riesz transform. We now observe that the equation of mass write under the following form in the unknown $\ln\rho$:
\begin{equation}
\p_{t}\log\rho+u.\n\log\rho+{\rm div}u=0.
 \label{g21}
\end{equation}
We now define $F$ and $G$ by the following expression:
$$
\begin{aligned}
&F=(2\mu+\lambda)(\log\rho+\D^{-1}{\rm
div}(\rho\,u)),\\
&G=(2\mu+\lambda){\rm
div}u-P(\rho)+\int_{\T^{N}}P(\rho)(\dot,x)dx.
\end{aligned}
$$
Here $F$ is comparable to a density unknown (more precisely we will obtain some informations on the $L^{\infty}$ norme of $\ln\rho$ via $L^{\infty}$ estimates on $F$) and $G$ is the famous \textit{effective pressure}. Moreover we shall denote in the sequel respectively by $\mathbb{P}$ and $\mathbb{Q}$ the projection on the space of divergence-free and curl-free vector fields.
Combining (\ref{20}) and (\ref{g21}), we obtain the following transport equation on $F$ with some additional terms.
\begin{equation}
\p_{t}F+u\cdot \n F+P(\rho)-\int_{\T^{N}}P(\rho)dx
=[u_{j},R_{i}R_{j}](\rho u_{i}).
 \label{24}
\end{equation}
Next by using the characteristic method, we define the flow $\Psi$ of $u$ by
\begin{equation}
\begin{cases}
 \begin{aligned}
 &\p_{t}\Psi(t,s,x)=u(t,\Psi(t,s,x)),\\
&\Psi_{/ t=s}=x,
 \end{aligned}
\end{cases}
\label{25}
\end{equation}
and by the characteristic method,  as $F$ check a transport equation we get:
\begin{equation}
\begin{aligned}
&F(t,\Psi(t,0,x))=F_{0}(x)-\int^{t}_{0}P(\rho(s,\Psi(s,0,x)))dxds+\int^{t}_{0}\int_{\T^{N}}P(\rho(s,\cdot))dxdt\\
&\hspace{7cm}+\int^{t}_{0}([u_{j},R_{i}R_{j}](\rho u_{i})(s,\Psi(s,0,x))ds.
\end{aligned}
 \label{26}
\end{equation}
Using the fact that $\rho(\cdot)\geq 0$, we obtain
\begin{equation}
\begin{aligned}
&F(t,x)\leq F_{0}(\Psi(0,t,x))+\int^{t}_{0}\int_{\T^{N}}P(\rho(s,\cdot))dxdt\\
&\hspace{6cm}+\int^{t}_{0}([u_{j},R_{i}R_{j}](\rho u_{i})(s,\Psi(s,t,x))ds.
\end{aligned}
 \label{27}
\end{equation}
It follows that
\begin{equation}
\begin{aligned}
\log(\rho(t,x))\lesssim&\log(\|\rho_{0}\|_{L^{\infty}})+\|(\D)^{-1}{\rm div} m_{0}\|_{L^{\infty}}+\|(\D)^{-1}{\rm div}(\rho u)(t,\cdot)\|_{L^{\infty}}\\
&\hspace{1cm}+\int^{t}_{0}\int_{\T^{N}}P(\rho(s,\cdot))dxdt+\int^{t}_{0}\|[u_{j},R_{i}R_{j}](\rho u_{i})(s,\cdot)\|_{L^{\infty}}ds.
\end{aligned}
 \label{g28}
\end{equation}
By using the Besov embedding as $B^{1}_{N+\e,\infty}\h L^{\infty}$ (see proposition \ref{interpolation}), we obtain
\begin{equation}
\begin{aligned}
\log(\rho(t,x))\lesssim&\log(\|\rho_{0}\|_{L^{\infty}})+\|(\D)^{-1}{\rm div} m_{0}\|_{L^{\infty}}+\|(\D)^{-1}{\rm div}(\rho u)\|_{L^{\infty}}\\
&\hspace{0,5cm}+\int^{t}_{0}\int_{\T^{N}}P(\rho(s,\cdot))dxdt+\int^{t}_{0}\|[u_{j},R_{i}R_{j}](\rho u_{i})(s,\cdot)\|_{B^{1}_{N+\e,\infty}}ds,
\end{aligned}
\label{29}
\end{equation}
with $\e>0$. 
In view of R. Coifman, P-.L. Lions, Y. Meyer and S. Semmes \cite{1}, the following map
\begin{equation}
\begin{aligned}
&W^{1,r_{1}}(\T^{N})^{N}\times L^{r_{2}}(\T^{N})\rightarrow W^{1,r_{3}}(\T^{N})^{N}\\
&\hspace{3cm}(a,b)\rightarrow[a_{j},R_{i}R_{j}]b_{i}
\end{aligned}
\label{a33}
\end{equation}
is continuous for any $N\geq2$ as soon as $\frac{1}{r_{3}}=\frac{1}{r_{1}}+\frac{1}{r_{2}}$.
\\
In dimension $N=3$ we have by H\"older inequalities and as by embedding $W^{1, 3+\frac{\e}{2}}\h B^{1}_{3+\frac{\e}{2},\infty}$ (see proposition \ref{interpolation}):
$$\|[u_{j},R_{i}R_{j}](\rho u_{i})(s,\cdot)\|_{B^{1}_{3+\frac{\e}{4+\frac{\e}{3}},\infty}}\leq \|\rho^{\frac{1}{6+\e}} u\|_{L^{6+\e}}\|\n u(s)\|_{L^{6}}(1+\|\rho(s)\|_{L^{\infty}}).$$
We finally obtain for a.e $(t,x)\in(0,+\infty)\times\mathbb{T}^{N}$ from the previous inequality and (\ref{29}):
\begin{equation}
\begin{aligned}
\log(\rho(t,x))\lesssim&\log(\|\rho_{0}\|_{L^{\infty}})+\|(\D)^{-1}{\rm div} m_{0}\|_{L^{\infty}}+\|(\D)^{-1}{\rm div}(\rho u)\|_{L^{\infty}}\\
&+\int^{t}_{0}\int_{\T^{N}}P(\rho(s,\cdot))dxdt+\|\rho u\|_{L^{\infty}(L^{6+\e})}\int^{t}_{0}\|\n u(s)\|_{L^{6}}ds.
\end{aligned}
\label{29imp}
\end{equation}
with $\e>0$.
We have a similar result for $N=2$. The estimate (\ref{29imp})  will be useful in the sequel in order to estimate $\rho_{L^{\infty}}$.
\section{A priori estimates for the velocity}
\label{section4}
\subsection{Gain of integrability on the velocity $u$}
We want here derive estimate of integrability on the velocity $u$. This idea has been successively used in different papers,
we refer in particularly to \cite{5H2} and \cite{5MV,5MV2}. To do that, we multiply the momentum equation by $u|u|^{p_{1}-2}$ and we apply integration by parts:
$$
\begin{aligned}
&\frac{1}{p_{1}}\int_{\T^{N}}\rho|u|^{p_{1}}(t,x)dx+\mu\int^{t}_{0}\int_{\T^{N}}\big(|u|^{p_{1}-2}|\n u|^{2}(s,x)+\frac{p_{1}-2}{4}
|u|^{p_{1}-4}|\n|u|^{2}|^{2}(s,x)\big)dxds\\
&\hspace{1,8cm}+\lambda\int^{t}_{0}\int_{\T^{N}}\big(({\rm div}u)^{2}|u|^{p_{1}-2}(s,x)+\frac{p_{1}-2}{2}
{\rm div}u\sum_{i}u_{i}\p_{i}|u|^{2}|u|^{p_{1}-4}(s,x)\big)dsdx\\
&\hspace{2,5cm}-\int^{t}_{0}\int_{\T^{N}} P(\rho)\big({\rm div}u|u|^{p_{1}-2}+(p_{1}-2)
\sum_{i,k}u_{i}u_{k}\p_{i}u_{k}|u|^{p_{1}-4}\big)(s,x)dsdx\\
&\hspace{11cm}\leq\int_{\T^{N}}\rho_{0}|u_{0}|^{p_{1}}dx.
\end{aligned}
$$
We have then by Young's inequality:
$$
\begin{aligned}
&\frac{\lambda(p_{1}-2)}{2}\int^{t}_{0}\int_{\T^{N}}
{\rm div}u\sum_{i}u_{i}\p_{i}|u|^{2}|u|^{p_{1}-4}(s,x)dsdx\\
&=\lambda\frac{p_{1}-2}{2}\int^{t}_{0}\int_{\T^{N}}
{\rm div}u\,u\cdot\n (|u|^{2})|u|^{p_{1}-4}(s,x)dsdx\leq\\
&\lambda \frac{p_{1}-2}{2}(\frac{\eta}{2}\int^{t}_{0}\int_{\T^{N}}
|{\rm div}u|^{2}|u|^{p_{1}-2}(s,x)dsdx+\frac{2}{\eta}\int^{t}_{0}\int_{\T^{N}}|\n |u|^{2}|^{2}|u|^{p_{1}-4}(s,x)dsdx)\\
\end{aligned}
$$
If we choose $\eta$ such that:
$$\lambda \frac{\eta(p_{1}-2)\lambda}{4}=s\mu+\lambda,$$
for some $s\in(0,\frac{1}{N})$, by the fact that $({\rm div}u)^{2}\leq N |\n u|^{2}$ we therefore obtain:
$$
\begin{aligned}
&\frac{1}{p_{1}}\int_{\T^{N}}\rho|u|^{p_{1}}(t,x)dx+A_{s}\int^{t}_{0}\int_{\T^{N}}|u|^{p_{1}-2}|\n u|^{2}(s,x)dsdx\\
&\hspace{1cm}+B_{s}\int^{t}_{0}\int_{\T^{N}}
|u|^{p_{1}-4}|\n|u|^{2}|^{2}(s,x)dxds\leq
\int^{t}_{0}\int_{\T^{N}} P(\rho)\big({\rm div}u|u|^{p_{1}-2}\\
&\hspace{3cm}+\frac{p_{1}-2}{2}
u\cdot\n(|u|^{2})|u|^{p_{1}-4}\big)(s,x)dsdx+\int_{\T^{N}}\rho_{0}|u_{0}|^{p_{1}}dx,
\end{aligned}
$$
with $A_{s}=\mu(1-s\,N)$ and $B(s)=\frac{p_{1}-2}{4}\mu-\frac{(p_{1}-2)^{2}\lambda^{2}}{16(s\mu+\lambda)}$.
By Young inequality, we get still:
$$
\begin{aligned}
&\frac{1}{p_{1}}\int_{\T^{N}}\rho|u|^{p_{1}}(t,x)dx+A_{s}\int^{t}_{0}\int_{\T^{N}}|u|^{p_{1}-2}|\n u|^{2}(s,x)dsdx\\
&+B_{s}\int^{t}_{0}\int_{\T^{N}}
|u|^{p_{1}-4}|\n|u|^{2}|^{2}(s,x)dxds\leq
C_{\e}\int^{t}_{0}\int_{\T^{N}} P(\rho)^{2}|u|^{p_{1}-2}(s,x)dsdx\\
&\hspace{10cm}+\int_{\T^{N}}\rho_{0}|u_{0}|^{p_{1}}dx,
\end{aligned}
$$
with $C_{\e}$  big enough.
\subsubsection*{Case $N=3$}
We now want to use the fact that $\n |u|^{\frac{p_{1}}{2}}\in L^{2}_{t}(L^{2})$, which implies that when $N=3$, $u\in L^{p_{1}}_{t}(L^{3p_{1}})$. More precisely we have:
$$\|u\|^{\frac{p_{1}}{2}}_{L^{p_{1}}_{t}(L^{3p_{1}})}\lesssim\|\n(|u|^{\frac{p_{1}}{2}})+\widetilde{|u|^{\frac{p_{1}}{2}}}(\cdot)\|_{L^{2}(L^{2})},$$
where $\widetilde{|u|^{\frac{p_{1}}{2}}}$  is the integral of $|u|^{\frac{p_{1}}{2}}$ over $\T^{N}$ which means:
$$\widetilde{|u|^{\frac{p_{1}}{2}}}(t)=\int_{\T^{N}}|u|^{\frac{p_{1}}{2}}(t,x)dx.$$
We have then by H\"older's inequalities with $\frac{p_{1}-2}{3p_{1}}+\frac{2(p_{1}+1)}{3p_{1}}=1$ and $\frac{p_{1}-2}{p_{1}}+\frac{2}{p_{1}}=1$:
$$
\begin{aligned}
 |\int^{t}_{0}\int_{\T^{N}} P(\rho)^{2}|u|^{p_{1}-2}(s,x)dsdx|&\leq
\|P(\rho)^{2}\|_{L_{t}^{\frac{p_{1}}{2}}(L^{\frac{3p_{1}}{2(p_{1}+1)}})}
\||u|^{p_{1}-2}\|_{L_{t}^{\frac{p_{1}}{p_{1}-2}}(L^{\frac{3p_{1}}{p_{1}-2}})},\\
&\leq\|P(\rho)\|_{L_{t}^{p_{1}}(L^{\frac{3p_{1}}{p_{1}+1}})}^{2}\|u\|_{L^{p_{1}}_{t}(L^{3p_{1}})}^{p_{1}-2},\\
&\leq C\|P(\rho)\|_{L_{t}^{p_{1}}(L^{\frac{3p_{1}}{p_{1}+1}})}^{2}(
\|\n( |u|^{\frac{p_{1}}{2}})+\widetilde{|u|^{\frac{p_{1}}{2}}}(\cdot)\|_{L_{t}^{2}(L^{2})})^{2-\frac{4}{p_{1}}}.
\end{aligned}
$$
Remarking that $\int_{\T^{N}}\rho_{0}dx=M\ne 0$, we can write as $\gamma\geq\frac{6}{5}$:
$$
\begin{aligned}
&\widetilde{|u|^{\frac{p_{1}}{2}}}(s)\leq\frac{1}{M}(\|\rho(s,\cdot)\|_{L^{\gamma}}\|\n |u|^{\frac{p_{1}}{2}}(s,\cdot)\|_{L^{2}}+\int_{\T^{N}}\rho(s,x)|u(s,x)|^{\frac{p_{1}}{2}}dx),\\
&\big(\int^{t}_{0}\widetilde{|u|^{\frac{p_{1}}{2}}}^{2}(s)ds\big)^{\frac{1}{2}}\leq
\frac{1}{M}(\|\rho\|_{L^{\infty}_{t}(L^{\gamma})}\|\n |u|^{\frac{p_{1}}{2}}\|_{L^{2}_{t}(L^{2})}+(Mt)^{\frac{1}{2}}\|\rho^{\frac{1}{p_{1}}}u\|_{L^{\infty}_{t}(L^{p_{1}})}^{\frac{p_{1}}{2}}).\\
\end{aligned}
$$
We obtain then:
$$
\begin{aligned}
&\|P(\rho)\|_{L_{t}^{p_{1}}(L^{\frac{3p_{1}}{p_{1}+1}})}^{2}(
\|\n( |u|^{\frac{p_{1}}{2}})+\widetilde{|u|^{\frac{p_{1}}{2}}}(\cdot)\|_{L_{t}^{2}(L^{2})})^{2-\frac{4}{p_{1}}}\leq\\
&C_{1}\|P(\rho)\|_{L_{t}^{p_{1}}(L^{\frac{3p_{1}}{p_{1}+1}})}^{2}\big(
\|\n( |u|^{\frac{p_{1}}{2}})\|_{L_{t}^{2}(L^{2})}^{\frac{2p_{1}-4}{p_{1}}}+
(Mt)^{\frac{p_{1}-2}{p_{1}}}\|\rho^{\frac{1}{p_{1}}}u\|_{L^{\infty}_{t}(L^{p_{1}})}^{p_{1}-2})
\end{aligned}
$$
By a standard application of Young inequality ($\frac{2p_{1}-4}{2p_{1}}+\frac{4}{2p_{1}}=1$), we obtain that:
\begin{equation}
\begin{aligned}
&\frac{1}{p_{1}}\int_{\T^{N}}\rho|u|^{p_{1}}(t,x)dx+A_{s}\int^{t}_{0}\int_{\T^{N}}|u|^{p_{1}-2}|\n u|^{2}(s,x)dsdx\\
&+B_{s}\int^{t}_{0}\int_{\T^{N}}
|u|^{p_{1}-4}|\n|u|^{2}|^{2}(s,x)dxds\leq C_{\e,t}^{2}\|P(\rho)\|_{L_{t}^{p_{1}}(L^{\frac{3p_{1}}{p_{1}+1}})}^{2p_{1}}+\frac{1}{p_{1}}\int_{\T^{N}}\rho_{0}|u_{0}|^{p_{1}}dx,\\
\end{aligned}
\label{againvitesse}
\end{equation}
where $C_{\e,t}$ is big enough and depend on the time $t$.
\subsubsection*{Case $N=2$}
By proceeding similarly, we have for all $q>1$ in particular $q$ large:
$$\|u\|_{L^{p_{1}}(L^{\frac{p_{1}q}{2}})}^{\frac{p_{1}}{2}}\leq C(\|\n |u|^{\frac{p_{1}}{2}}\|_{L^{2}(L^{2})}+\bar{u}_{\frac{p_{1}}{2}}).$$
We have then by H\"older's inequalities with $\frac{2(p_{1}-2)}{qp_{1}}+\frac{(q-2)p_{1}+4)}{qp_{1}}=1$ and $\frac{p_{1}-2}{p_{1}}+\frac{2}{p_{1}}=1$:
$$
\begin{aligned}
& |\int^{t}_{0}\int_{\T^{N}} P(\rho)^{2}|u|^{p_{1}-2}(s,x)dsdx|\leq
\|P(\rho)^{2}\|_{L_{t}^{\frac{p_{1}}{2}}(L^{\frac{qp_{1}}{(q-2)p_{1}+4}})}
\||u|^{p_{1}-2}\|_{L_{t}^{\frac{p_{1}}{p_{1}-2}}(L^{\frac{qp_{1}}{2(p_{1}-2)}})},\\
&\hspace{5,4cm}\leq\|P(\rho)\|_{L_{t}^{p_{1}}(L^{\frac{2qp_{1}}{(q-2)p_{1}+4}})}^{2}\|u\|_{L^{p_{1}}_{t}(L^{qp_{1}})}^{p_{1}-2},\\
&\hspace{3,5cm}\leq C\|P(\rho)\|_{L_{t}^{p_{1}}(L^{\frac{2qp_{1}}{(q-2)p_{1}+4}})}^{2}(
\|\n( |u|^{\frac{p_{1}}{2}})+\widetilde{|u|^{\frac{p_{1}}{2}}}(\cdot)\|_{L_{t}^{2}(L^{2})})^{2-\frac{4}{p_{1}}},
\end{aligned}
$$
and so:
\begin{equation}
\begin{aligned}
&\frac{1}{p_{1}}\int_{\T^{N}}\rho|u|^{p_{1}}(t,x)dx+A_{s}\int^{t}_{0}\int_{\T^{N}}|u|^{p_{1}-2}|\n u|^{2}(s,x)dsdx\\
&+B_{s}\int^{t}_{0}\int_{\T^{N}}
|u|^{p_{1}-4}|\n|u|^{2}|^{2}(s,x)dxds\leq C_{\e,t}^{2}\|P(\rho)\|_{L_{t}^{p_{1}}(L^{\frac{2qp_{1}}{(q-2)p_{1}+4}})}^{2}+\frac{1}{p_{1}}\int_{\T^{N}}\rho_{0}|u_{0}|^{p_{1}}dx.\\
\end{aligned}
\end{equation}
\subsection{Gain of derivatives on the velocity $u$}
\label{gainvitesse}
In this section we deal with the case $N=3$. The case $N=2$ follows the same lines. In the sequel we will follow the procedure
developed in \cite{CCK5} and \cite{5H2} to get some energy inequalities. The main idea compared with the results in \cite{CCK2, CCK3,5CK2} is to obtain energy inequalities which depend only on the control on $\rho\in L^{\infty}$. It implies that we have to be careful to not introduce some derivatives on the density in the goal to ``kill'' the coupling between velocity and pressure. To do this the notions of \textit{effective pressure} and \textit{effective velocity} will play a crucial role. One of the main difference with the work of Desjardins in \cite{CCK5} will correspond to take in consideration the gain of integrability on the velocity that we can obtain when the pressure is enough integrable. Furthermore the fact to obtain energy inequalities by multiplying by $f(t)=\min(1,t)$ will allow us to get weak-strong solutions with very critical initial data for the scaling of the equations. Finally compared with \cite{5H2}, as the initial data are large, it will add many technical difficulties in particular in the obtention of bootstrap estimates.\\
Multiplying first the equation of conservation of momentum by $f(t)\p_{t}u$ and integrating over $(0,T)\times\T^{N}$, we deduce that:
\begin{equation}
 \begin{aligned}
  &\int^{t}_{0}\int_{\T^{N}}f(s)\rho|\p_{t}u|^{2}dxds+\frac{1}{2}\int_{\T^{N}}f(t)\big(\mu|\n u(t,x)|^{2}+(\lambda+\mu)|{\rm div}u|^{2}(t,x)\big)dx\\
&+\int^{t}_{0}\int_{\T^{N}}\n P(\rho)\cdot f(s)\p_{t}u dxds\leq \int^{t}_{0}\int_{\T^{N}}\frac{f^{'}(s)}{2}\big(\mu|\n u(s,x)|^{2}+(\lambda+\mu)|{\rm div}u|^{2}(s,x)\big)dxds\\
&\hspace{2cm}+\int^{t}_{0}\|\sqrt{f(s)\rho}\,\p_{t}u\|_{L^{2}(\T^{N})}
(\|\sqrt{f(s)\rho}\,(u\cdot\n)u\|_{L^{2}(\T^{N})}+\|\sqrt{\rho}g\|_{L^{2}(\T^{N})})ds.
 \end{aligned}
\label{a74}
\end{equation}
Next we use the equation of mass conservation to write:
$$
\begin{aligned}
&\int_{\T^{N}}\n P(\rho)\cdot f(t)\p_{t}u dx=-\int_{\T^{N}}P(\rho)f(t)\p_{t}{\rm div}u \,dx,\\
&=-\p_{t}\int_{\T^{N}}f(t)P(\rho){\rm div}u \,dx+\int_{\T^{N}}\p_{t}(f(t)P(\rho)){\rm div}u\, dx\\
&=-\p_{t}\int_{\T^{N}}f(t)P(\rho){\rm div}u\, dx-\int_{\T^{N}}f(t)\big[{\rm div}(P(\rho)u){\rm div}u+(\rho P^{'}(\rho)-P(\rho)){\rm div}u\big] dx\\
&\hspace{9cm}+\int_{\T^{N}}f^{'}(t)P(\rho){\rm div}u\, dx,\\
&=-\p_{t}\int_{\T^{N}}f(t)P(\rho){\rm div}u \,dx+\frac{1}{2\mu+\lambda}\int_{\T^{N}}f(t)P(\rho)u\cdot\n (G+P(\rho)) dx\\
&-\frac{1}{(2\mu+\lambda)^{2}}\int_{\T^{N}}f(t)(\rho P^{'}(\rho)-P(\rho))\big(G^{2}-P(\rho)^{2}+2(\lambda+2\mu)P(\rho){\rm div}u\big)dx \\
&\hspace{9,5cm}+\int_{\T^{N}}f^{'}(t)P(\rho){\rm div}u \,dx,
\end{aligned}
$$
We now set:
$$\Pi_{f}(s)=s(\int^{s}_{0}\frac{f(z)}{z^{2}}dz),$$
with $f$ a $C^{\infty}$ function. We have then by using the mass equation the following equality:
$$\p_{t}\Pi_{f}(\rho)+{\rm div}(\Pi_{f}(\rho)u)+f(\rho){\rm div}u=0.$$
By this fact we obtain that:
$$P(\rho)\big(u\cdot\n P(\rho)-2(\rho P^{'}(\rho)-P(\rho)){\rm div}u\big)=P(\rho)(-\p_{t}(P(\rho))-3P^{'}(\rho)\rho{\rm div}u+2P(\rho){\rm div}u).$$
We have then:
$$
\begin{aligned}
&\int_{\T^{N}} P(\rho)\big(u\cdot\n P(\rho)-2(\rho P^{'}(\rho)-P(\rho)){\rm div}u\big)dx=\\
&-\frac{1}{2}\int_{\T^{N}}\p_{t}(P(\rho)^{2})dx+3\int_{\T^{N}}\p_{t}\Pi_{3P^{'}(\rho)\rho}dx-2\int_{\T^{N}}\p_{t}\Pi_{P(\rho)}dx.
\end{aligned}
$$
Next by integration by parts, we have:
$$\Pi_{3P^{'}(s)s}=\frac{1}{2}P(s)^{2}+\frac{s}{2}\int^{s}_{0}\frac{P(z)^{2}}{z^{2}}dz.$$
So:
$$
\begin{aligned}
&\int_{\T^{N}} P(\rho)\big(u\cdot\n P(\rho)-2(\rho P^{'}(\rho)-P(\rho)){\rm div}u\big)dx=
\int_{\T^{N}}\p_{t}(P(\rho)^{2}-\Pi_{P(\rho)})dx.
\end{aligned}
$$
We set $k(s)=P(s)^{2}-\frac{1}{2}\Pi_{P(s)}$. Let us observe that in the $P_{a,\gamma}$ case, we have  $k(s)=
a^{2}s^{2\gamma}((2\gamma-\frac{3}{2})/(2\gamma-1))$.\\
We obtain finally:
$$
\begin{aligned}
&\int_{\T^{N}}\n P(\rho)\cdot f(t)\p_{t}u dx=-\p_{t}\int_{\T^{N}}f(t)P(\rho){\rm div}u dx+\frac{1}{\lambda+2\mu}\p_{t}\int_{\T^{N}}f(t)k(\rho)dx\\
&+\frac{1}{2\mu+\lambda}\int_{\T^{N}}f(t)P(\rho)u\cdot\n G dx+\frac{1}{(2\mu+\lambda)^{2}}\int_{\T^{N}}f(t)P^{2}(\rho)(\rho P^{'}(\rho)-P(\rho)) dx\\
&\hspace{1,8cm}-\frac{1}{(2\mu+\lambda)^{2}}\int_{\T^{N}}f(t)G^{2}(\rho P^{'}(\rho)-P(\rho)) dx-\frac{1}{\lambda+2\mu}\int_{\T^{N}}f^{'}(t)k(\rho)dx\\
&\hspace{9cm}+\int_{\T^{N}}f^{'}(t)P(\rho){\rm div}u dx.
\end{aligned}
$$
Inserting the above inequality in (\ref{a74}) and by Young's inequality, we obtain:
\begin{equation}
 \begin{aligned}
  \int^{t}_{0}\int_{\T^{N}}f(s)\rho|\p_{t}u|^{2}dxds+\frac{1}{2}\int_{\T^{N}}f(t)\big(\mu|\n u(t,x)|^{2}+(\lambda+\mu)({\rm div}u(t,x))^{2}\big)dx&\\
+\frac{1}{(2\mu+\lambda)^{2}}\int^{t}_{0}\int_{\T^{N}}f(s)P(\rho)^{2}(\rho P^{'}(\rho)-P(\rho))dx ds&\\
+\frac{1}{\lambda+2\mu}\int_{\T^{N}}f(t)k(\rho(t,x))dx
\leq C+\int_{\T^{N}}f(t)P(\rho(t,x)){\rm div}u(t,x)dx&\\
+\frac{1}{\lambda+2\mu}\int^{t}_{0}\int_{\T^{N}}f^{'}(s)k(\rho)dxds
-\int^{t}_{0}\int_{\T^{N}}f^{'}(t)P(\rho){\rm div}u dx&\\
+C\int^{t}_{0}\int_{\T^{N}}f(s)(|\rho P^{'}(\rho)-P(\rho)|G^{2}+|P(\rho)u||\n G|+|\sqrt{\rho}u\cdot\n u|^{2}&\\
+|\sqrt{\rho}g|^{2})dxds.&
 \end{aligned}
\label{79}
\end{equation}
In the sequel we set:
$$
\begin{aligned}
&A(t)=\int^{t}_{0}\int_{\T^{N}}f(s)\rho|\p_{t}u|^{2}dxds+\frac{1}{2}\int_{\T^{N}}f(t)\big(\mu|\n u(t,x)|^{2}+(\lambda+\mu)({\rm div}u(t,x))^{2}\big)dx\\
&+\frac{1}{(2\mu+\lambda)^{2}}\int^{t}_{0}\int_{\T^{N}}f(s)P(\rho)^{2}(\rho P^{'}(\rho)-P(\rho))dx ds+\frac{1}{\lambda+2\mu}\int_{\T^{N}}f(t)k(\rho(t,x))dx.
\end{aligned}
$$
We obtain finally by using H\"older inequalities:
\begin{equation}
 \begin{aligned}
&A(t)\leq C+C_{t}\|\rho\|_{L^{\infty}}+C\int^{t}_{0}\int_{\T^{N}}f(s)(|\rho P^{'}(\rho)-P(\rho)|G^{2}+|P(\rho)u||\n G|\\
&\hspace{8cm}+|\sqrt{\rho}u\cdot\n u|^{2}+|\sqrt{\rho}g|^{2})dxds.\\
&\leq C+C\int^{t}_{0}(\|\rho P^{'}(\rho)-P(\rho)\|_{L^{\infty}}\|\sqrt{f(s)}G\|_{L^{2}}^{2}+f(s)\|g(\rho)(s,\cdot)\|_{L^{\infty}}
\|\sqrt{\rho}u\|_{L^{2}}\\
&\hspace{2cm}\times\|\n G\|_{L^{2}}+\|\sqrt{\rho}u\|_{L^{4}}^{2}\|\sqrt{f(s)}\n u\|_{L^{4}}^{2}+\|\rho\|_{L^{\infty}}\|\sqrt{f(s)}g\|_{L^{2}}^{2})dxds,\\
&\leq C+C\int^{t}_{0}(\|h(\rho(s,\cdot)\|_{L^{\infty}}\|\sqrt{f(s)}\n u\|_{L^{2}}^{2}+f(s)\|i(\rho(s,\cdot)\|_{L^{\infty}}+f(s)\|\n G\|_{L^{2}}
\\
&\hspace{1cm}\times \|k(\rho(s,\cdot))\|_{L^{\infty}}+\|\sqrt{\rho}u\|_{L^{4}}^{2}\|\sqrt{f(s)}\n u\|_{L^{4}}^{2}+\|\rho\|_{L^{\infty}}\|\sqrt{f(s)}g\|_{L^{2}}^{2})dxds,
 \end{aligned}
\label{a81}
\end{equation}
where $k(\rho)=\frac{P(\rho)}{\sqrt{\rho}}$, $h(s)=|s P^{'}(s)-P(s)|$ and $i(s)=h(s)P(s)^{2}$.
\subsubsection*{Estimates on $\mathbb{P}u$ and $G$}
We now want  to obtain bounds on $\mathbb{P}u$ and $G$, assuming that $\rho$ is a priori bounded in $L^{\infty}(\T^{N})$.
Indeed we want to show that the control of $A(t)$ in (\ref{a81}) depends only on a control on $\|\rho\|_{L^{\infty}}$.\\
To do this, we use once more the equation of conservation of momentum in order to take in account the ellipticity of the momentum equation. To cancel out the coupling between the pressure and the velocity, we will consider the unknown $\mathbb{P}u$ and $G$. By applying the operator $\mathbb{P}u$ and by studying the effective pressure, we obtain the following equations:
\begin{equation}
\mu\D \mathbb{P}u=\mathbb{P}(\rho\dot{u})-\mathbb{P}(\rho g),
 \label{83}
\end{equation}
\begin{equation}
\n G=  \mathbb{Q}(\rho\dot{u})- \mathbb{Q}(\rho g),
 \label{84}
\end{equation}
where $\dot{u}=\p_{t}u+u\cdot\n u$. We recall here that $ \mathbb{P}$ is the projector on free-divergence vector field and
$ \mathbb{Q}$ is the projector on gradient vector field.
Therefore we have:
\begin{equation}
\begin{aligned}
&\|\n G\|_{L^{2}}+\|\D\mathbb{P}u\|_{L^{2}}\leq C\|\rho(s,\cdot)\|^{\frac{1}{2}}_{L^{\infty}}\big(\|\sqrt{\rho}\p_{s}u(s\cdot)\|_{L^{2}}
+\|\sqrt{\rho}u\cdot\n u(s,\cdot)\|_{L^{2}}\\
&\hspace{5cm}+\|\rho(s,\cdot)\|^{\frac{1}{2}}_{L^{\infty}}\|g(s,\cdot)\|_{L^{2}}\big).
\end{aligned}
 \label{86}
\end{equation}
\subsubsection*{The case $N=3$}
For simplicity, we will treat only the case of the dimension $3$. We recall that for all $1<p<+\infty$ by Calderon-Zygmund theory we have:
$$\|\n u\|_{L^{p}}\leq C(\|\n \mathbb{P}u\|_{L^{p}}+\|R G\|_{L^{p}}+\|R(P(\rho))\|_{L^{p}}),$$
where $R$ is a pseudo differential operator of order $0$ such that for all $f\in H^{1}(\T^{N})$ $\int_{\T^{N}}Rfdx=0$.
We want now to recall the Gagliardo-Nirenberg's theorem that we will use frequently:
$$\forall f\in H^{1}(\T^{N})\;\;\mbox{such that}\;\;\int_{\T^{N}}f dx=0,\;\;\|f\|^{2}_{L^{4}(\T^{N})}\leq C\|f\|^{\frac{1}{2}}_{L^{2}(\T^{N})}
\|\n f\|_{L^{2}(\T^{N})}^{\frac{3}{2}}.$$
We deduce that from Gagliardo-Nirenberg's inequality, Young's inequalities and (\ref{86}):
\begin{equation}
\begin{aligned}
&\|\sqrt{\rho}u\|^{2}_{L^{4}}\|\sqrt{f(s)}\n u\|^{2}_{L^{4}}\leq Cf(s)\|\sqrt{\rho}u\|^{2}_{L^{4}}(\|R(P(\rho))\|^{2}_{L^{4}}
+\|\n \mathbb{P}u\|_{L^{4}}+\|R G\|_{L^{4}}\big)\\
&\leq C\|\sqrt{\rho}u\|^{2}_{L^{4}}\big(f(s)\|P(\rho)\|^{2}_{L^{4}}+\big(\sqrt{f(s)}(\|\n u\|_{L^{2}}+\|P(\rho)\|_{L^{2}}))^{\frac{1}{2}}\\
&\hspace{6cm}\times(\sqrt{f(s)}(\|\D \mathbb{P}u\|_{L^{2}}+\|\n G\|_{L^{2}})\big)^{\frac{3}{2}}\big)\\
&\leq C\big(f(s)\|\sqrt{\rho}u\|^{2}_{L^{4}}\|P(\rho)\|^{2}_{L^{4}}+\|\rho(s,\cdot)\|^{\frac{3}{4}}_{L^{\infty}}\|\sqrt{\rho}u\|^{2}_{L^{4}}\big(\sqrt{f(s)}(\|\n u\|_{L^{2}}+\|P(\rho)\|_{L^{2}}))^{\frac{1}{2}}\\
&\times(\sqrt{f(s)}(\|\sqrt{\rho}\p_{s}u(s\cdot)\|_{L^{2}}
+\|\sqrt{\rho}u\cdot\n u(s,\cdot)\|_{L^{2}}+\|\rho(s,\cdot)\|^{\frac{1}{2}}_{L^{\infty}}\|f(s,\cdot)\|_{L^{2}}\big)^{\frac{3}{2}}\big)\\[3mm]
&\leq C\big(f(s)\|\sqrt{\rho}u\|^{2}_{L^{4}}\|P(\rho)\|^{2}_{L^{4}}+\frac{1}{\e}f(s)(\|\n u\|^{2}_{L^{2}}+\|P(\rho)\|^{2}_{L^{2}})\|\rho(s,\cdot)\|^{3}_{L^{\infty}}\|\sqrt{\rho}u\|^{8}_{L^{4}}\\[3mm]
&+\e f(s)(\|\sqrt{\rho}\p_{s}u(s\cdot)\|^{2}_{L^{2}}+\|\sqrt{\rho}u\cdot\n u(s,\cdot)\|^{2}_{L^{2}}+\|\rho(s,\cdot)\|_{L^{\infty}}\|g(s,\cdot)\|^{2}_{L^{2}})\big).\\
\end{aligned}
 \label{99}
\end{equation}
Hence we obtain by Young inequality from (\ref{86}):
\begin{equation}
\begin{aligned}
&f(s)\|k(\rho(s,\cdot)\|_{L^{\infty}}\|\n G(s,\cdot)\|_{L^{2}}\leq \frac{C}{\e}\|k(\rho(s,\cdot)\|^{2}_{L^{\infty}}\|\rho(s,\cdot)\|_{L^{\infty}}f(s)\\
&+\e(\|\sqrt{f(s)\rho}\p_{s}u(s\cdot)\|^{2}_{L^{2}}
+\|\sqrt{\rho f(s)}u\cdot\n u(s,\cdot)\|^{2}_{L^{2}}+f(s)\|\rho(s,\cdot)\|_{L^{\infty}}\|g(s,\cdot)\|^{2}_{L^{2}}).
\end{aligned}
\label{100}
\end{equation}
By adding (\ref{99}) and (\ref{100}), we obtain:
\begin{equation}
\begin{aligned}
&\|\sqrt{\rho}u\|^{2}_{L^{4}} \|\sqrt{f(s)}\n u\|^{2}_{L^{4}} +f(s)\|k(\rho(s,\cdot)\|_{L^{\infty}}\|\n G(s,\cdot)\|_{L^{2}}\leq\\
&C\big(f(s)\|\sqrt{\rho}u\|^{2}_{L^{4}} \|P(\rho)\|^{2}_{L^{4}}+\frac{1}{\e}f(s)(\|\n u\|^{2}_{L^{2}}+\|P(\rho)\|^{2}_{L^{2}})\|\rho(s,\cdot)\|^{3}_{L^{\infty}}\|\sqrt{\rho}u\|^{8}_{L^{4}}\\
&+\e f(s)(\|\sqrt{\rho}\p_{s}u(s\cdot)\|^{2}_{L^{2}}+\|\rho(s,\cdot)\|_{L^{\infty}}\|g(s,\cdot)\|^{2}_{L^{2}})\big)+\frac{C}{\e}f(s)\|k(\rho(s,\cdot)\|^{2}_{L^{\infty}}\\
&\hspace{10cm}\times\|\rho(s,\cdot)\|_{L^{\infty}}.
\end{aligned}
\label{aa81}
\end{equation}
Here at the difference of B. Desjardins in \cite{CCK5} we will use the gain of integrability on the velocity obtained in inequality (\ref{againvitesse}) to control the term $\|\sqrt{\rho}u\|_{L^{\infty}_{T}(L^{4})}$. In  \cite{CCK5} B. Desjardins estimate this quantity via the control on $\|\n u\|_{L^{\infty}_{T}(H^{1})}$. In particular, it explains why in dimension $N=3$, he can not estimate the regularizing effects on $u$ only by a control of $\rho$ in norm $L^{\infty}$. By (\ref{againvitesse}), we have then:
\begin{equation}
\|\sqrt{\rho}u\|_{L^{\infty}_{T}(L^{4})}\leq C\|\rho\|^{\frac{1}{4}}_{L^{\infty}_{T}(L^{\infty})}(1+\|P(\rho)\|_{L^{4}(L^{3})}^{2}).
\label{againvitesse1}
\end{equation}
Therefore we have by using inequality (\ref{againvitesse1}), (\ref{a81}) and (\ref{aa81}):
\begin{equation}
 \begin{aligned}
&A(t)\leq C+C\int^{t}_{0}(\|h(\rho(s,\cdot))\|_{L^{\infty}}+\frac{1}{\e}\|\rho(s,\cdot)\|_{L^{\infty}}^{5}
\|\rho^{\frac{1}{4}}u\|_{L^{4}})f(s)\|\n u\|_{L^{2}}^{2}\\
&+\frac{1}{\e}f(s)\|P(\rho(s,\cdot))\|_{L^{2}}^{2}\|\rho(s,\cdot)\|_{L^{\infty}}^{5}
\|\rho^{\frac{1}{4}}u\|_{L^{4}}+f(s)\|P(\rho(s,\cdot))\|_{L^{\infty}}^{2}\|\rho(s,\cdot)\|_{L^{\infty}}^{\frac{1}{2}}\\\
&\hspace{2cm}\times|\rho^{\frac{1}{4}}u\|_{L^{4}}^{2}+\|\rho(s,\cdot)\|_{L^{\infty}}\|g(s,\cdot)_{L^{2}}^{2}+\phi(\|\rho(s,\cdot)\|_{L^{\infty}})ds.
 \end{aligned}
\end{equation}
We have then by using (\ref{againvitesse1})
\begin{equation}
 \begin{aligned}
&A(t)\leq C+C(1+(\int^{t}_{0}\|P(\rho(s,\cdot))\|_{L^{\infty}}^{4}ds)^{4})\int^{t}_{0}\phi(\|\rho(s,\cdot)\|_{L^{\infty}})f(s)\big((1+\|\n u\|_{L^{2}}^{2}\\
&\hspace{5cm}+\|P(\rho(s,\cdot))\|^{2}_{L^{2})}+\|g(s,\cdot)\|_{L^{2}}^{2+\alpha}\big)ds,\\
&\leq C+C(1+\int^{t}_{0}\|P(\rho(s,\cdot))\|_{L^{\infty}}^{16}ds)\int^{t}_{0}\phi(\|\rho(s,\cdot)\|_{L^{\infty}})(A(s)+f(s))+\|g(s,\cdot)\|_{L^{2}}^{2+\alpha}\big)ds,
 \end{aligned}
\end{equation}
where $C$ depends on the time $t$ and $\alpha>0$. Here $\phi$ is in $C^{0}(\R_{+},\R_{+}^{*})\cap C^{1}(0,\infty)$ such that $\phi(s)\leq M+Cs^{\beta}$ for some positive $M,C>0$ and $\beta>1$. We define here by ${\cal F}$ the space endowed with this type of function. In the sequel as we will use a function $\phi_{\beta}$ with $\beta\in\mathbb{N}$ it will means that $\phi_{\beta}\in{\cal F}$.
Gronwall's lemma provides the following bound:
\begin{equation}
\begin{aligned}
&A(t)\leq C(1+\int^{t}_{0}P^{16}(\|\rho(s,\cdot)\|_{L^{\infty}})ds\int^{t}_{0}\phi(\|\rho(s,\cdot)\|_{L^{\infty}}ds
+\int^{t}_{0}\phi_{1}(\|\rho(s,\cdot)\|_{L^{\infty}}ds)\\
&\times\exp(C(1+\int^{t}_{0}P^{16}(\|\rho(s,\cdot)\|_{L^{\infty}})ds)\int^{t}_{0}\phi_{2}(\|\rho(s,\cdot)\|_{L^{\infty}})ds)),
\end{aligned}
\label{5impa1}
\end{equation}
where $\phi_{1}\in C^{0}(\R_{+},\R_{+}^{*})\cap C^{1}(0,\infty)$ such that $\phi(s)\geq \e_{0}s$ for some positive $s$.\\
Next we have for $\phi_{3}\in{\cal F}$ with $\alpha$ enough big:
$$
\begin{aligned}
 \int^{t}_{0}P^{4}(\|\rho(s,\cdot)\|_{L^{\infty}})ds)\int^{t}_{0}\phi_{2}(\|\rho(s,\cdot)\|_{L^{\infty}})ds
&\leq(\int \phi_{3}(\|\rho(s,\cdot)\|_{L^{\infty}})ds)^{2},\\
&\leq \int \phi^{2}_{3}(\|\rho(s,\cdot)\|_{L^{\infty}})ds.
\end{aligned}
$$
Finally from (\ref{5impa1}), we obtain:
\begin{equation}
\begin{aligned}
&A(t)\leq C\exp(C\int^{t}_{0}\phi_{4}(\|\rho(s,\cdot)\|_{L^{\infty}})ds)).
\end{aligned}
\label{5impa2}
\end{equation}
\subsubsection*{Control of $\sup_{0<t\leq T}f^{N}(t)\int\rho|\dot{u}|^{2}(t,x)dx+\int\int f^{N}(s)|\n\dot{u}|^{2}dxds$}
In the sequel, we want to obtain estimate on $\n u$ in $L^{1}_{T}(BMO)$, that's why we need of additional regularity estimates. In fact more precisely we want get more regularity on $v_{1}$ the \textit{effective velocity} define in \cite{H3}. We recall here briefly the definition of $v_{1}$, the idea is that new variable check an heat equation with additional source terms. To achieve it, we need to include the pressure term in the study of the linearized equation of the momentum equation as in \cite{H3}. For that, we will try to express the gradient of the pressure as a Laplacian term, so we have to solve:
$$\D v=\n P(\rho).$$
Let ${\cal E}$ be the fundamental solution of the Laplace operator. We will set in the sequel: $v=\n{\cal E}*\big(P(\rho)-P(\bar{\rho})\big)=\n\big({\cal E}*[P(\rho)-P(\bar{\rho})]\big)$ ( $*$ here means the operator of convolution). We verify next that:
$$
\begin{aligned}
\n{\rm div}v=\n\D \big({\cal E}*[P(\rho)-P(\bar{\rho})]\big)=\D\n\big({\cal E}*[P(\rho)-P(\bar{\rho})]\big)=\D v=\n P(\rho).
\end{aligned}
$$
By this way we can now rewrite the momentum equation of (\ref{1}) as:
$$\p_{t}u+u\cdot \n u-\frac{\mu}{\rho}\D\big(u-\frac{1}{\nu}v\big)-\frac{\lambda+\mu}{\rho}\n{\rm div}\big(u-\frac{1}{\nu}v\big)=f,$$
with $\nu=2\mu+\lambda$. We now want to calculate $\p_{t}v$, by the transport equation we get:
$$\p_{t}v=\n{\cal E}*\p_{t}P(\rho)=-\n {\cal E}*\big(P^{'}(\rho){\rm div}(\rho u)\big).$$
By setting  $v_{1}=u-\frac{1}{\nu}v$ which is called the effective velocity, we can now rewrite the system (\ref{1}) as follows:
\begin{equation}
\begin{cases}
\begin{aligned}
&\p_{t}q+(v_{1}+\frac{1}{\nu}v)\cdot\n q+\frac{1}{\nu}P^{'}(1)q =-(1+q){\rm div}v_{1}\\
&\hspace{3cm}-\frac{1}{\nu}(P(\rho)-P(1)-P^{'}(1))-\frac{1}{\nu}q(P(\rho)-P(1)),\\
&\p_{t}v_{1}-\frac{1}{1+q}{\cal A}v_{1}=f-u\cdot\n u+\frac{1}{\nu}\n(\D)^{-1}\big(P^{'}(\rho){\rm div}(\rho u)\big),\\
&q_{/ t=0}=a_{0},\;(v_{1})_{/ t=0}=(v_{1})_{0},
\end{aligned}
\end{cases}
\label{0.7}
\end{equation}
By this way the coupling between $v_{1}$ and the pressure disappears. Your goal now is to prove some regularity results on $v_{1}$. More precisely we want to prove that $\n v_{1}\in L^{1}_{T}(B^{\N+\e}_{2,\infty})$ such that as  $L^{1}_{T}(B^{\N+\e}_{2,\infty})\h L^{1}_{T}(L^{\infty})$ by proposition \ref{interpolation}, $v_{1}$ is Lipschitz. We recall that for obtaining strong solution, it is enough to obtain that $u$ is Lipschitz. We will see this point in the sequel.\\
\\
To obtain such estimates on $v_{1}$ we follow the ideas of the proof of D. Hoff in \cite{5H2}. For the completeness of the proof, we would like recall the main arguments used by D. Hoff in  \cite{5H2}.  We have then to derive estimates for the terms $f(t)^{2}\int_{\T^{N}}|\dot{u}|^{2}(t,x)dx$ and $\int^{t}_{0}\int_{\T^{N}}f{N}(s)|\n\dot{u}|^{2}dxds$. First we rewrite the momentum equation on the following form:
$$\rho\dot{u}-\mu\D u-(\lambda+\mu)\n{\rm div}u+\n P(\rho)=\rho g.$$
We apply to the momentum equation the operator $\frac{d}{dt}=\p_{t}+u\cdot\n$, we recall some elementary computations:
$$
\begin{aligned}
\frac{d}{dt}\rho\dot{u}^{j}&=\rho \frac{d}{dt}\dot{u}^{j}+(\p_{t}\rho )\dot{u}^{j}+\dot{u}^{j}\sum_{k}\p_{k}\rho \;u^{k},\\
&=\rho \frac{d}{dt}\dot{u}^{j}-\rho{\rm div}u \dot{u}^{j} ,
\end{aligned}
$$
We have next:
$$
\begin{aligned}
\mu\frac{d}{dt}\D u^{j}&=\mu\p_{t}\D u^{j}+\sum_{k}\p_{k}\D u^{j} u^{k},\\
&=\mu(\p_{t}\D u^{j}+{\rm div}(\D u^{j} u)-\D u^{j}{\rm div}u),
\end{aligned}
$$
and (where $D={\rm div}u$):
$$
\begin{aligned}
(\lambda+\mu)\frac{d}{dt}\p_{j}{\rm div}u&=(\lambda+\mu)(\p_{t}\p_{j}D+{\rm div}(\p_{j}D u)-\p_{j}D{\rm div}u),
\end{aligned}
$$
We finally obtain:
\begin{equation}
\begin{aligned}
&\rho\frac{d}{dt}\dot{u}^{j}+\p_{j}\p_{t}P(\rho)+{\rm div}\p_{j}P(\rho)=\mu(\p_{t}\D u^{j}+{\rm div}(\D u^{j} u))\\
&\hspace{5cm}+
(\lambda+\mu)(\p_{t}\p_{j}D+{\rm div}(\p_{j}D u))+\rho\frac{d}{dt}g^{j}.
\end{aligned}
\label{b2.10}
\end{equation}
We have then by summing over $j$:
\begin{equation}
\begin{aligned}
&\frac{1}{2}f(t)^{N}\int_{\T^{N}}\rho(t,x)|\dot{u}(t,x)|^{2}dx=\int^{t}_{0}\int_{\T^{N}}Nf^{N-1}(s)
f^{'}(s)\rho|\dot{u}|^{2}
dxds\\
&+\int^{t}_{0}\int_{\T^{N}} f(s)^{N}\dot{u}^{j}\big[-(\p_{j}\p_{t}P+{\rm div}(\p_{j}Pu)+\mu[\D\p_{t}u^{j}+{\rm div}(\D u^{j}u)]
\\
&\hspace{4cm}+(\lambda+\mu)[\p_{j}\p_{t}D+{\rm div}(\p_{j}D u)]+\rho\dot{g}^{j}\big]dxds.
\end{aligned}
\label{12.11}
\end{equation}
Since $f^{N-1}(s)
f^{'}(s)\leq f(s)$ we can apply (\ref{5impa2}) to bound the first term on the right.
Next by integrations by part we get:
$$
\begin{aligned}
-\int_{0}^{t}\!\!\int_{\T^{N}} f(s)^{N}\dot{u}^{j}(\p_{j}\p_{t}P+{\rm div}(\p_{j}Pu))dsdx& =
\int_{0}^{t}\!\!\int_{\T^{N}} f(s)^{N}(\p_{j}\dot{u}^{j}\p_{t}P+\p_{k}\dot{u}^{j}\p_{j}P u^{k}dsdx,\\
&=\int\!\!\int f(s)^{N}P^{'}(\p_{j}\dot{u}^{j}\p_{t}\rho+\p_{k}\dot{u}^{j}\p_{j}\rho u^{k})dsdx
\end{aligned}
$$
$$
\begin{aligned}
&\int_{0}^{t}\int_{\T^{N}} f(s)^{N}P^{'}[-\p_{j}\dot{u}^{j}(\rho\p_{k}u^{k}+\p_{k}\rho u^{k})+\p_{k}\dot{u}^{j}\p_{j}\rho u^{k}]dxds,\\
&=-\int_{0}^{t}\int_{\T^{N}} f(s)^{N}[P^{'}\rho D\p_{j}\dot{u}^{j}+\p_{k}Pu^{k}\p_{j}\dot{u}^{j}-\p_{j}Pu^{k}\p_{k}\dot{u}^{j}]dxds,\\
&=-\int_{0}^{t}\int_{\T^{N}}  f(s)^{N}[P^{'}\rho D\p_{j}\dot{u}^{j}-P(D\p_{j}\dot{u}^{j}-\p_{j}u^{k}\p_{k}\dot{u}^{j})]dxds.
\end{aligned}
$$
We bound therefore therefore the previous term by using H\"older inequalities and Young's inequality:
$$
\begin{aligned}
&C\big(\int^{t}_{0}\int_{\T^{N}}P(\rho)f(s)^{N}|\n u|^{2}dxds\big)^{\frac{1}{2}}\big(\int^{t}_{0}\int_{\T^{N}}f(s)^{N}|\n \dot{u}|^{2}dxds\big)^{\frac{1}{2}}\\
&\hspace{1.8cm}\leq C_{\e}\int^{t}_{0}\int_{\T^{N}} P(\rho)f(s)^{N}|\n u|^{2}(s,x)dxds+\e\int^{t}_{0}\int_{\T^{N}} f(s)^{N}|\n \dot{u}(s,x)|^{2}dxds,
\end{aligned}
$$
with $\e$ small enough for applying a bootstrap argument in the sequel. We are now interested in writing precisely the contribution of the third term on the right-hand side of (\ref{12.11}), we have then:
$$
\begin{aligned}
&-\mu\int^{t}_{0}\int_{\T^{N}} f^{N}[\n\dot{u}^{j}\cdot\n u^{j}_{t}+(\n\dot{u}^{j}\cdot u)\D u^{j}]dxds\\
&\hspace{1cm}=-\mu\int^{t}_{0}\int_{\T^{N}} f^{N}[\n\dot{u}^{j}\cdot(\n u^{j}_{t}+\n(\n\dot{u}^{j}\cdot u))+\dot{u}^{j}_{k}(u^{k}u^{j}_{ll}-(u^{j}_{l}u^{l})_{k})]dxds,\\
&\hspace{1cm}=-\mu\int^{t}_{0}\int_{\T^{N}}f^{N}[\n\dot{u}^{j}\cdot(\n u^{j}_{t}+\n(\n\dot{u}^{j}\cdot u))+\dot{u}^{j}_{k}(u^{k}u^{j}_{ll}-u^{j}_{lk}u^{l}-u^{j}_{l}u^{l}_{k})]dxds,\\
&\hspace{1cm}\leq -\mu\int^{t}_{0}\int_{\T^{N}} f^{N}|\n\dot{u}^{j}|^{2}+M\int\int f^{2}|\n u|^{2}(|\n\dot{u}|+|\dot{D}|)dxds.
\end{aligned}
$$
The last term on the right-hand side of (\ref{12.11}) may be treated as follows:
$$
\begin{aligned}
&-(\lambda+\mu)\int^{t}_{0}\int_{\T^{N}}f^{N}(\p_{j}\p_{t}D+{\rm div}(\p_{j}D u) dxds\leq
-(\lambda+\mu)\int^{t}_{0}\int_{\T^{N}}f^{N} \dot{D}^{2}dxds\\
&\hspace{7cm}+M\int^{t}_{0}\int_{\T^{N}} f^{N}|\n u|^{2}(|\n\dot{u}|+|\dot{D}|)dxds.
\end{aligned}
$$
From (\ref{12.11}) and by the preceding inequalities, it then follows by Young's inequalities that:
\begin{equation}
\begin{aligned}
&f(t)^{N}\int_{\T^{N}}\rho|\dot{u}|^{2}(t,x)dx+\int^{t}_{0}\int_{\T^{N}}f^{N}(s)(\mu|\n\dot{u}|^{2}+(\lambda+\mu)|\dot{D}|^{2})dxds\\
&\leq M\big[C_{0}+C_{\e}C_{0}\|P(\rho)\|_{L^{\infty}}+
\int^{t}_{0}\int_{\T^{N}}f(s)
\rho|\dot{u}|^{2}dxds\\
&\hspace{3cm}+\int^{t}_{0}\int_{\T^{N}}f^{N}(s)|\n u|^{4}dxds]+\int^{t}_{0}\int_{\T^{N}}f^{N}(s)\rho|\dot{g}|^{2}dxds,\\
&\leq M\big[C_{0}+C_{\e}C_{0}\|P(\rho)\|_{L^{\infty}}+
A(t)+\int^{t}_{0}\int_{\T^{N}}f^{N}(s)|\n u|^{4}dxds\\
&\hspace{7cm}+\int^{t}_{0}\int_{\T^{N}}f^{N}(s)\rho|\dot{g}|^{2}dxds].
\end{aligned}
\label{aa1}
\end{equation}
We recall here that $\omega={\rm curl}u$. Next from the momentum equation, we obtain as in the works of D. Hoff in \cite{5H2} by applying the operator ${\rm curl}$:
$$\mu|\n \w|^{2}=\mu({\rm div}(\w\n \w)+\p_{j}(\rho \w\dot{u}^{k})-\p_{k}(\rho w\dot{u}^{j})+\rho(\dot{u}^{j}\p_{k}\w-\dot{u}^{k}\p_{j}\w).$$
Integrating on $(0,T)\times\T^{N}$ and multiplying by $f(t)$ we get:
\begin{equation}
\begin{aligned}
\int^{t}_{0}\int_{\T^{N}}f(s)|\n \w|^{2}dxds&\leq M(\int^{t}_{0}\int_{\T^{N}} f(s)\rho^{2}|\dot{u}|^{2}(s,x)dsdx\\
&\leq M\|\rho\|_{L^{\infty}}A(t).
\end{aligned}
\label{aa2}
\end{equation}
Similarly we have:
\begin{equation}
\begin{aligned}
&\sup_{0<t\leq T}f(t)^{N}\int |\n\omega|^{2}dx\lesssim \int_{\T^{N}} f^{N}(t)\rho|\dot{u}|^{2}(t,x)dsdx+\int_{\T^{N}}f(t)^{N}\rho|w|^{2}(t,x)dx
\\
&\hspace{4cm}\lesssim \int_{\T^{N}} f^{N}(t)\rho|\dot{u}|^{2}(t,x)dsdx+C\|\rho\|_{L^{\infty}},\\
&\lesssim C_{0}+C_{\e}C_{0}\|P(\rho)\|_{L^{\infty}}+
A(t)+\int^{t}_{0}\int_{\T^{N}}f^{N}(s)|\n u|^{4}dxds+C\|\rho\|_{L^{\infty}}.
\end{aligned}
\label{aa3}
\end{equation}
To complete the estimates (\ref{aa1}) and (\ref{aa3}), we will need to estimate the following term
$\int^{t}_{0}\int_{\T^{N}}f^{N}(s)|\n u|(s,x)^{4}dxdt$. After we will come back to the inequality  (\ref{5impa2}) and we will prove the existence part of theorem \ref{theo1}.
\subsubsection*{Control of $\int^{t}_{0}\int_{\T^{N}}f^{N}(s)|\n u|(s,x)^{4}dxdt$}
We have then by splitting  the term $\n u$ as follows:
$$\n u=\n(\D)^{-1}{\rm div}\omega+\frac{1}{2\mu+\lambda}[\n(\D)^{-1}\n G+\n(\D)^{-1}\n P(\rho)].$$
we have then by Calderon-Zygmund theory:
\begin{equation}
\int^{t}_{0}\int_{\T^{N}} f^{N}|\n u|^{4}dsdx\lesssim \int^{t}_{0}\int_{\T^{N}} f^{N}\big((RG)^{4}+\omega^{4})(s,x)+f^{N}(s)RP(\rho)^{4}(s,x)\big)dxds,
\label{vorti}
\end{equation}
with $R$ a pseudo-differential operator of order $0$. Let focus us on the case $N=3$. We can apply Gagliardo-Nirenberg, we have then:
\begin{equation}
\begin{aligned}
&\int^{t}_{0}\int_{\T^{N}}f^{3}(s)(RG)^{4}(s,x)dsdx\leq\int^{t}_{0}f(s)^{3}\big(\int_{\T^{N}} G^{2}dx\big)^{\frac{1}{2}}\big(\int_{\T^{N}}|\n G|^{2}dx\big)^{\frac{3}{2}}dt,\\
&\hspace{2cm}\leq\|\sqrt{f(t)} G\|_{L^{\infty}_{t}(L^{2})}\|f(t)^{\frac{3}{2}} \n G\|_{L^{\infty}_{t}(L^{2})}
\int_{0}^{t}\int_{\T^{N}} f(s)|\n G|(s,x)^{2}dsdx.
\end{aligned}
\label{22imp1}
\end{equation}
By the definition of $G$ we have easily:
$$
\begin{aligned}
f(t)\int_{T^{N}} G^{2}(t,x)dx\leq M[\|P(\rho)\|_{L^{\infty}(L^{2})}+A(t)].
\end{aligned}
$$
Moreover as $(\lambda+2\mu)\D G={\rm div}(\rho (\dot{u}+g))$, we have by classical estimates on the elliptic system:
\begin{equation}
\begin{aligned}
f(t)^{3}\int_{\T^{N}}|\n G|^{2}dx&\leq M f(t)^{3}\|\rho\|^{\frac{1}{2}}_{L^{\infty}}(\int_{\T^{N}}\rho|\dot{u}|^{2}(t,x)dx\\
&\hspace{4cm}+
\|\rho\|^{\frac{1}{2}}_{L^{\infty}}\int_{\T^{N}}|g|^{2}(t,x)dx).
\end{aligned}
 \label{aaa1}
\end{equation}
We finally get by using  (\ref{22imp1}) and (\ref{aaa1}):
\begin{equation}
\begin{aligned}
&\int^{t}_{0}\int_{\T^{N}}f^{3}(s)(RG)^{4}(s,x)dsdx\leq M \|\rho\|^{\alpha}_{L^{\infty}}(1+A(t)^{2})
\end{aligned}
\label{22imp11}
\end{equation}
with $\alpha>0$.\\
A similar argument may be applied to the vorticity term, so that we have:
\begin{equation}
\int^{t}_{0}\int_{\T^{N}}f^{3}(s)|\n u|^{4}(s,x)dsdx\leq M\|\rho\|_{L^{\infty}}^{\alpha}(1+A^{2}(t)+\|g\|_{L^{\infty}(L^{2})}).
\end{equation}
From (\ref{5impa2}), (\ref{aa1}) and (\ref{aa3}), we can conclude that:
\begin{equation}
\begin{aligned}
&B(t)f(t)^{2}\int_{\T^{N}}\rho|\dot{u}|^{2}(t,x)dx+\int^{t}_{0}\int_{\T^{N}}f^{2}(s)(\mu|\n\dot{u}|^{2}+(\lambda+\mu)|\dot{D}|^{2})dxds\\
&\hspace{4cm}\leq M\big[C_{0}+C_{\e}C_{0}\|P(\rho)\|_{L^{\infty}}+CC_{e}A(t)(1+A(t)^{9})].\\[2mm]
&\sup_{0<t\leq T}f(t)^{2}\int |\n\omega|^{2}dx
\leq M(C_{0}+C_{\e}C_{0}\|P(\rho)\|_{L^{\infty}}+
A(t)C\|\rho\|_{L^{\infty}}\\
&\hspace{9cm}+\e B(t)+CC_{e}A(t)^{10}).
\end{aligned}
\label{aa33}
\end{equation}
We can remark that all the inequalities (\ref{5impa2}), (\ref{aa2}) and (\ref{aa33}) only depend on the control of $\|\rho\|_{L^{\infty}}$. In the following subsection, we are going to explain how we can control the $L^{\infty}$ norm of the density  $\rho$ in finite time.
\subsubsection*{Conclusion}
We will treat for simplicity only the case $N=3$.
We want here to explain how to obtain  $L^{\infty}$ estimate on the density $\rho$ in finite time.
From (\ref{29}), we have:
\begin{equation}
\begin{aligned}
\log(\rho(t,x))\leq&C\big(\log(\|\rho_{0}\|_{L^{\infty}})+\|(\D)^{-1}{\rm div} m_{0}\|_{L^{\infty}}+\|(\D)^{-1}{\rm div}(\rho u)\|_{L^{\infty}}\\
&+C \int^{t}_{0}\|[u_{j},R_{i}R_{j}](\rho u_{i})(s,\cdot)\|_{B^{1}_{N+\e,1}}ds\big).
\end{aligned}
\label{denscrucial}
\end{equation}
We have then by Sobolev embedding with $\e>0$, (\ref{againvitesse}) and let $p\geq 6+\e$ with $\e>0$ and $q=\frac{\gamma p}{p-1}$ such that $\frac{1}{p}+\frac{1}{q}<\frac{1}{3}$ which means $\gamma\geq 6$ :
\begin{equation}
\begin{aligned}
\|(\D)^{-1}{\rm div}(\rho u)\|_{L^{\infty}_{T}(L^{\infty})}&\leq\|\rho u\|_{L^{\infty}_{T}(L^{3+\e})},\\
&\leq \|\rho\|_{L^{\infty}_{T}(L^{\gamma})}^{\frac{p-1}{p}}\|\rho^{\frac{1}{p}}u\|_{L^{\infty}_{T}(L^{p})}\\
&\leq C\|\rho\|_{L^{\infty}_{T}(L^{\gamma})}^{\frac{p-1}{p}}(1+\int^{t}_{0}\|P(\rho(s,\cdot))\|^{6+\e}_{L^{\infty}}ds).
\end{aligned}
\label{333final}
\end{equation}
We can proceed similarly in the case $N=2$ and $\gamma\geq 1$ suffices.\\
We now want to bounded $\|\n u\|_{L^{6}}$ to control in the sequel the last term on the right hand side in (\ref{denscrucial}). We have then by Calderon Zygmund theory and the fact that $\D G={\rm div}(\rho\dot{u})+{\rm div}(\rho g)$, by using classical estimates for elliptic equation we have:  
$$
\begin{aligned}
\|\n u\|_{L^{6}}&\lesssim\|G\|_{L^{6}}+\|P(\rho)\|_{L^{\infty}}+\|\w\|_{L^{6}}\\
&\leq \|\rho\|_{L^{\infty}}^{\frac{1}{2}}\|\rho^{\frac{1}{2}}\dot{u}\|_{L^{2}}+\|\rho\|_{L^{\infty}}\|g\|_{L^{2}}+
\|P(\rho)\|_{L^{\infty}}+\|\n w\|_{L^{2}}.
\end{aligned}
$$
We have then by interpolation for $\alpha>0$  and  $\theta$ small enough:
$$
\begin{aligned}
&f(t)^{\frac{1-\theta}{2}}\|\n u\|_{L^{6-\alpha}}\leq\|\n u\|_{L^{2}}^{\theta}(\|\rho\|_{L^{\infty}}^{\frac{1}{2}}\sqrt{f(t)}\|\rho^{\frac{1}{2}}
\dot{u}\|_{L^{2}}+\|\rho\|_{L^{\infty}}\|g\|_{L^{2}}\\
&\hspace{5cm}+\|P(\rho)\|_{L^{\infty}}+\sqrt{f(t)}\|\n\w\|_{L^{2}})^{1-\theta}
\end{aligned}
$$
We obtain then from (\ref{29imp}) with $\alpha>0$ small enough,  by Young inequality with $p=\frac{2}{1-\theta}$ and $q=\frac{2}{\theta}$ and H\"older inequalities ($\frac{p-1}{\gamma p}+\frac{1}{p}<\frac{1}{6}$ as $\gamma>\frac{6(p-1)}{p-6}$):
$$
\begin{aligned}
&\int^{t}_{0}\|[u_{j},R_{i}R_{j}](\rho u_{i})(s,\cdot)\|_{B^{1}_{N+\e,1}}ds\leq\int^{t}_{0}\|\rho u(s,\cdot)\|_{L^{6+\e}}\|\n u(s,\cdot)\|_{L^{6-\alpha}}ds\\
&\leq \int^{t}_{0}\|\rho u(s,\cdot)\|_{L^{6+\e}}\frac{1}{f(s)^{\frac{1-\theta}{2}}}\|\n u\|_{L^{2}}^{\theta}\big(\|\rho\|_{L^{\infty}}^{\frac{1}{2}}\sqrt{f(t)}\|\rho^{\frac{1}{2}}
\dot{u}\|_{L^{2}}+\|\rho\|_{L^{\infty}}\|g\|_{L^{2}}\\
&\hspace{5cm}+\|P(\rho)\|_{L^{\infty}}+\sqrt{f(t)}\|\n\w\|_{L^{2}}\big)^{1-\theta}dx\\
&\lesssim \|\rho\,u(s,\cdot)\|_{L^{\infty}_{t}(L^{6+\e})}\|\frac{1}{f(s)^{\frac{1-\theta}{2}}}\|_{L^{\frac{2}{1-\theta}-\alpha}}\,A(t)^{\frac{1-\theta}{2}}(\int\|\n u(s,\cdot)\|^{2}ds)^{\frac{\theta}{2}},\\
&\lesssim\|\rho\, u(s,\cdot)\|^{\frac{2}{\theta}}_{L^{\infty}_{t}(L^{6+\e})}+A(t)\lesssim\big(\|\rho^{1-\frac{1}{p}}\|_{L^{\infty}(L^{\frac{\gamma p}{p-1}})}\|\rho^{\frac{1}{p}} u(s,\cdot)\|\big)^{\frac{2}{\theta}}_{L^{\infty}_{t}(L^{6+\e})}+A(t),\\
&\lesssim\|\rho^{\frac{1}{p}} u(s,\cdot)\|\big)^{\frac{2}{\theta}}_{L^{\infty}_{t}(L^{6+\e})}+A(t)\lesssim 1+\int^{t}_{0}\phi(\|\rho(s)\|_{L^{\infty}})ds+A(t).
\end{aligned}
$$
By using the previous inequality and (\ref{29}), (\ref{5impa2}) and (\ref{333final}) we conclude that:
$$
\begin{aligned}
\log\rho(t,x)&\leq C+\int^{t}_{0}\phi(\|\rho(s)\|_{L^{\infty}})ds+A(t),\\
&\leq C_{t}+C\exp(\int^{t}_{0}\phi_{1}(\|\rho(s,\cdot)\|_{L^{\infty}})ds))
\end{aligned}
$$
with $\phi\in{\cal F}$.
We obtain then:
\begin{equation}
 |\rho(t,x)|\leq C\exp(\int^{t}_{0}\phi_{1}(\|\rho(s)\|_{L^{\infty}})ds).
\label{115}
\end{equation}
Denoting by $w(t)$ the right-hand side of (\ref{115}), we conclude that:
$$\frac{d}{dt}w(t)\leq C \phi_{2}(\|\rho(s)\|_{L^{\infty}})w(t)\leq C w(t)(M+C_{1} w^{\beta}(t)),$$
because $\phi_{2}\in{\cal F}$. We have then as $w(t)\geq C$:
$$\frac{\frac{d}{dt}w(t)}{w^{1+\beta}(t)}\leq M_{1}$$
so that there exists $T_{0}$ such that for all $T<T_{0}$
$$\|\rho\|_{L^{\infty}((0,T)\times\T^{N})}\leq C_{T}.$$
\subsubsection*{Proof that $\frac{1}{\rho}$ is in $L^{\infty}$ when $\rho$ is in $L^{\infty}$}
We npw want to check that if we assume that $\rho_{0}$ is bounded away from zero then on $(0,T_{0})$ the density $\rho$ remains bounded away from zero. From (\ref{26}), we have:
\begin{equation}
\begin{aligned}
&\log\rho(t,x)\geq\ln\rho_{0}(x)-\|(\D)^{-1}{\rm div}m_{0}\|_{L^{\infty}}-\|(\D)^{-1}{\rm div}m_{0}\|_{L^{\infty}_{T_{0}}(L^{\infty})}\\
&-\int^{t}_{0}\|P(\rho(s,\dot)\|_{L^{\infty}}ds+\int^{t}_{0}\int_{\T^{N}}P(\rho(s,\cdot))dxdt-\int^{t}_{0}\|[u_{j},R_{i}R_{j}](\rho u_{i})(s,\cdot)\|_{L^{\infty}}ds.
\end{aligned}
 \label{2226}
\end{equation}
As  we assume here that $\frac{1}{\rho_{0}}$ belongs to $L^{\infty}$ and we have shown that all the terms in the right hand side
of (\ref{2226}) are bounded, we can conclude that $\frac{1}{\rho}$ is in $L^{\infty}_{T_{0}}(L^{\infty})$.
\subsubsection*{Proof of (\ref{b1.22}) and (\ref{b1.23})}
We now want to get supplementary estimates on the solution, and more precisely we would like to prove a control of the gradient of the velocity $\n u$ in $L^{1}(BMO(\T^{N}))$. In the sequel we will call such solution weak-strong solution, indeed to control $\n u$ in $L^{1}(BMO(\T^{N}))$ is not enough in order to obtain strong solution but we will see that by adding  slightly some regularity hypothesis on the initial density it will be enough. For more details we refer to the next section . \\
To obtain this additional regularity on the velocity, we need of new estimates on $u$.
To do this, we are widely inspired by  the technics introduced by D. Hoff in \cite{Hoffn2} that we will recall for the completeness of the proof. It can be also seen as an extension of \cite{H3}, indeed we will crucially use some interpolation results in order to cancel out the coupling between pressure and velocity. In \cite{H3} we are working directly with the \textit{effective velocity} which also allows us to obtain regularizing effects on $v_{1}$. In some sense the methods used in \cite{H3} and \cite{Hoffn2} are completely dual, in \cite{H3} the use of effective pressure avoids to use a direct argument interpolation. Furthermore in \cite{H3}, the obtained results are really in the scaling of the system.\\
As explained previously  in \cite{H3}, we are able to prove that $v_{1}$ is Lipschitz with initial velocity critical for the scaling of the equations. More precisely we have $\n v_{1}$ in $\widetilde{L}^{1}_{T_{0}}(B^{\NN}_{p,1})\h L^{1}_{T_{0}}(L^{\infty})$ when $u_{0}$ and $\rho_{0}$
are critical for the scaling of the equation. More precisely $u_{0}$ is in $ B^{\frac{N}{p_{1}}-1}_{p_{1},1}$ and $\rho_{0}$ in $B^{\NN}_{p,1}$ with $1\leq p_{1}\leq p<+\infty$ and some extra conditions on $(p_{1},p)$. This method is maybe more direct but however is probably more technical because the intensive use of Besov space with very sharp hypothesis on the paraproduct laws. For more details we refer to \cite{H3}. \\
 We now would like to explain why in our context it seems more advantageous to use the technics of \cite{Hoffn2}. Indeed here for proving (\ref{b1.22}) and (\ref{b1.23}) in theorem \ref{theo1}, we want to have minimal hypothesis on the initial density, it means that we assume only $\rho_{0}\in L^{\infty}$ and $\frac{1}{\rho_{0}}\in L^{\infty}$. That is why we can not directly use  \cite{H3} where we need to suppose that $\rho_{0}$ in $B^{\NN}_{p,1}$  (what is a bit more regular). In return as in \cite{Hoffn2}, we have to suppose that $u_{0}$ is in $H^{\N-1+\e}$ what is subcritical in terms of scaling of the equation (in comparison the initial velocity is really critical in \cite{H3}).\\
For the completeness of the proof we would recall these different technics. We would also point out that an other difficulty  compared with  \cite{Hoffn2} is  to be very careful and very accurate inasmuch as we work with large initial data. \\
First , we mollify initial data satisfying the conditions of theorem \ref{theo1} and then uses  the result of \cite{H3}
to obtain a solution $(\rho,u)$ defined at least for small time. We now want to derive estimates on the solution $(\rho,u)$ which do not depend of the mollifier process but  depend only on the initial data.
\\
In the sequel we will treat only by simplicity the case $N=3$.
Fixing the local in time solution $(\rho,u)$ described above on the interval $[0,T_{0}]$ with $T_{0}>0$, we therefore assume throughout this section that $C^{-1}\leq \rho\leq C$ with $C>0$. A crucial point is also to check that these estimates depend only of the condition $\rho$ in $L^{\infty}$.\\
We define a differential operator ${\cal L}$ acting on functions $w:[0,T]\times\T^{N}\rightarrow\T^{N}$ by
$${\cal L}w=\p_{t}(\rho w)+{\rm div}(\rho u\otimes w)-\mu\D w-\lambda\n{\rm div}w,$$
and we define $w_{1}$ and $w_{2}$ by:
\begin{equation}
\begin{aligned}
&{\cal L}w_{1}=0,\;\;\;\;\;\;\;\;\;\;\;\;\;\;{\cal L}w_{2}=-\n P(\rho),\\
&(w_{1})_{/t=0}=u_{0},\;\;\;\;\;(w_{2})_{/t=0}=0.
\end{aligned}
 \label{g2.1}
\end{equation}
We observe here that by uniqueness $w_{1}+w_{2}=u$. By energy estimates we obtain:
\begin{equation}
\sup_{0\leq t\leq T_{0}}\int_{\T^{N}}\rho(t,x)|w_{1}(t,x)|^{2}dx+\int^{T_{0}}_{0}\int_{\T^{N}}|\n w_{1}|^{2}dxdt\leq C\int_{\T^{N}}\rho_{0}|u_{0}|^{2}dx
 \label{g2.2}
\end{equation}
and:
\begin{equation}
\sup_{0\leq t\leq T_{0}}\int_{\T^{N}}\rho(t,x)|w_{2}(t,x)|^{2}dx+\int^{T_{0}}_{0}\int_{\T^{N}}|\n w_{2}|^{2}dxdt\leq CT \sup_{0\leq t}|P(\rho(t,\cdot)|^{2}.
 \label{g2.3}
\end{equation}
We shall derive (\ref{b1.22}) and (\ref{b1.23}) simultaneously as consequences of estimates for the following quantities (we can observe that this method is also to rely to the  Kato spaces for incompressible Navier-Stokes, see \cite{Kato}) :
$$\sup_{0\leq t\leq T_{0}}t^{1-k}\int_{\T^{N}}|\n w_{1}(t,x)|^{2}dx+\int^{T_{0}}_{0}\int_{\T^{N}} t^{1-k}\rho|\dot{w}_{1}|^{2}dxdt,$$
for $k=0,1$ and,
$$\sup_{0\leq t\leq1}\int|\n w_{2}(t,x)|^{2}dx+\int^{1}_{0}\int\rho|\dot{w}_{2}|^{2}dxdt.$$
To derive these bounds, we multiply equations (\ref{g2.1}) by $\dot{w}_{1}$ and $\dot{w}_{2}$, respectively and integrate. The details which are nearly identical to those in the previous section are left to the reader. The essential point is only to see that these estimates depends only of $\rho$ in $L^{\infty}$. But as we have seen previously that $\frac{1}{\rho}$ is bounded in $L^{\infty}$ if $\frac{1}{\rho_{0}}$ and $\rho$ are respectively in $L^{\infty}$ and $L^{\infty}_{T_{0}}(L^{\infty}(\T^{N})$, we can consider in the sequel that $C^{-1}\leq \rho\leq C$.  And this last point do not depend on the mollifying process but only on the initial data and the fact that $\rho$ is in $L^{\infty}$. That is why in the sequel we do not care from the constant coming from $\|\frac{1}{\rho}\|_{L^{\infty}}$.\\
We obtain more precisely: 
\begin{equation}
\begin{aligned}
 &\frac{1}{2}\mu t^{k}\int_{\T^{N}}|\n w_{1}(\tau,x)|^{2}dx+\int^{t}_{0}\int_{\T^{N}}s^{k}|\dot{w}_{1}|^{2}dxds\leq\\
&\hspace{2cm}\frac{1}{2}\mu k\int^{t}_{0}\int_{\T^{N}}s^{k-1}|\n w_{1}|^{2}dxds+\int^{t}_{0}\int_{\T^{N}}s^{k}|\n w_{1}|^{2}|\n u|dxds,
\end{aligned}
 \label{g2.4}
\end{equation}
and
\begin{equation}
\begin{aligned}
 &\frac{1}{2}\mu\int_{\T^{N}}|\n w_{2}(t,x)|^{2}dx+\int^{t}_{0}\int_{\T^{N}}|\dot{w}_{2}|^{2}dxds\leq \int_{\T^{N}}P(\rho(t,x)){\rm div}w_{2}(t,x)dx\\
&\hspace{6cm}+\int^{t}_{0}\int_{\T^{N}}(|\n w_{2}|^{2}|\n u|+
|\n w_{2}||\n u|)dxd\tau .
\end{aligned}
 \label{g2.5}
\end{equation}
By proceeding exactly as in the previous section, we obtain the following results:
\begin{equation}
\sup_{0\leq t\leq T_{0}}\int_{\T^{N}}|\n w_{1}(t,x)|^{2}dx+\int^{T_{0}}_{0}\int_{\T^{N}}|\dot{w}_{1}|^{2}dxdt\leq C\|u_{0}\|_{H^{1}}^{2},
 \label{g2.7}
\end{equation}
\begin{equation}
\sup_{0\leq t\leq T_{0}}t\int_{\T^{N}}|\n w_{1}(t,x)|^{2}dx+\int^{T_{0}}_{0}\int t|\dot{w}_{1}|^{2}dxdt\leq C\|u_{0}\|_{L^{2}}^{2},
 \label{g2.8}
\end{equation}
\begin{equation}
\sup_{0\leq t\leq T_{0}}\int_{\T^{N}}|\n w_{2}(t,x)|^{2}dx+\int^{T_{0}}_{0}\int_{\T^{N}}|\dot{w}_{2}|^{2}dxdt\leq C C_{0},
 \label{g2.9}
\end{equation}
where $C_{0}$ depends only of the initial data and of $\|\rho\|_{L_{T_{0}}^{\infty}}$. Now since the solution operator $u_{0}\longrightarrow w_{1}(t,\cdot)$ is linear, we can apply a standard Riesz-Thorin interpolation argument to deduce from (\ref{g2.7}) and (\ref{g2.8}) that:
\begin{equation}
\sup_{0\leq t\leq T_{0}}t^{1-\beta}\int_{\T^{N}}|\n w_{1}(t,x)|^{2}dx+\int^{T_{0}}_{0}\int_{\T^{N}} t^{1-\beta}|\dot{w}_{1}|^{2}dxdt\leq C\|u_{0}\|_{H^{\beta}}^{2}.
 \label{g2.10}
\end{equation}
As $u=w_{1}+w_{2}$, we then conclude from (\ref{g2.9}) and (\ref{g2.10}) that:
\begin{equation}
\sup_{0\leq t\leq T_{0}}t^{1-\beta}\int_{\T^{N}}|\n u(t,x)|^{2}dx+\int^{T_{0}}_{0}\int_{\T^{N}} t^{1-\beta}|\dot{u}|^{2}dxdt\leq C
C_{0}\|u_{0}\|_{H^{\beta}}^{2}.
 \label{g2.11}
\end{equation}
The next step consists in obtaining bounds for the terms
$$\sup_{0\leq t\leq T_{0}}t^{2-\beta}\int_{\T^{N}}|\n u(t,x)|^{2}dx+\int^{T_{0}}_{0}\int_{\T^{N}} t^{2-\beta}|\n \dot{u}|^{2}dxdt$$
appearing in (\ref{b1.22}). To do this, we multiply the momentum equation of (\ref{1})  by $t^{2-\beta}\dot{u}$ and integrate.
The details are exactly as in the previous section, except now we apply the $\beta$ dependent smoothing rates established in (\ref{g2.10}). Combining
these bounds with (\ref{g2.2}), (\ref{g2.3}) and (\ref{g2.10}), we then obtain (\ref{b1.22}) for times $t\leq T_{0}$.\\
\\ 
To prove (\ref{b1.23}), we observe that for $k=0,1$,
$$\sup_{0\leq t\leq 1}\|w_{1}(t,\cdot)\|_{H^{k}}\leq C\|u_{0}\|_{H^{k}},$$
by (\ref{g2.2}) and (\ref{g2.7}). Thus:
$$\sup_{0\leq t\leq 1}\|w_{1}(t,\cdot)\|_{H^{\beta}}\leq C\|u_{0}\|_{H^{k}},$$
for $\beta\in [0,1]$. As $u=w_{1}+w_{2}$, and applying (\ref{g2.9}) we obtain that:
$$\sup_{0\leq t\leq 1}\|w_{1}(t,\cdot)\|_{H^{\beta}}\leq C C_{0},$$
and then for $r\in(2,\frac{6}{3-2\beta})$ in the case that $\beta>0$, that:
$$\sup_{0\leq t\leq 1}\|u(r,\cdot)-\widetilde{u}\|_{L^{r}}\leq C C_{0}.$$
This proves (\ref{b1.23}).
\subsubsection*{Regularity on the gradient of the velocity $u$ in $L^{1}_{T_{0}}\big((W^{1,\alpha}+BMO)(\T^{N})\big)$ with $\alpha>N$.}
\label{v1reg}
Here we want to examine the regularity of the gradient of the velocity and to prove that $\n u$ is in $L^{1}_{T_{0}}(BMO)$ for showing (\ref{c1.23}). To show some regularity estimates on $\n u$, we begin with verifying that the new variable $v_{1}$ introduced in \cite{H3} called ``effective velocity'' belongs to $L^{1}_{T_{0}}(W^{2,\alpha})$ with $\alpha>N$ which implies that $\n v_{1}\in L^{1}_{T_{0}}(L^{\infty})$. We recall here the definition of $v_{1}$ introduced in \cite{H3}. The idea of \cite{H3} was to introduce a variable $v_{1}$ which allows to cancel out the coupling between the velocity and the pressure in the momentum equation of (\ref{1}). In this goal, we need to integrate the pressure term in the study of the linearized equation of the momentum equation. To do this, we will try to express the gradient of the pressure as a Laplacian term, so we set:
$$\D v=\n P(\rho).$$
We have then $v=(\D)^{-1}\n P(\rho)$ with $(\D)^{-1}$ the inverse Laplacian with zero value on $\T^{N}$. In the sequel we will set:
$$v_{1}=u-\frac{1}{\lambda+2\mu}v.$$
We check then easily that:
$$
\begin{aligned}
&{\rm div}v_{1}=\frac{1}{\lambda+2\mu}G\;\;\mbox{and}\;\;{\rm curl}v_{1}=\omega.
\end{aligned}
$$
We have then:
\begin{equation}
 \begin{aligned}
\D u=&\n{\rm div}u+{\rm div}\omega=\n{\rm div}v_{1}+{\rm div}{\rm curl}v_{1}+\frac{1}{2\mu+\lambda}\n (P(\rho)),\\
=&\D v_{1}+(2\mu+\lambda)^{-1}\n (P(\rho)).
\end{aligned}
\label{11.22}
\end{equation}
We can easily show that $\int^{T_{0}}_{0}\|\n v_{1}\|_{L^{\infty}}dt<+\infty$ if (\ref{b1.22}) holds.
To see this, we apply standard elliptic theory 
on $v_{1}$ (indeed  $\D v_{1}=\frac{1}{\lambda+2\mu}\n G+{\rm div}\omega$) combined with the fact that $\D G ={\rm div}(\rho\dot{u}-\rho g)$ and $\mu\D\mathbb{P}u=\mathbb{P}(\rho\dot{u}-\rho g)$. By simplicity, we will consider only the case $N=3$. The case $N=2$ follows the same lines. For some $\alpha>3$ and $\e>0$ determined by $\alpha$ we have then by Sobolev embedding and Gagliardo-Nirenberg estimates:
\begin{equation}
\begin{aligned}
&\|\n v_{1}\|_{L^{\infty}}\leq C(\|\n G\|_{L^{\alpha}}+\|\n\omega\|_{L^{\alpha}}), \\
&\lesssim(\|\rho\dot{u}(t,\cdot)\|_{L^{\alpha}}+\|\rho g\|_{L^{\alpha}}+\|\n\omega\|_{L^{\alpha}}),\\
&\lesssim(\|\rho\|_{L^{\infty}}\|\dot{u}(t,\cdot)-\widetilde{\dot{u}}\|^{\frac{1-\e}{2}}_{L^{2}}\|\n\dot{u}(t,\cdot)\|^{\frac{1+\e}{2}}_{L^{2}}+
\|\rho\|_{L^{\infty}}\|g\|_{L^{\alpha}}+\|\n\omega\|_{L^{\alpha}}),\\
&\lesssim(\|\rho\|_{L^{\infty}}(\|\dot{u}(t,\cdot)\|^{\frac{1-\e}{2}}_{L^{2}}+\widetilde{\dot{u}}^{\frac{1-\e}{2}})\|\n\dot{u}(t,\cdot)\|^{\frac{1+\e}{2}}_{L^{2}}+
\|\rho\|_{L^{\infty}}\|g\|_{L^{\alpha}}+\|\n\omega\|_{L^{\alpha}}).
\end{aligned}
\label{1213}
\end{equation}
We recall here that $\widetilde{\rho\dot{u}}=\widetilde{\rho g}$, we have then as $\frac{1}{\rho}\geq C$ on $[0,T_{0}]$ and as $\rho$ belongs to $L^{\infty}$:
$$\widetilde{\dot{u}}\leq\frac{1}{\min_{(0,T_{0})\times\T^{N}}\rho(t,x)}\widetilde{\rho g}\lesssim \widetilde{g}.$$
So that by using (\ref{condig}), (\ref{b1.22}), (\ref{1213}) and one time more that $\D G ={\rm div}(\rho\dot{u}-\rho g)$ and $\mu\D\mathbb{P}u=\mathbb{P}(\rho\dot{u}-\rho g)$, we finally obtain:
\begin{equation}
\begin{aligned}
\int^{T_{0}}_{0}\|\n v_{1}(t)\|_{L^{\infty}}dt&\lesssim \int^{T_{0}}_{0}t^{\beta}(t^{1-s}\int_{\T^{N}}|\dot{u}|^{2}dx)^{\frac{1-\e}{4}}(t^{\sigma}\int_{\T^{N}}|\n\dot{u}|^{2}dx)^{\frac{1+\e}{4}}dt+1.
\end{aligned}
\label{4567}
\end{equation}
with $s=\N+\e-1$ ($\e>0$) and where $4\beta=(s-1)(1-\e)-(\sigma+\e)$.\\
From (\ref{11.22}) and the fact that $\frac{1}{\rho}$ is in $L^{\infty}$, we have that
$(t^{1-s}\int_{\T^{N}}|\dot{u}|^{2}dx)^{\frac{1-\e}{4}}(t^{\sigma}\int_{\T^{N}}|\n\dot{u}|^{2}dx)^{\frac{1+\e}{4}}$ is in $L^{2}_{T_{0}}$. We can then conclude from (\ref{4567}) that:
$$\int^{T_{0}}_{0}\|\n v_{1}(t)\|_{L^{\infty}}dt\leq C(\int^{T_{0}}_{0}t^{2\beta}dt)^{\frac{1}{2}}+C_{0}.$$
The above integral is therefore finite as $2\beta>-1$. 
A similar result result holds for $N=2$ with $s>0$. Thus for the solution constructed in the previous section, $\int^{T}_{0}\|\n v_{1}(t,\cdot)\|_{L^{\infty}}dt$ is finite if (\ref{b1.22}) holds with the additional conditions in comparison of the previous sections, $\inf\rho_{0}\geq c>0$ and $u_{0}\in H^{\N+\e-1}$ with $\e>0$.\\
More precisely we have proved in fact that:
\begin{equation}
 \n v_{1}\in L^{1}_{T}(W^{1,\,\alpha})\hookrightarrow L^{1}_{T}(B^{1+\e}_{N,1}).
\label{gainv1}
\end{equation}
with $\alpha=N+2\e$ where $\e>0$.\\
As $P(\rho)\in L^{\infty}$ we deduce from (\ref{11.22}) and the results of Calderon-Zygmund, that:
$$\n u \in L^{1}_{T_{0}}\big(BMO(\T^{N})\big).$$
\section{Lipschitz estimates on the velocity $u$}
\label{section6}
In the theorem \ref{corollaire1}, we have assumed additional hypothesis on the initial density, in particular the fact that $\rho_{0}\in B^{\e}_{\infty,\infty}$ with $\e>0$. The goal here will be to show that this $B^{\e}_{\infty,\infty}$ regularity on the density is conserved on $(0,T_{0})$. In particular it will allow us to prove that the velocity is Lipschitz. It will be then enough to obtain the uniqueness of the solution constructed in theorem \ref{theo1}. Furthermore we will observe that for preserving the Lipschitz estimate on the velocity it will be enough to control the norm $L^{\infty}$ of the density. It will give us then
the blow-up result of the theorem \ref{corollaire1}.
\subsubsection{Control of $\rho\in L^{\infty}_{T_{0}}(B^{\e}_{\infty,\infty})$ and of $\n u\in L^{1}_{T_{0}}(B^{\e}_{\infty,\infty})$ }
\label{section62}
We now want to estimate $\rho$ in $L^{\infty}(B^{\e}_{\infty,\infty})$ and to prove that this control depends only on the norm $L^{\infty}$ of the density. In view of proposition \ref{transport} where in our case $h(\rho)=P(\rho)$, $\sigma=\e$, $p=p_{1}=r=+\infty$, we have for all $t\in[0,T]$ and $0<\e<1$:
\begin{equation}
 \|\rho\|_{\widetilde{L}^{\infty}_{t}(B^{\e}_{\infty,\infty})}\leq e^{CV(t)}(\|\rho_{0}\|_{B^{\e}_{\infty,\infty}}+\int^{t}_{0}C \|\rho(\tau,\cdot)\|_{L^{\infty}}\|{\rm div}v_{1}(\tau,\cdot)\|_{B^{\e}_{\infty,\infty}}d\tau\big)),
\label{38}
\end{equation}
where $V(t)=\int^{t}_{0}\big(\|\n u(\tau)\|_{B^{0}_{\infty,\infty}\cap L^{\infty}}+\|{\rm div}v_{1}(\tau)\|_{B^{\e}_{\infty,\infty}}+\|\rho(\tau)\|^{s}_{L^{\infty}}+1\big)d\tau$,
where $s$ the smallest integer such that $P^{'}\in W^{s,\infty}$. We have seen by (\ref{gainv1}) that $\n v_{1}\in L^{1}(0,T, B^{\e}_{\infty,\infty})$ with $\e>0$ small enough. By Besov embedding we recall that $B^{0}_{\infty,\infty}$ belongs to $L^{\infty}$. The main difficulty now is to control $\n u\in L^{1}(0,T,L^{\infty})$, for this we recall that by definition of the effective velocity and by proposition \ref{interpolation}:
\begin{equation}
\begin{aligned}
\|\n u\|_{L^{1}_{T}(L^{\infty})}&\lesssim \|\n v_{1}\|_{L^{1}_{T}(B^{\e}_{\infty,\infty})}+\|P(\rho)\|_{L^{1}_{T}(B^{0}_{\infty,1})},\\
&\leq \|\n v_{1}\|_{L^{1}_{T}(B^{\e}_{\infty,\infty})}+\|\rho\|^{s}_{L^{\infty}_{T}(L^{\infty})}\|\rho\|_{L^{1}_{T}(B^{0}_{\infty,1})}.
\end{aligned}
\label{562}
\end{equation}
The idea now is to use logarithmic estimates that we will inject in the inequality (\ref{38}), by using proposition \ref{interpolationlog} we have:
$$
\begin{aligned}
\|\rho(t)\|_{B^{0}_{\infty,1}}\leq C\|\rho(t)\|_{B^{0}_{\infty,\infty}}\log (e+\frac{\|\rho(t)\|_{B^{\e}_{\infty,\infty}}}{\|\rho(t)\|_{B^{0}_{\infty,\infty}}}),
\end{aligned}
$$
and we recall the following inequality:
$$\forall x>0, \;\forall\delta>0,\;\log(e+\frac{\delta}{x})\leq\log(e+\frac{1}{x})(1+\log\delta).$$
We obtain then from the previous inequality:
\begin{equation}
\|\rho(t)\|_{B^{0}_{\infty,1}}\leq \|\rho(t)\|_{B^{0}_{\infty,\infty}}\big(1+\log(\|\rho(t)\|_{B^{\e}_{\infty,\infty}})\big)\log (e+\frac{1}{\|\rho(t)\|_{B^{0}_{\infty,\infty}}}),
\label{4354}
\end{equation}
Let $X(t)=\int^{t}_{0}\|\rho(s)\|_{B^{0}_{\infty,1}}ds$, we have then from (\ref{562}), proposition \ref{interpolation} and from the fact that $\rho\in L^{\infty}(L^{\infty})$:
\begin{equation}
V(t)\leq C\big(1+ X(t)+\int^{t}_{0}\big(\|\n v_{1}(\tau)\|_{B^{\e}_{\infty,\infty}}+\|{\rm div}v_{1}(\tau)\|_{B^{\e}_{\infty,\infty}})d\tau\big).
\label{43541}
\end{equation}
Combining (\ref{38}) and the previous inequality leads to:
$$
\begin{aligned}
&X(t)\leq \int^{t}_{0}\|\rho(s)\|_{B^{0}_{\infty,\infty}}\biggl(1+CV(t)+\\
&\hspace{1cm}\log\big(\|\rho_{0}\|_{B^{\e}_{\infty,\infty}}+\int^{t}_{0}C \|\rho(\tau,\cdot)\|_{L^{\infty}}\|{\rm div}v_{1}(\tau,\cdot)\|_{B^{\e}_{\infty,\infty}}d\tau\big)\biggl)\log (e+\frac{1}{\|\rho(s)\|_{B^{0}_{\infty,\infty}}})ds,\\
&\leq \int^{t}_{0}\|\rho(s)\|_{B^{0}_{\infty,\infty}}\biggl(1+CX(t)+C\int^{t}_{0}\big(\|\n v_{1}(\tau)\|_{B^{\e}_{\infty,\infty}}+\|{\rm div}v_{1}(\tau)\|_{B^{\e}_{\infty,\infty}}\big)d\tau\\
&+\log\big(\|\rho_{0}\|_{B^{\e}_{\infty,\infty}}+\int^{t}_{0}C \|\rho(\tau,\cdot)\|_{L^{\infty}}\|{\rm div}v_{1}(\tau,\cdot)\|_{B^{\e}_{\infty,\infty}}d\tau\big)\biggl)\log (e+\frac{1}{\|\rho(s)\|_{B^{0}_{\infty,\infty}}})ds.
\end{aligned}
$$
Applying Gr\"onwall inequality and inequality (\ref{gainv1}) shows that:
$$
\begin{aligned}
X(t)&\leq C_{t,0}\exp(C\int^{t}_{0}\|\rho(s)\|_{B^{0}_{\infty,\infty}} \log (e+\frac{1}{\|\rho(s)\|_{B^{0}_{\infty,\infty}}})ds),\\
&\leq C_{t,0}\exp(C\int^{t}_{0}(1+\|\rho(s)\|_{L^{\infty}})ds),
\end{aligned}
$$
where $C_{t,0}$ depends only of the time $t$ and the initial data. As $\rho\in L^{\infty}_{t}(L^{\infty})$, we conclude that $X(t)\leq C_{t}$ and by this way we have proved that:
\begin{equation}
 \|\rho\|_{L^{\infty}_{T}(B^{\e}_{\infty,\infty})}\leq C_{0,T},
\label{gainclas}
\end{equation}
where $C_{t,0}$ depends only of the time $T$ and the initial data. We want to point out that we have in fact proved the following assertion:
\begin{equation}
\n u\in \widetilde{L}^{1}(B^{\e}_{\infty,\infty}).
\label{regugradientvi}
\end{equation}
\subsubsection{Control of $\rho\in \widetilde{L}^{\infty}(B^{1}_{N,1})$ for $N=3$ and  $\rho\in \widetilde{L}^{\infty}(B^{1}_{N+\e,1})$ for $N=2$ (with $\e>0$ arbitrary small) when $P(\rho)=K\rho$ with $K>0$}
\label{section61}
We will treat only the case $N=3$, the case $N=2$ follows exactly the same lines in the proof. In this case, we need to show for the sequel that $\rho\in L^{\infty}(B^{1}_{N,1})$, and for this we proceed exactly as previous.
Indeed as $\rho_{0}\in B^{1}_{N,1}$, by proceeding as in the previous section we can show that $\rho$ is in $L^{\infty}_{T}(B^{1}_{N,1})$. It suffices in particular to apply the proposition \ref{transport}. In the general case $P(\rho)=a\rho^{\gamma}$ with $\gamma>1$ the fact to control
$\rho\in \widetilde{L}^{\infty}(B^{1}_{N,1})$ will be crucial to obtain the uniqueness of the solutions and in particular for using the results of P. Germain (see \cite{PG}).
\section{Proof of theorem \ref{theo1} and \ref{corollaire1}}
\label{section51}
\subsection{Existence of weak solutions for theorem \ref{theo1} and \ref{corollaire1}}
In the sequel for simplicity, we will treat only the case $N=3$, the case $N=2$ follows the same arguments.
\subsubsection*{Existence of weak solutions for theorem  \ref{theo1}}
\label{section5}
The above arguments of section \ref{section3}, \ref{section4} and \ref{section6} are not rigorous, since we have to assume that $(\rho,u)$ is a solution of system (\ref{1}) (but it is exactly what we want to prove). Furthermore we need that this solution $(\rho,u)$ is enough regular to apply the different estimates which use crucially integration by parts in particularly. That is why to overcome this difficulty,  we need to smooth out the data in order to get a sequence of local solutions $(\rho^{n},u^{n})_{n\in\mathbb{N}}$ on $[0,T_{n}]$ to (\ref{1}) by using the results  \cite{DW} or \cite{H3}. $T^{n}$ ihere corresponds to the lifespan of  the solution $(\rho_{n},u_{n})$. More precisely we assume that:
$$\frac{1}{\rho_{0}^{n}}\in L^{\infty},\;\rho_{0}^{n}\in B^{1+\e}_{N,1},\,u^{n}_{0}\in B^{0}_{N,1}\;\;\mbox{and}\;\;f^{n}\in \widetilde{L}^{1}(B^{0}_{N,1}),$$
and that $(\rho_{0}^{n},u^{n}_{0})\rightarrow_{n\rightarrow +\infty} (\rho_{0},u_{0})$ in the norm of the spaces in what belong the initial data $(\rho_{0},u_{0})$. To do that, it suffices to  smooth out the data as follows:
$$\rho_{0}^{n}=S_{n}\rho_{0},\;\;u_{0}^{n}=S_{n}u_{0}\;\;\;\mbox{and}\;\;\;f^{n}=S_{n}f.$$
The main difficulty now is to prove that $T_{n}$ goes not to $0$ when $n$ goes to infinity. To do this we can observe from section \ref{section3} and \ref{section4} that there exists a time $T_{0}>0$ independent on $n$ such that $(\rho^{n},u{n})_{n\in\mathbb{N}}$ verify uniformly in function of $n$ the estimates
(\ref{1.21}), (\ref{b1.22}), (\ref{b1.23}) and (\ref{c1.23}). Furthermore $(\rho{n})_{n\in\mathbb{N}}$ verifies uniformly in function of $n$ on the interval $(0,T_{0})$ the following control:
\begin{equation}
\|\rho^{n}\|_{L_{T_{0}}^{\infty}(L^{\infty}(\T^{N})}\leq C\;\;\;\mbox{and}\;\;\;
\|\frac{1}{\rho^{n}}\|_{L_{T_{0}}^{\infty}(L^{\infty}(\T^{N})}\leq C.
\label{solutionsapprochees}
\end{equation}
Now we suppose by the absurd that $T^{n}\rightarrow_{n\rightarrow +\infty}0$. It means that for $n$ enough big $T^{n}<T_{0}$, then on $(0,T^{n})$  $(\rho_{n},u_{n})_{n\in\mathbb{N}}$ verify uniformly in function of $n$ the estimates
(\ref{1.21}), (\ref{b1.22}), (\ref{b1.23}) and (\ref{c1.23}) and (\ref{solutionsapprochees}). From section \ref{section6} we can also prove as $\rho_{0}^{n}\in B^{1+\e}_{N,1}$ that for all $n$:
\begin{equation}
\rho^{n}\in \widetilde{L}^{\infty}(B^{1+\e}_{N,1})\;\;\;\mbox{and}\;\;\;\n u^{n}\in L^{1}_{T_{n}}(L^{\infty}).
\label{solutionsapprochees1}
\end{equation}
From the continuation theorem proved in \cite{DW} when we have a solution of (\ref{1}) which verifies (\ref{solutionsapprochees}) and  (\ref{solutionsapprochees1}) then we can extend the solution. It is then a contradiction with the fact that $T_{n}$ is the lifespan of the solution. We have then for all $n\in\mathbb{N}$, $T^{n}\geq T_{0}$.\\
\\
We have then prove that the solutions $(\rho^{n},u^{n})_{n\in\mathbb{N}}$ of system (\ref{1}) with initial data $(\rho_{0}^{n},u_{0}^{n})_{n\in\mathbb{N}}$ exist on the interval $(0,T_{0})$ with $T_{0}>0$. Furthermore $(\rho_{0}^{n},u_{0}^{n})_{n\in\mathbb{N}}$  verify uniformly in $n$ estimates
(\ref{1.21}), (\ref{b1.22}), (\ref{b1.23}) and (\ref{c1.23}) and \ref{solutionsapprochees}. It is easy by using the results of Lions in \cite{13}, of Feireisl et al in \cite{5F3} or of Novotn\`y and stra\v{s}kraba in \cite{NS} to conclude that $(\rho^{n},u^{n})_{n\in\mathbb{N}}$ goes to a solution $(\rho,u)$ of system {\ref{1}) and that $(\rho,u)$ checks (\ref{1.21}), (\ref{b1.22}), (\ref{b1.23}) and (\ref{c1.23}) and $\rho\in L^{\infty}_{T_{0}}(L^{\infty}(\T^{N}))$.\\
The only point where we need to be careful is the case when $\gamma<\NN$ (in fact it is possible only when $N=2$), indeed in this case we can not use \cite{5F3} but with all the estimates (\ref{1.21}), (\ref{b1.22}), (\ref{b1.23}) and (\ref{c1.23})  uniformly verifying by $(\rho^{n},u^{n})_{n\in\mathbb{N}}$ and (\ref{solutionsapprochees}) it is an easy exercise to conclude.
\subsubsection*{Existence of weak solutions for theorem  \ref{corollaire1}}
The proof in this case follows exactly the same lines than in the previous section.
\subsection{Uniqueness for theorem \ref{corollaire1} when $P(\rho)=a\rho^{\gamma}$ with $\gamma>1$}
We now discuss the uniqueness of the solutions of theorem \ref{theo1}. For this we want to use the result of P. Germain \cite{PG} which is a result of weak-strong uniqueness. In the sequel we will note $(\rho_{1},u_{1})$ the solution of the theorem \ref{theo1} which exits on the time interval $[0,T_{0}]$.
We have shown that our solution check $\rho\in L^{\infty}(L^{\infty})$. By theorem \ref{theo1}, we obtain that our solution verify the following inequalities:
\begin{equation}
 \begin{aligned}
&\sup_{0<t\leq+\infty}\int_{\T^{N}}[\frac{1}{2}\rho(t,x)|u(t,x)|^{2}+|P(\rho(t,x))|+\sigma(t)|\n u(t,x)|^{2}dx\\
&+\sup_{0<t\leq +\infty}\int_{\T^{N}}[\frac{1}{2}\rho(t,x)f(t)^{N}(\rho|\dot{u}(t,x)|^{2}+|\n \omega(t,x)|^{2})dx\\
&+\int^{+\infty}_{0}\int_{\T^{N}}[|\n u|^{2}+f(s)\rho|\dot{u}|^{2}+|\omega|^{2})+\sigma^{N}|\n \dot{u}|^{2}]dxdt\\
&\hspace{9cm}\leq C(C_{0}+C_{f})^{\theta},
 \end{aligned}
\label{51.21}
\end{equation}
and we obtain moreover:
\begin{equation}
\begin{cases}
 \begin{aligned}
&\sqrt{\rho}\p_{t}u\in L^{2}_{t}(L^{2}(\T^{N})),\\
&\sqrt{t}{\cal P}u\in L^{2}_{T}(H^{2}(\T^{N})),\\
&\sqrt{t}G=\sqrt{t}[(\lambda+2\mu){\rm div}u-P(\rho)]\in L^{2}_{T}(H^{1}(\T^{N})),\\
&\sqrt{t}\n u\in L^{\infty}_{T}(L^{2}(\T^{N})),
 \end{aligned}
\end{cases}
\label{513}
\end{equation}
Now we assume that there exists two solutions of system (\ref{1}) $(\rho,u)$ and $(\rho_{1},u_{1})$ in the class of the solution of theorem \ref{corollaire1} with the same initial data $(\rho_{0},u_{0})$. Furthermore $(\rho_{0},u_{0})$ verify the conditions of the theorem \ref{corollaire1}. We now want to prove that
$(\rho,u)=(\rho_{1},u_{1})$ on $[0,T_{0}]$. To do this, we will use the result of P. Germain in \cite{PG}. To see this we have just to verify that $(\rho_{1},u_{1}$ verify the conditions of the theorem 2.2 of \cite{PG}.
For simplicity we will prove only the result for $N=3$. Before proving this assertion we would like to recall the theorem 2.2 of \cite{PG}.
\begin{theorem}
\label{theoPG}
Take initial data such that:
$$\rho_{0}\in L^{\infty}(\T^{N})\;\;\mbox{and}\;\;\sqrt{\rho_{0}}u_{0}\in L^{2}(\T^{N}).$$
A solution $(\rho_{1},u_{1})$ is unique on $[0,T_{0}]$ in the set of solutions$(\rho,u)$ such that:
$$\sqrt{\rho}u\in L^{\infty}(L^{2}),\;\n u\in L^{2}_{T_{0}}(L^{2})\;\;\mbox{and}\;\;\rho\in L^{\infty}_{T_{0}}(L^{\infty}),$$
provided that:
\begin{itemize}
\item If $N=2$:
$$
\begin{aligned}
&\n\rho_{1}\in L^{\infty}_{T_{0}}(L^{p}),\;\n u_{1}\in L^{1}_{T_{0}}(L^{\infty})\;\;\mbox{and}\;\;\sqrt{t}\dot{u_{1}}\in L^{2}_{T_{0}}(L^{p}),
\end{aligned}
$$
with $p>2$.
\item If $N\geq 3$:
$$
\begin{aligned}
&\n\rho_{1}\in L^{\infty}_{T_{0}}(L^{N}),\;\n u_{1}\in L^{1}_{T_{0}}(L^{\infty})\;\;\mbox{and}\;\;\sqrt{t}\dot{u_{1}}\in L^{2}_{T_{0}}(L^{N}).
\end{aligned}
$$
\end{itemize}
\end{theorem}
As we know from section \ref{section62} and \ref{section61} that $\n\rho_{1}$ is in $L^{\infty}_{T_{0}}(B^{\e}_{N,1})\hookrightarrow L_{T_{0}}^{\infty}(L^{N})$  and that $\n u_{1}$ is in $L^{1}_{T_{0}}(L^{\infty})$; it suffices to prove that $\sqrt{t} \dot{u}_{1}$ belongs to $L^{2}_{T_{0}}(L^{N})$.\\
We recall then that by Gagliardo-Nirenberg inequalities we have:
$$\sqrt{t}\|\dot{u}_{1}\|_{L^{3}}\leq(t^{\frac{1}{4}-\frac{\e}{2}}\|\dot{u}_{1}-\widetilde{\dot{u}_{1}}\|_{L^{2}})^{\frac{1}{2}}(t^{\frac{3}{2}-\frac{\e}{2}}
\|\n \dot{u}_{1}\|_{L^{2}})^{\frac{1}{2}}t^{\frac{\e}{2}}.$$
From the inequalities (\ref{b1.221}), we deduce that
$$(t^{\frac{1}{4}-\frac{\e}{2}}\|\dot{u}_{1}\|_{L^{2}})^{\frac{1}{2}}\in L^{4}_{t_{0}}(L^{4})\quad\hbox{and}\quad
(t^{\frac{3}{4}-\frac{\e}{2}}\|\n \dot{u}_{1}\|_{L^{2}})^{\frac{1}{2}}\in L^{4}_{T_{0}}(L^{4})$$
 which means that $\sqrt{t}\dot{u}_{1}\in L^{2}_{T_{0}}(L^{3})$.\\
We have then proved then that by using the theorem \ref{theoPG} of \cite{PG}, we have $(\rho,u)=(\rho_{1},u_{1})$  on $[0,T_{0}]$ for all $T_{0}>0$ which conclude the proof.
\subsection{Uniqueness for theorem \ref{corollaire1}  when $P(\rho)=K\rho$ with $K>0$}
In this case, we do not need of any condition of type $\rho_{0}\in B^{1}_{N,1}$, indeed in this specific case we would like to use the results of D. Hoff in \cite{Hoffuni}. To see how to proceed in this case we refer to \cite{H2}.
\subsection{Condition of blow-up for theorem \ref{corollaire1} }
We want to show that by assuming only that $\rho$ is in $L^{\infty}_{T_{0}}$ then we can extend the strong solutions constructed in theorem \ref{corollaire1}. In fact we have proved in section \ref{section4} that the regularity properties (\ref{1.21}), (\ref{b1.22}), (\ref{b1.23}) and (\ref{c1.23}) hold as long as:
\begin{equation}
\sup_{t\in[0,T_{0}]}\|\rho\|_{L^{\infty}_{t}(L^{\infty}(\T^{N}))}<+\infty.
\label{11.23}
\end{equation}
Furthermore using the regularity properties (\ref{1.21}), (\ref{b1.22}), (\ref{b1.23}), (\ref{c1.23}) we have shown in section \ref{section4} that the effective velocity verifies  $\n v_{1}\in L^{1}_{T_{0}}(B^{1+\e}_{N,1})$. In subsection \ref{section62} combining the facts that $\n v_{1}$ and $\rho$ belong respectively to $L^{1}_{T_{0}}(B^{1+\e}_{N,1})$ and  $L^{\infty}_{T_{0}}(\T^{N})$, we have shown that $\n u\in  L^{1}_{T_{0}}(B^{\e}_{\infty,\infty})$. In particular we have $\n u\in  L^{1}_{T_{0}}(L^{\infty})$. This lat condition is then classical to prove that we can extend the strong solutions of theorem \ref{corollaire1}. For more details we refer to \cite{DW} or \cite{H1}.
\section{Proof of theorem \ref{theo3}}
\label{section7}
\subsection{How to obtain a regularizing effect on $v_{1}$ when $\rho\in L^{\infty}(L^{q})$}
We now want to work with the \textit{effective velocity} $v_{1}$ introduced in the previous sections to obtain new estimates this last when we assume that  $\rho$ belongs to $L^{\infty}(L^{q})$. We will give more details on the value of $q$ in the sequel of the proof.  We can now rewrite the momentum equation of system (\ref{1}). We obtain then the following equation where we have set $\nu=2\mu+\lambda$:
$$\rho\p_{t}u+\rho u\cdot \n u-\mu\D\big(u-\frac{1}{\nu}v\big)-(\lambda+\mu)\n{\rm div}\big(u-\frac{1}{\nu}v\big)=\rho g,$$
where we recall that $v=(\D)^{-1}(\n P(\rho))$ with $(\D)^{-1}$ the inverse Laplacian with zero mean value on $\T^{N}$.
As $v_{1}=u-\frac{1}{\nu}v$  we have:
$$\rho\p_{t}v_{1}+\rho u\cdot \n u-\mu\D v_{1}-(\lambda+\mu)\n{\rm div}v_{1}=\rho g-\frac{1}{\nu}\rho\p_{t}v.$$
As ${\rm div}v=P(\rho)-\int_{\T^{N}}P(\rho)dx$, from the transport equation we obtain:
$$
\begin{aligned}
{\rm div}\p_{t}v&=-P^{'}(\rho)\rho{\rm div}u-\n P(\rho)\cdot u+\widetilde{P^{'}(\rho)\rho{\rm div}u}+\widetilde{\n P(\rho)\cdot u}\\
&=-{\rm div}(P(\rho)u)+(P(\rho)-\rho P^{'}(\rho)){\rm div}u-\widetilde{P(\rho){\rm div}u}+\widetilde{\rho P^{'}(\rho)){\rm div}u}.
\end{aligned}
$$
In the sequel we will need to use the Bogovskii operator that we note $\Lambda^{-1}$ (see \cite{NS} p168 for a definition), we obtain then:
\begin{equation}
 \p_{t}v=\Lambda^{-1}\big(-{\rm div}(P(\rho)u)+(P(\rho)-\rho P^{'}(\rho)){\rm div}u-\widetilde{P(\rho){\rm div}u}+\widetilde{\rho P^{'}(\rho){\rm div}u}\big).
\label{Bogo}
\end{equation}
We get finally:
\begin{equation}
\rho\p_{t}v_{1}-\mu\D v_{1}-(\lambda+\mu)\n{\rm div}v_{1}=\rho g-\rho u\cdot \n u-\frac{1}{\nu}\rho\p_{t}v.
 \label{hchaleur}
\end{equation}
We set $f(t)=\min(t,1)$ and we remark that $f(0)=0$. We multiply then (\ref{hchaleur}) by $f(t)\p_{t}v_{1}$ and integrating on $(0,t)\times\T^{N}$(with $0<t<T$) we obtain then:
\begin{equation}
\begin{aligned}
&\int^{t}_{0}\int_{\T^{N}}f(s)\rho|\p_{s}v_{1}|^{2}dxds+\frac{1}{2}\int_{\T^{N}}f(t)\big(\mu|\n v_{1}(t,x)|^{2}+\xi({\rm div} v_{1})^{2}(t,x)\big)dx\leq\\
&\hspace{3cm}\int^{t}_{0}\int_{\T^{N}}f^{'}(s)(\mu|\n v_{1}(t,x)|^{2}+\xi({\rm div} v_{1})^{2}(t,x)\big)dxds\\
&\hspace{1cm}+\int^{t}_{0}\int_{\T^{N}}
\rho u\cdot\n u \, f(s)\p_{t}v_{1}dxds+\int^{t}_{0}\int_{\T^{N}}(\rho g-\frac{1}{\nu}\rho\p_{s}v)f(s)\p_{s}v_{1}dxds,
\end{aligned}
 \label{h74}
\end{equation}
where $\xi=\mu+\lambda$. We have then as for all $t\in(0,T)$:
\begin{equation}
\|P(\rho)\|_{L^{2}_{t}(L^{2})}\leq C,
\label{45imp1}
\end{equation}
and by Young's inequality:
\begin{equation}
\begin{aligned}
&\int^{t}_{0}\int_{\T^{N}}f(s)\rho|\p_{s}v_{1}|^{2}dxds+\frac{1}{2}\int_{\T^{N}}f(t)\big(\mu|\n v_{1}(t,x)|^{2}+\xi({\rm div} v_{1})^{2}(t,x)\big)dx\leq\\
&\hspace{2cm}C\big(1+\int^{t}_{0}\int_{\T^{N}}
f(s)\rho |u\cdot\n v_{1}|^{2}dxds+\int^{t}_{0}\int_{\T^{N}}
f(s)\rho |u\cdot\n v|^{2}dxds\\
&\hspace{7cm}+\int^{t}_{0}\int_{\T^{N}}f(s)\rho (|g|^{2}+|\p_{s}v|^{2})dxds\big).
\end{aligned}
 \label{h741}
\end{equation}
We have next to control the terms on the right hand side of (\ref{h741}). In the sequel for simplicity, we will  treat only the case $N=3$. The case $N=2$ follows exactly similar lines. We can now recall that from the works of A. Mellet and A. Vasseur in \cite{5MV2}, we are able to control the velocity $u$ in $L^{\infty}(L^{\infty})$ if we suppose that $\rho$ is in $L^{\infty}(L^{3\gamma+\e})$ with $\e>0$. In fact from the inequality (\ref{againvitesse}), we can obtain a gain of integrability on the velocity, i.e $\rho^{\frac{1}{p}}u\in L^{\infty}(L^{p})$ with $p$ arbitrary big if $P(\rho)\in L^{p}(L^{\frac{3p}{p+1}})$. In our case it will be the case as we assume at least that $P(\rho)$ is in $ L^{\infty}(L^{3})$, i.e:
\begin{equation}
\|P(\rho)\|_{L_{t}^{\infty}(L^{3})}\leq C.
\label{45imp11}
\end{equation}
For the simplicity of the calculus we will assume that $u\in L_{t}^{\infty}(L^{\infty})$.
In fact we have only a control on  $\rho^{\frac{1}{p}}u$ in $L^{\infty}(L^{p})$ for $p$ arbitrarily large, but in the sequel all the expressions to treat will be of the form $\rho u$. It would suffice to apply the H\"older's inequalities with $\rho^{1-\frac{1}{p}}(\rho^{\frac{1}{p}}u)$. 
\subsubsection*{Regularizing effect on $\D v_{1}$}
We want here to use the regularizing effect on $v_{1}$ and proceed in the sequel by bootstrap. To do it we use the momentum equation (\ref{hchaleur}) and we have:
\begin{equation}
\begin{aligned}
&\mu\D v_{1}+(\lambda+\mu)\n{\rm div}v_{1}=\rho\p_{t}v_{1}+\rho u\cdot \n v_{1}+\rho u\cdot\n(\D)^{-1}\n( P(\rho))\\
&\hspace{9cm}-\rho g+\frac{\rho}{\nu}\p_{t}v.
\end{aligned}
 \label{bchaleur1}
\end{equation}
We now want to take in consideration  the ellipticity of (\ref{bchaleur1}), in this goal we would like to recall that from (\ref{h741}) we can hope only a control of $\sqrt{\rho f(t)}\p_{t}v_{1}$ in $L^{2}_{t}(L^{2})$. We set then $\frac{1}{p}=\frac{1}{2}+\frac{1}{2q}$ and we have:
\begin{equation}
\begin{aligned}
&\|\D v_{1}\|_{L^{p}}\leq \|\rho\|_{L^{q}}^{\frac{1}{2}}\|\sqrt{\rho}\p_{t}v_{1}\|_{L^{2}}+\|\rho\|_{L^{q}}^{\frac{1}{2}}\|\sqrt{\rho}\p_{t}v_{1}\|_{L^{2}}+\|\rho\|_{L^{q}}\|u\|_{L^{\infty}}
\|\rho\|_{L^{q}}^{\gamma}\\
&\hspace{8cm}+\|\rho\|_{L^{q}}(\|\p_{t}v\|_{L^{q_{1}}}+\|g\|_{L^{q_{1}}}),
\end{aligned}
\label{effetreg}
\end{equation}
with the following conditions: $\frac{1}{q}+\frac{1}{q_{1}}\leq\frac{1}{p}$, $\frac{1}{q_{1}}=\frac{1}{2}-\frac{1}{2q}$ and $\frac{1}{q}+\frac{\gamma}{q}\leq\frac{1}{p}$ (let $q\geq 2\gamma+1$).
Next by Gagliardo-Nirenberg, we have
\begin{equation}
\|\n v_{1}\|_{L^{q_{2}}}\leq \|\n v_{1}\|_{L^{2}}^{\e}\|\D v_{1}\|_{L^{p}}^{1-\e},
\label{Gagli1}
\end{equation}
with: $\frac{1}{q_{2}}=\frac{\e}{2}+\frac{1-\e}{p}-\frac{1-e}{3}=\frac{1}{2}+\frac{1}{2q}-\frac{1}{3}+\frac{\e}{3}-\frac{\e}{2q}=\frac{1}{6}+\frac{1}{2q}(1-\e)$.  We can now estimate the following term $\int^{t}_{0}\|\sqrt{f(s)\rho}u\cdot\n v_{1}\|_{L^{2}}^{2}ds$.
\subsubsection*{Estimate on the term $\int^{t}_{0}\|\sqrt{f(s)\rho}u\cdot\n v_{1}\|_{L^{2}}^{2}ds$}
We have then by using (\ref{Gagli1}) where $\e$ is defined such that:
\begin{equation}
\frac{1}{q_{1}}\geq \frac{1}{2}+\frac{1}{2q}-\frac{1}{3}+\frac{\e}{3}   \Leftrightarrow  \frac{1}{q}-\frac{\e}{2q}\leq\frac{1}{3}  \Leftrightarrow q>3,
\label{45imp3}
\end{equation}
and with (\ref{effetreg}) we have:
$$
\begin{aligned}
&\|\sqrt{f(s)\rho}\;u\cdot\n v_{1}\|^{2}_{L^{2}}\leq \|\rho\|_{L^{q}}\|u\|^{2}_{L^{\infty}}\|\sqrt{f(s)} \n v_{1}\|^{2}_{L^{q_{1}}},\\
&\hspace{3,3cm}\leq \|\rho\|_{L^{q}}\|u\|^{2}_{L^{\infty}}f(s)\|\n v_{1}\|_{L^{2}}^{2\e}\|\D v_{1}\|_{L^{p}}^{2(1-\e)},\\
&\leq C\|\rho\|_{L^{q}}\|u\|^{2}_{L^{\infty}}f(s)\|\n v_{1}\|_{L^{2}}^{2\e}\,\big(\|\rho\|_{L^{q}}^{\frac{1}{2}}\|\sqrt{\rho}\p_{t}v_{1}\|_{L^{2}}+|\rho\|_{L^{q}}^{\frac{1}{2}}\|\sqrt{\rho}\; u\cdot\n v_{1}\|_{L^{2}}
\\
&\hspace{3cm}+\|\rho\|_{L^{q}}\|u\|_{L^{\infty}}
\|\rho\|_{L^{q}}^{\gamma}+\|\rho\|_{L^{q}}(\|\p_{t}v\|_{L^{q_{1}}}+\|g\|_{L^{q_{1}}})\big)^{2(1-\e)},\\
\end{aligned}
$$
We have then:
$$
\begin{aligned}
&\|\sqrt{f(s)\rho}\;u\cdot\n v_{1}\|^{2}_{L^{2}}\leq \\
&C\|u\|^{2}_{L^{\infty}}f(s)^{\e}\|\n v_{1}\|_{L^{2}}^{2\e}\,\|\rho\|_{L^{q}}^{2-\e}\big(\|\sqrt{f(s)\rho}\;\p_{t}v_{1}\|_{L^{2}}+
\|\sqrt{f(s)\rho}\; u\cdot\n v_{1}\|_{L^{2}}\big)^{2(1-\e)}\\
&\hspace{0,5cm}+C\|\rho\|_{L^{q}}\|u\|^{2}_{L^{\infty}}f(s)\|\n v_{1}\|_{L^{2}}^{2\e}(\|\rho\|_{L^{q}}\|u\|_{L^{\infty}}
\|\rho\|_{L^{q}}^{\gamma}+\|\rho\|_{L^{q}}(\|\p_{t}v\|_{L^{q_{1}}}+\|g\|_{L^{q_{1}}})\big)^{2(1-\e)},\\
\end{aligned}
$$
Next by Young inequality with $\frac{1}{1-\e}$ and $\frac{1}{\e}$, we get:
\begin{equation}
\begin{aligned}
&\|\sqrt{f(s)\rho}u\cdot\n v_{1}\|^{2}_{L^{2}}\leq C_{\alpha}f(s)\|\rho\|^{\frac{2}{\e}-1}_{L^{q}}\|u\|^{\frac{2}{\e}}_{L^{\infty}}\|\n v_{1}\|_{L^{2}}^{2}+\alpha\|\sqrt{f(s)\rho}\p_{t}v_{1}\|^{2}_{L^{2}}\\
&+\alpha\|\sqrt{f(s)\rho}\; u\cdot\n v_{1}\|_{L^{2}}^{2}+C_{\alpha}f(s)\|\rho\|^{\frac{1}{\e}}_{L^{q}}\|u\|^{\frac{2}{\e}}_{L^{\infty}}\|\n v_{1}\|_{L^{2}}^{2}+\alpha\|\rho\|^{2}_{L^{q}}\|u\|^{2}_{L^{\infty}}
\|\rho\|_{L^{q}}^{2\gamma}\\
&\hspace{4cm}+\alpha\|\rho\|^{2}_{L^{q}}(\|\sqrt{f(s)}\p_{t}v\|^{2}_{L^{q_{1}}}+
\|g\|^{2}_{L^{q_{1}}}),\\
\end{aligned}
\label{superimp}
\end{equation}
Here $\alpha$ is very small and $C_{\alpha}$ can be very big.
From (\ref{superimp}) we obtain:
\begin{equation}
\begin{aligned}
&(1-\alpha)\|\sqrt{f(s)\rho}u\cdot\n v_{1}\|^{2}_{L^{2}}\leq C_{\alpha}f(s)\|\rho\|^{\frac{2}{\e}-1}_{L^{q}}\|u\|^{\frac{2}{\e}}_{L^{\infty}}\|\n v_{1}\|_{L^{2}}^{2}\\
&+\alpha\|\sqrt{f(s)\rho}\p_{t}v_{1}\|^{2}_{L^{2}}+C_{\alpha}f(s)\|\rho\|^{\frac{1}{\e}}_{L^{q}}\|u\|^{\frac{2}{\e}}_{L^{\infty}}\|\n v_{1}\|_{L^{2}}^{2}+\alpha\|\rho\|^{2}_{L^{q}}\|u\|^{2}_{L^{\infty}}
\|\rho\|_{L^{q}}^{2\gamma}\\
&\hspace{6cm}+\alpha\|\rho\|^{2}_{L^{q}}(\|\sqrt{f(s)}\p_{t}v\|^{2}_{L^{q_{1}}}+
\|g\|^{2}_{L^{q_{1}}}),\\
\end{aligned}
\label{superimp1}
\end{equation}
We now inject inequality (\ref{superimp1}) in ( \ref{h741}) and we obtain by choosing $\alpha$ enough small for doimg a bootstrap:
\begin{equation}
\begin{aligned}
&\int^{t}_{0}\int_{\T^{N}}f(s)\rho|\p_{s}v_{1}|^{2}dxds+\frac{1}{2}\int_{\T^{N}}f(t)\big(\mu|\n v_{1}(t,x)|^{2}+\xi({\rm div} v_{1})^{2}(t,x)\big)dx\leq\\
&\hspace{2cm}C\biggl(1+(\|\rho\|^{\frac{2}{\e}-1}_{L^{\infty}_{t}(L^{q})}+\|\rho\|^{\frac{1}{\e}}_{L^{\infty}_{t}(L^{q})})
\|u\|^{\frac{2}{\e}}_{L^{\infty}_{t}(L^{\infty})}\int^{t}_{0}f(s)\|\n v_{1}(s,\cdot)\|_{L^{2}}^{2}ds\\
&+t\|\rho\|^{2+2\gamma}_{L^{\infty}_{t}(L^{q})}\|u\|^{2}_{L^{\infty}_{t}(L^{\infty})}
+\|\rho\|^{2}_{L^{\infty}_{t}(L^{q})}\int^{t}_{0}(\|\sqrt{f(s)}\p_{t}v(s,\cdot)\|^{2}_{L^{q_{1}}}+
\|g(s,\cdot)\|^{2}_{L^{q_{1}}})ds\\
&\hspace{2cm}+\int^{t}_{0}\int_{\T^{N}}
f(s)\rho |u\cdot\n v|^{2}dxds+\int^{t}_{0}\int_{\T^{N}}f(s)\rho (|g|^{2}+|\p_{s}v|^{2})dxds\big).
\end{aligned}
 \label{h741a2}
\end{equation}
It reminds to bound the last term on the right hand side of (\ref{h741a2}). We need in particular to prove that $\sqrt{f(s)}\p_{t}v\in L^{2}_{t}(L^{q_{1}})$.\\
Now from  (\ref{Bogo}), we have:
$$
\begin{aligned}
&\int^{t}_{0}\|\sqrt{f(s)}\p_{t}v(s,\cdot)\|^{2}_{L^{q_{1}}}ds\leq C(1+ \int^{t}_{0}\|\sqrt{f(s)}P(\rho)(s,\cdot)u(s,\cdot)\|^{2}_{L^{q_{1}}}ds\\
&\hspace{5cm}+\int^{t}_{0}\|\sqrt{f(s)}\Lambda^{-1}\big((P(\rho)-\rho P^{'}(\rho)){\rm div}u\big)(s,\cdot)\|^{2}_{L^{q_{1}}}ds.
\end{aligned}$$
We have then:
$$\int^{t}_{0}\|\sqrt{f(s)}P(\rho)(s,\cdot)u(s,\cdot)\|^{2}_{L^{q_{1}}}ds\leq C\|u\|^{2}_{L^{\infty}_{t}(L^{\infty})}\|\rho\|_{L^{2\gamma}(L^{\gamma q_{1}})}.$$
It means that we need that:
\begin{equation}
q\geq \gamma q_{1}  \Leftrightarrow q\geq \frac{2q\gamma}{q-1}  \Leftrightarrow q(q-(2\gamma+1))\geq 0   \Leftrightarrow q\geq 2\gamma+1.
\label{45imp5}
\end{equation}
Next we have as $(P(\rho)-\rho P^{'}(\rho)){\rm div}u$ belongs to $L^{2}(L^{p_{2}})$ with $p_{2}=\frac{2q}{q+2\gamma}$ ($\frac{1}{p_{2}}=\frac{1}{2}+\frac{\gamma}{q}$) and the properties on the Bosvskii operator we have
\begin{equation}
\int^{t}_{0}\|\sqrt{f(s)}\Lambda^{-1}\big(P(\rho)-\rho P^{'}(\rho)){\rm div}u\big)\|^{2}_{L^{q_{1}}}\leq
C\|\rho\|_{L^{\infty}_{t}(L^{q})}^{2\gamma}\|{\rm div}u\|^{2}_{L^{2}_{t}(L^{2})},
\label{imp44x}
\end{equation}
where we need to assume that
\begin{equation}
\frac{\gamma}{q}+\frac{1}{2}-\frac{1}{3}\leq\frac{1}{q_{1}}  \Leftrightarrow q\geq 3\gamma+\frac{3}{2}.
\label{45imp6}
\end{equation}
Now from (\ref{h741a2}) and the previous inequalities we have:
\begin{equation}
\begin{aligned}
&\int^{t}_{0}\int_{\T^{N}}f(s)\rho|\p_{s}v_{1}|^{2}dxds+\frac{1}{2}\int_{\T^{N}}f(t)\big(\mu|\n v_{1}(t,x)|^{2}+\xi({\rm div} v_{1})^{2}(t,x)\big)dx\leq\\
&\hspace{6cm}C\big(1+\int^{t}_{0}\int_{\T^{N}}f(s)\rho (|g|^{2}+|\p_{s}v|^{2})dxds\big).
\end{aligned}
 \label{hh7412}
\end{equation}
By proceeding as in the previous terms we can easily bound the last terms on the right hand side of (\ref{hh7412}).  We conclude finally that if $q$ check (\ref{45imp1}), (\ref{45imp11}), (\ref{45imp3}), (\ref{45imp5}) and (\ref{45imp6}), we have then :
\begin{equation}
\begin{aligned}
&\int^{t}_{0}\int_{\T^{N}}f(s)\rho|\p_{s}v_{1}|^{2}dxds+\frac{1}{2}\int_{\T^{N}}f(t)\big(\mu|\n v_{1}(t,x)|^{2}+\xi({\rm div} v_{1})^{2}(t,x)\big)dx\leq C_{t},
\end{aligned}
\end{equation}
where $C_{t}$ depends on $t$.
\subsection*{Control on the norm $\rho\in L^{\infty}$}
Now we can come back to equation (\ref{g28}) in order to get a control in norm $L^{\infty}$ on the density. More precisely we have:
\begin{equation}
\begin{aligned}
&\log(\rho(t,x))\leq\log(\|\rho_{0}\|_{L^{\infty}})+C\|(\D)^{-1}{\rm div} m_{0}\|_{L^{\infty}}+C\|(\D)^{-1}{\rm div}(\rho u)\|_{L^{\infty}}\\
&+C\int^{t}_{0}\int_{\T^{N}}P(\rho(s,x)ds\,dx+C \int^{t}_{0}\frac{1}{\sqrt{f(s)}}\|[\sqrt{f(s)}u_{j},R_{i}R_{j}](\rho u_{i})(s,\cdot)\|_{L^{\infty}}ds.
\end{aligned}
 \label{28}
\end{equation}
Easily, we have by Sobolev embedding and (\ref{againvitesse}) with $\e>0$:
$$
\begin{aligned}
\|(\D)^{-1}{\rm div}(\rho u)\|_{L^{\infty}}&\leq\|\rho u\|_{L^{\infty}(L^{3+\e})}\\
&\leq\|\rho\|_{L^{\infty}(L^{q})}\|u\|_{L^{\infty}(L^{\infty})},
\end{aligned}
$$
when we assume that:
\begin{equation}
\|\rho\|_{L^{\infty}_{t}(L^{q})}\leq C\;\;\;\mbox{with}\;\;q>3.
\label{45imp7}
\end{equation}
We recall then from the previous section that $\sqrt{f(s)}\D v_{1}\in L^{2}(L^{p})$ with $\frac{1}{p}=\frac{1}{2}+\frac{1}{2q}$, by Gagliardo-Nirenberg inequality we have:
$$\|\n v_{1}\|_{L^{q_{2}}}\leq \|\D v_{1}\|^{1-\e}_{L^{p}}\|\n v_{1}\|^{\e}_{L^{2}},$$
with $\frac{1}{q_{2}}=\frac{1}{2}+\frac{1}{2q}-\frac{1}{3}+\frac{\e}{3}-\frac{\e}{2q}$. We have then
by using the results of R. Coifman et al in \cite{1}:
\begin{equation}
\begin{aligned}
&(f(s))^{\frac{1}{2}-\frac{\e}{2}}\|[(v_{1})_{j},R_{i}R_{j}](\rho u_{i})(s,\cdot)\|_{W^{1,\alpha}}\leq \|\n v_{1}\|_{L^{2}}^{\e}\|\sqrt{f(s)} \D v_{1}\|_{L^{p}}^{1-\e}\\
&\hspace{9cm}\times\|\rho\|_{L^{q}}\|u\|_{L^{\infty}}.
\end{aligned}
\label{infinal1}
\end{equation}
In the sequel we will need that $\alpha>3$. Indeed by Sobolev embedding we will prove that $[(v_{1})_{j},R_{i}R_{j}](\rho u_{i})(s,\cdot)$ is in $L^{1}_{t}(L^{\infty})$. That is why in the sequel we have to assume that when $\e$ is chosen arbitrary small:
\begin{equation}
\frac{1}{q}+\frac{1}{q_{2}}=\frac{1}{2}+\frac{3}{2q}-\frac{1}{3}+\frac{\e}{3}-\frac{\e}{2q}<\frac{1}{3}\Leftrightarrow q>9
\label{45imp8}
\end{equation}
We have finally by using hypothesis (\ref{45imp8}) and  (\ref{infinal1}) as  $\frac{1}{f(s)^{\frac{1}{2}-\frac{\e}{2}}}\in L^{2}_{t}$ and as we have shown that$(f(s))^{\frac{1}{2}-\frac{\e}{2}}\|[(v_{1})_{j},R_{i}R_{j}](\rho u_{i})(s,\cdot)\|_{W^{1,\alpha}}$ in $L^{2}_{t}(L^{\infty})$ by Sobolev embedding (because here $\alpha>3$) :
$$
\begin{aligned}
&|\int^{t}_{0}\frac{1}{f(s)^{\frac{1}{2}-\frac{\e}{2}}}f(s)^{\frac{1}{2}-\frac{\e}{2}}\|[(v_{1})_{j},R_{i}R_{j}](\rho u_{i})(s,\cdot)\|_{L^{\infty}}ds|\leq C \|\rho\|_{L^{\infty}_{t}(L^{q})}\|u\|_{L^{\infty}_{t}(L^{\infty})}\\
&\hspace{8cm}\|\n v_{1}\|_{L^{2}(L^{2})}^{\e}\|\sqrt{f(s)} \D v_{1}\|_{L^{2}(L^{p})}^{1-\e}.
\end{aligned}
$$
We proceed similarly for the term $\|[(v_{j},R_{i}R_{j}](\rho u_{i})(s,\cdot)\|_{L^{1}_{s}(L^{\infty})}$. This term is the most important because it decides of the value of $q$ that we must choose. Indeed we have the results of R. Coifman et al in \cite{1}:
$$\|[(\D)^{-1}\n(P(\rho))_{j},R_{i}R_{j}](\rho u_{i})(s,\cdot)\|_{W^{1,\beta}}\leq\|\rho\|_{L^{q}}^{\gamma}\|\|u\|_{L^{\infty}}
\|\rho\|_{L^{q}}.$$
We need that $\beta>3$ for using Sobolev embedding and then to prove by this way  that $[(\D)^{-1}\n(P(\rho))_{j},R_{i}R_{j}](\rho u_{i})(s,\cdot)$ is in $L^{1}_{t}(L^{\infty})$. We assume then that:
\begin{equation}
\frac{1}{q}+\frac{\gamma}{q}<\frac{1}{3}\Leftrightarrow q>3(\gamma+1).
\label{45imp9}
\end{equation}
If we summarize all the inequalities on $q$, i.e (\ref{45imp1}), (\ref{45imp11}), (\ref{45imp3}), (\ref{45imp5}), (\ref{45imp6}), (\ref{45imp7}), (\ref{45imp8}) and (\ref{45imp9}), we need that:
\begin{equation}
\begin{aligned}
\rho\in L^{\infty}_{t}(L^{3\gamma}(\T^{N}))\cap L^{\infty}_{t}(L^{3\gamma+\frac{3}{2}}(\T^{N}))\cap L^{\infty}(L^{9+\e}(\T^{N}))\cap L^{\gamma+1}_{t}(L^{3\gamma+3}(\T^{N}))
\end{aligned}
\label{final1111}
\end{equation}
Finally under theses conditions we control
$|\log\rho|1_{\{\rho\geq1\}}\in L^{\infty}_{t}(L^{\infty}(\T^{N})$. From theorem \ref{corollaire1}, we have seen that we can control $\n u$ in $L^{1}_{t}(L^{\infty})$ and that we can extend beyond $T$ the solutions constructed in theorem \ref{corollaire1}. It achieves the proof of theorem \ref{theo3}.
\null{\hfill $\Box$}
\section{Further comments, results and open problems}
\label{section8}
We now want to describe the additional problems when we consider variable viscosity coefficients. In particular we will mention the specific case of shallow-water system. We will point out also that we can get in some specific cases on the choice of the variable viscosity coefficients (which includes in particular the case of the shallow-water system, see \cite{5BD1}) a
gain of integrability on the pressure. More precisely we assume that the viscosity coefficients verify:
\begin{equation}
\lambda(\rho)=2(\rho\mu^{'}(\rho)-\mu(\rho)). \label{5coeff}
\end{equation}
In particular in this we are able to verify the condition (\ref{egcrucial}). Unfortunately we will explain why it seems difficult to apply the theorem \ref{theo3} with this choice on the capillarity coefficients.
\subsection*{When the viscosity coefficients are variable}
We  briefly want to remind some results of global weak solution when the
viscosity coefficients are variables. Furthermore we would explain why it seems complicated to obtain similar result than theorem \ref{theo3} in this case. Indeed one of the main reason is the loss of the so called structure of \textit{effective pressure} or of \textit{effective velocity}. We will give more details on this in the sequel.\\
In \cite{5BD1} Bresch and Desjardins showed  a result of global
stability of weak solutions for the non isothermal Navier-Stokes
system assuming density dependence on the viscosity coefficients $\mu$ and $\lambda$, considering perfect gas law with some cold pressure close to the vacuum, and the algebraic relation (\ref{5coeff}).\\
The key point in this paper is to
show that the structure of the diffusion term provides some regularity for the
density thanks to a new mathematical entropy inequality. This one has been discovered in \cite{5BD2}, we call it the BD entropy. More precisely they are able to obtain a gain of derivative on the density, i.e $\sqrt{\rho}\n\va(\rho)\in L^{\infty}(L^{2})$ when $\sqrt{\rho_{0}}\n\va(\rho_{0})\in L^{2}$.
Here we have set  $\va^{'}(\rho)=\frac{2\mu^{'}(\rho)}{\rho}$.\\
Mellet and Vasseur by using the BD entropy, get in \cite{5MV} a very interesting new stability result when the pressure is considered barotropic, i.e $P(\rho)=a\rho^{\gamma}$ with $\gamma\geq 1$ and $a>0$. Note that the main difficulty is to establish the compactness of
$\sqrt{\rho}u$ in $L^{2}$ strong, and the key ingredient to achieve this is an additional
estimate which bounds $\rho|u|^2$ in a space better than $L^{\infty}(0, T ; L^{1}(\T^{N}))$. Indeed
the BD viscosity coefficients vanish on the vacuum set, so it means that we loss the control $\n u$ in $L^{2}(L^{2}(\T^{N}))$, that is why it is not so clear to obtain strong convergence on the terms of type $\rho u\otimes u$.\\
Unfortunately, the construction of approximate solutions satisfying: energy estimates, BD mathematical entropy and Mellet-Vasseur estimates is far from being proven except in dimension one or with symmetry assumptions, see \cite{aMV1}, \cite{aJJX}, \cite{aGJX}. Note that approximate solutions construction process has been proposed in \cite{aBD} satisfying energy estimates and BD mathematical entropy. This leads to global existence of weak solutions,
however only if we assume that extra terms or cold pressure are present.\\
\subsection{The BD entropy and the theorem \ref{theo3}}
In theorem \ref{theo3}, we prove a blow-up criterion for strong solution $(\rho,u)$ on $(0,T)$ who say us that if in dimension three:
\begin{equation}
\rho\in L^{\gamma+1}_{T}(L^{(N+1+\e)\gamma}(\T^{N}))\;\;\mbox{and}\;\;\rho\in L^{\infty}_{T}( L^{9+\e}(\T^{N})\cap L^{3\gamma+\frac{3}{2}}(\T^{N})),
\label{egcrucial1}
\end{equation}
then the solution $(\rho,u)$ can extend beyond $T$.
 We now can motivate our assumption (\ref{egcrucial1}) in the light of the BD entropy (see \cite{5BD}). Indeed by choosing some viscosity coefficients verifying the equality (\ref{5coeff}), we are able to control $\sqrt{\rho}\n\va(\rho)$ in $L^{\infty}(L^{2})$. In particular we can obtain the relation (\ref{egcrucial1}) by Sobolev embedding and by choosing $\mu(\rho)=\mu\rho^{\alpha}$ with $\alpha$ big enough. It means that if we would be able to extend theorem \ref{theo3} to the case of the BD viscosity coefficients, we could prove the existence of global strong solution in dimension $N=3$ for compressible Navier-Stokes with this choice of viscosity coefficient. However this type of viscosity coefficient \textit{kills} the structure of \textit{effective pressure} or of \textit{effective velocity} and we can not apply our proof. In particular it appears not so clear how to obtain estimates on $\p_{t}u$ by multiplying the momentum equation by $\p_{t}u$.\\
 It shows that the structure of the viscosity coefficients plays a crucial role for compressible Navier-Stokes system and that the structure changes completely and depends crucially on the choice of the viscosity coefficients.\\
 An interesting open problem would be to extend theorem \ref{theo3} to the case of variable viscosity coefficients and in particular BD coefficients.
 \subsection*{Importance of the regularity conditions on the source term $g$}
  V. A. Waigant has built in \cite{CCK25} explicit solutions for which the maximal integrability of the density corresponds to $L^{q}(0,1,L^{q})$ with $q=\frac{\gamma(3N+2)-N}{2N}$. It means that
(\ref{egcrucial1}) fails in this case except that the force term $g$ introduced in \cite{CCK25} is less regular than what the theorem \ref{theo3} requires. It means that the regularity of $g$ is crucial to get strong solutions and it shows in particular that it is quite hard to obtain a gain of integrability on the pressure. In passing we recall that it is possible to prove that we have a gain of integrability on the density in $L^{q}_{T}(L^{q})$ with $q=\gamma+\frac{2}{N}\gamma-1$. The paper of  V. A. Waigant was a 	counter-example to prove that in dimension $N=2,3$ it was not possible to hope a control $L^{2}_{loc}$ on the density. It indicates the necessity to use other arguments than the theory of renormalized solutions to obtain global weak solution when $N=2,3$.
\subsection*{On some extension of the theorem \ref{theo3} and some questions of  scaling}
In theorem \ref{theo3} we need to assume that $\lambda=0$ to get a control $L^{\infty}$ on the velocity  $u$ as in the article of A. Mellet and A. Vasseur in \cite{5MV2}. We recall that in this paper they need of a control on $P(\rho)\in L^{\infty}(L^{3+\e})$ with $\e>0$ for $N=3$ to obtain a control $L^{\infty}$ on the velocity when the viscosity coefficients are constant. To do this, they use some De Giorgi technics (see also  \cite{V6}  where Vasseur reproves  the so-called result of Caffarelli-Kohn and Nirenberg in \cite{CKN} with this type of arguments). For De Giorgi methods we refer to \cite{Degiorgi}, where De Giorgi proves the so called XIX Hilbert's problem, which can be reduced to the regularity of weak solutions to nonlinear elliptic equations or systems.\\
As in \cite{5MV2}, the pressure plays an important role, and in some sense the pressure is the good unknown to control in order to get global strong solutions. In fact, the pressure plays the role of the velocity for the incompressible Navier-Stokes .\\
Furthermore  we could extend the theorem \ref{theo3} to the case where $\lambda(\rho)$ is non null. However in this case we would need of stronger assumption on the control of $P(\rho)$ in $L^{\infty}(L^{q})$ with $q$ bigger. Indeed as in \cite{5MV2} the gain of integrability on the velocity $u$ in $L^{\infty}(L^{p})$ depends on the ratio between $\mu$ and $\lambda$.\\
We would like to point out that in the theorem \ref{theo3} we need to assume that the initial data $u_{0}$ is in $L^{\infty}$ in order to obtain a minimal condition of blow-up in term of integrability on the density. We just observe that in terms of scaling we ask one derivative in more on the initial velocity, but in return we are "not so far" to obtain a blow-up criterion on the pressure with one derivative in less compared with the scaling of the equations. More precisely the fact as in \cite{5MV2} to ask a control on $P(\rho)$ in $L^{\infty}(L^{3+\e})$ in dimension $N=3$ to control $u$ in $L^{\infty}$ is a condition of type "one derivative in less for the scaling". Indeed we recall that the scaling of the pressure (when $\frac{1}{\rho}$ and $\rho$ are in $L^{\infty}$) is $L_{T}^{\infty}(W^{1,N}(\T^{N})$. In our case, we are "not so far" to get a similar result as we ask that $P(\rho)$ is in $ L^{\infty}(L^{N+1+\e}(\T^{N}))$.\\
\\
\textit{ In some sense we claim that the derivative in more for the scaling on the velocity is transferred in one derivative in less for the blow-up condition on the pressure.}
\section{Appendix}
\label{section9}
This section is devoted to the proof of commutator estimates which have been used in section $2$ and $3$. They are based on
paradifferentiel calculus, a tool introduced by J.-M. Bony in \cite{5BJM}. The basic idea of paradifferential calculus is that
any product of two distributions $u$ and $v$ can be formally decomposed into:
$$uv=T_{u}v+T_{v}u+R(u,v)=T_{u}v+T^{'}_{v}u$$
where the paraproduct operator is defined by $T_{u}v=\sum_{q}S_{q-1}u\D_{q}v$, the remainder operator $R$ by
$R(u,v)=\sum_{q}\D_{q}u(\D_{q-1}v+\D_{q}v+\D_{q+1}v)$ and $T^{'}_{v}u=T_{v}u+R(u,v)$.
We would like to remind a basic lemma on the commutator (for more details see \cite{BCD} p 110).
\begin{lemme}
Let $\theta$ be a $C^{1}$ function on $\T^{N}$ such that $\hat{\theta}\in L^{1}$. There exists a constant $C$ such that for any Lipschitz function $a$ and any function $b$ in $L^{p}$ with $p\in [1,+\infty]$, we have:
$$\forall \lambda>0,\;\;\|[\theta(\lambda^{-1}D),a]b\|_{L^{p}}\leq C\lambda^{-1}\|\n a\|_{L^{\infty}}\|b\|_{L^{p}}.$$
\label{d297}
\end{lemme}
{\bf Proof:} In order to prove this lemma, it suffice to rewrite $\theta(\lambda^{-1}D)$ as a convolution operator. More precisely we have:
$$
\begin{aligned}
&\big([\theta(\lambda^{-1}D),a]b\big)(x)=\theta(\lambda^{-1}D)(ab)(x)-a(x)\theta(\lambda^{-1}D)b(x),\\
&\hspace{2cm}=\lambda^{N}\int_{\T^{N}}h(\lambda(x-y))(a(y)-a(x))b(y)dy\;\;\;\mbox{with}\;\;h={\cal F}^{-1}\theta.
\end{aligned}
$$
As the function $a$ is Lipschitz, we have therefore
$$\big|\big([\theta(\lambda^{-1}D),a]b\big)(x)\big|\leq\lambda^{N}\|\n a\|_{L^{\infty}}\int_{\T^{N}}|h(\lambda(x-y))||y-x|b(y)dy.
$$
By Young's inequality, it implies:
$$\|\big([\theta(\lambda^{-1}D),a]b\|_{L^{p}}\leq\lambda^{-1}\|\,|\cdot|h\|_{L^{1}}|\n a\|_{L^{\infty}}\|b\|_{L^{p}}.$$
This conclude the proof of the lemma.
\null{\hfill $\Box$}\\
\\
Inequality (\ref{abc13}) is a consequence of the classical following lemma \ref{aKlemme3}  which is also proved in \cite{BCD} p 112. For the completeness of the proof we will give here his proof.
\begin{lemme}
\label{aKlemme3}
Let $\sigma\in\R$, $1\leq r\leq+\infty$ and $1\leq p_{1}\leq p\leq+\infty$. Let $u$ a vector field over $\R^{N}$. Assume that:
$$\sigma >-N\inf(\frac{1}{p^{'}},\frac{1}{p_{1}})\;\;\mbox{or}\;\;\sigma>-1-N\inf(\frac{1}{p^{'}},\frac{1}{p_{1}})\;\;\mbox{if}\;\;{\rm div}u=0.$$
Denote $R_{q}=[u\cdot\n,\D_{q}]a$. There exists a constant $C$ depending on $p,p_{1},\sigma$ and $N$ such that:
\begin{equation}
\|(2^{q\sigma}\|R_{q}\|_{L^{p}})_{q}\|_{l^{r}}\leq C\|a\|_{B^{\sigma}_{p,r}}\|u\|_{B^{\NNN}_{p_{1},\infty}\cap L^{\infty}}\;\;\;\mbox{if}\;\;\sigma<\NNN+1.
\label{57}
\end{equation}
and in the critical case:
\begin{equation}
\|(2^{q\sigma}\|R_{q}\|_{L^{p}})_{q}\|_{l^{r}}\leq C\|a\|_{B^{\sigma}_{p,r}}\|u\|_{B^{\NNN}_{p_{1},1}}\;\;\;\mbox{if}\;\;\sigma=\NNN+1.
\label{58}
\end{equation}
In the limit case $\sigma=-N\inf(\frac{1}{p^{'}},\frac{1}{p_{1}})$ (or  $\sigma=-1-N\inf(\frac{1}{p^{'}},\frac{1}{p_{1}})$ for ${\rm div}u=0$), then we have:
\begin{equation}
\sup_{q}2^{q\sigma}\|R_{q}\|_{L^{p}}\leq C\|a\|_{B^{\sigma}_{p,\infty}}\|\n u\|_{B^{\NNN}_{p_{1},1}}.
\label{59}
\end{equation}
\end{lemme}
{\bf Proof:}
In order to show that only the gradient part of $u$ is involved in the estimates, we  shall decompose into low and high frequencies $u=S_{0}u+u_{1}$. We remind that there exists a constant $C$ such that:
\begin{equation}
\forall p_{2}\in [1,+\infty], \;\;\|S_{0}\n u\|_{L^{p_{2}}}\leq C\|\n u\|_{L^{p_{2}}}\;\;\mbox{and}\;\;
\|\n u_{1}\|_{L^{p_{2}}}\leq C\|\n u_{1}\|_{L^{p_{2}}}.
\label{d257}
\end{equation}
Furthermore as $u_{1}$ is spectrally supported away form the origin, by the so-called Bernstein«s inequalities, we obtain that:
\begin{equation}
\forall p_{2}\in [1,+\infty], \;\;
\frac{2^{q}}{C}  \|\D_{q} u_{1}\|_{L^{p_{2}}}\leq\|\D_{q}\n u_{1}\|_{L^{p_{2}}}\leq C  2^{q}\|\D_{q} u_{1}\|_{L^{p_{2}}}.
\label{d258}
\end{equation}
Now we can write $R_{q}$ under the following form:
$$
\begin{aligned}
R_{q}&=u\cdot\n\D_{q}a-\D_{q}(u\cdot\n a),\\
&=[u_{1}^{k},\D_{q}]\p_{k}a+[S_{0}u^{k},\D_{q}]\p_{k}a.
\end{aligned}
$$
The proof  is based on Bony's decomposition which enables us to split $R_{q}$
into $R_{q}=\sum_{i=1}^{8}R^{i}_{q},$
where:
$$
\begin{aligned}
&R^{1}_{q}=[T_{u_{1}^{k}},\D_{q}]\p_{k}a,\hspace{3cm}R^{2}_{q}=T_{\p_{k}\D_{q}}u_{1}^{k},\\
&R^{3}_{q}=-\D_{q}T_{\p_{k}a}u_{1}^{k},\hspace{3,1cm}R^{4}_{q}=\p_{k}R(u_{1}^{k},\D_{q}a),\\
&R^{5}_{q}=-R({\rm div}u_{1},\D_{q}a),\hspace{2,3cm}R^{6}_{q}=-\p_{k}\D_{q}R(u_{1},a),\\
&R^{7}_{q}=\D_{q}R({\rm div}u_{1},a),\hspace{2,6cm}R^{8}_{q}=[S_{0}u^{k},\D_{q}]\p_{k}a.\\
\end{aligned}
$$
In the sequel we will denote $(c_{q})_{q\geq-1}$ a sequence such that $\|(c_{q})\|_{l^{r}}\leq 1$.
\subsubsection*{Bounds for $2^{q\sigma}\|R^{1}_{q}\|_{L^{p}}$}
Using proposition \ref{d210} we obtain:
$$R^{1}_{q}=\sum_{|q-q^{'}|\leq 4}[S_{q^{'}-1}u_{1}^{k},\D_{q},]\p_{k}\D_{q^{'}}a$$
Hence according to lemma \ref{d297} and inequality (\ref{d257}) we obtain:
\begin{equation}
\begin{aligned}
2^{q\sigma}\|R^{1}_{q}\|_{L^{p}}&\leq C\|\n v\|_{L^{\infty}}\sum_{|q-q^{'}|\leq 4}2^{q^{'}\sigma}\|\D_{q^{'}}a\|_{L^{p}},\\
&\leq Cc_{j}\|\n v\|_{L^{\infty}}\|a\|_{B^{\sigma}_{p,r}}.
\end{aligned}
\label{d259}
\end{equation}
\subsubsection*{Bounds for $2^{q\sigma}\|R^{2}_{q}\|_{L^{p}}$}
From proposition \ref{d210}, we have:
$$
R^{2}_{q}=\sum_{q^{'}\geq q-3}S_{q^{'}-1}\p_{k}\D_{q}a \;\D_{q^{'}}u_{1}^{k}
$$
Hence using (\ref{d257}) and (\ref{d258}) we have:
\begin{equation}
2^{q\sigma}\|R^{2}_{q}\|_{L^{p}}\leq Cc_{j}\|\n v\|_{L^{\infty}}\|a\|_{B^{\sigma}_{p,r}}.
\label{d260}
\end{equation}
\subsubsection*{Bounds for $2^{q\sigma}\|R^{3}_{q}\|_{L^{p}}$}
One has:
\begin{equation}
\begin{aligned}
R^{3}_{q}&=-\sum_{|q-q^{'}|\leq 4}\D_{q}(S_{q^{'}-1}\p_{k}a\, \D_{q^{'}}u_{1}^{k}),\\
&=-\sum_{|q-q^{'}|\leq 4, q^{''}\leq q^{'}-2}\D_{q}(\D_{q^{''}}\p_{k}a\, \D_{q^{'}}u_{1}^{k}).
\end{aligned}
\label{d261}
\end{equation}
Therefore denoting $\frac{1}{p_{3}}=\frac{1}{p}-\frac{1}{p_{1}}$ and using (\ref{d257}) and (\ref{d258}),
$$
\begin{aligned}
2^{q\sigma}\|R^{3}_{q}\|_{L^{p}}&\leq C\sum_{|q-q^{'}|\leq 4, q^{''}\leq q^{'}-2}2^{q\sigma} \|\D_{q^{''}}\p_{k}a\|_{L^{p_{3}}}\, \|\D_{q^{'}}u_{1}^{k}\|_{L^{p_{1}}},\\
&\leq C\sum_{|q-q^{'}|\leq 4, q^{''}\leq q^{'}-2}2^{(q-q^{''})(\sigma-1-\NNN)} \|\D_{q^{''}}a\|_{L^{p}}2^{q^{'}\NNN}\|\D_{q^{'}}\n u\|_{L^{p_{1}}},
\end{aligned}
$$
Hence if $\sigma<1+\NNN$,
\begin{equation}
2^{q\sigma}\|R^{3}_{q}\|_{L^{p}}\leq Cc_{j}\|\n v\|_{B^{\NNN}_{p_{1},\infty}}\|a\|_{B^{\sigma}_{p,r}}.
\label{d263}
\end{equation}
\subsubsection*{Bounds for $2^{q\sigma}\|R^{4}_{q}\|_{L^{p}}$}
Letting $\widetilde{\D}_{q^{'}}=\D_{q^{'}-1}+\D_{q^{'}}+\D_{q^{'}+1}$, we have:
$$R^{4}_{q}=\sum_{|q-q^{'}|\leq 2}\p_{k}(\D_{q^{'}}u_{1}^{k}\D_{q}\widetilde{\D}_{q^{'}}a).$$
Hence by virtue of ({\ref{d258}), we get:
\begin{equation}
2^{q\sigma}\|R^{4}_{q}\|_{L^{p}}\leq Cc_{j}\|\n v\|_{B^{\NNN}_{p_{1},\infty}}\|a\|_{B^{\sigma}_{p,r}}.
\label{d265}
\end{equation}
We follow the same lines for bounding  $2^{q\sigma}\|R^{5}_{q}\|_{L^{p}}$.
\subsubsection*{Bounds for $2^{q\sigma}\|R^{6}_{q}\|_{L^{p}}$ and $2^{q\sigma}\|R^{7}_{q}\|_{L^{p}}$}
We will begin with treating the case $\frac{1}{p}+\frac{1}{p_{1}}\leq 1$. Let $\frac{1}{p_{4}}=\frac{1}{p}+\frac{1}{p_{1}}$. Then if $\sigma>-1-\NNN$ proposition \ref{produit1} and proposition \ref{interpolation} yield:
\begin{equation}
2^{q\sigma}\|R^{6}_{q}\|_{L^{p}}\leq Cc_{j}\|\n u_{1}\|_{B^{\NNN}_{p_{1},\infty}}\|a\|_{B^{\sigma}_{p,r}}.
\label{d266}
\end{equation}
Now if  $\frac{1}{p}+\frac{1}{p_{1}}> 1$, the previous argument has to be applied with $p^{'}$ instead of $p_{3}$ and one still gets (\ref{d266}) provided that $\sigma>-1-\frac{N}{p^{'}}$. One has then:
\begin{equation}
2^{q\sigma}\|R^{6}_{q}\|_{L^{p}}\leq Cc_{j}\|\n u_{1}\|_{B^{\NNN}_{p_{1},\infty}}\|a\|_{B^{\sigma}_{p,r}}.
\label{d267}
\end{equation}
Note that in the limit case $\sigma=-1-\min(\NNN,\frac{N}{p^{'}})$, proposition \ref{produit3} yields:
\begin{equation}
\sup_{q}2^{q\sigma}\|R^{6}_{q}\|_{L^{p}}\leq Cc_{j}\|\n u_{1}\|_{B^{\NNN}_{p_{1},1}}\|a\|_{B^{\sigma}_{p,r}}.
\label{d268}
\end{equation}
Similar arguments allows us to obtain:
\begin{equation}
\begin{aligned}
&2^{q\sigma}\|R^{7}_{q}\|_{L^{p}}\leq Cc_{j}\|\n u_{1}\|_{B^{\NNN}_{p_{1},\infty}}\|a\|_{B^{\sigma}_{p,r}}\;\;\;\;\mbox{if}\;\;\sigma>-\min((\NNN,\frac{N}{p^{'}}),\\
&2^{q\sigma}\|R^{7}_{q}\|_{L^{p}}\leq Cc_{j}\|\n u_{1}\|_{B^{\NNN}_{p_{1},1}}\|a\|_{B^{\sigma}_{p,r}}\;\;\;\;\mbox{if}\;\;\sigma=-\min((\NNN,\frac{N}{p^{'}})\;\;r=+\infty.
\end{aligned}
\label{d269}
\end{equation}
Finally let mention that if $\sigma>-1$ then the standard continuity results for the remainder combined with the embedding $L^{\infty}\h B^{0}_{\infty,\infty}$ yield:
\begin{equation}
2^{q\sigma}\|R^{6}_{q}\|_{L^{p}}\leq Cc_{j}\|\n u\|_{L^{\infty}}\|a\|_{B^{\sigma}_{p,r}}.
\label{d271}
\end{equation}
The same inequality holds true for $R^{7}_{q}$ if $\sigma>0$.
\subsubsection*{Bounds for $2^{q\sigma}\|R^{8}_{q}\|_{L^{p}}$}
We have then:
$$R^{8}_{q}=-\sum_{|q-q^{'}|\leq 1}[\D_{q},\D_{-1}u]\cdot\n \D_{q^{'}}a.$$
From lemma \ref{d297} we obtain:
\begin{equation}
\begin{aligned}
2^{q\sigma}\|R^{8}_{q}\|_{L^{p}}&\leq C\sum_{|q-q^{'}|\leq 1}\|\n \D_{-1}u\|_{L^{\infty}}2^{q^{'}\sigma}
 \|\D_{q^{'}}a\|_{L^{p}},\\
& \leq Cc_{j}\|\n u\|_{L^{\infty}}\|a\|_{B^{\sigma}_{p,r}}.
\end{aligned}
\label{d272}
\end{equation}
Combining inequalities (\ref{d259}), (\ref{d260}), (\ref{d263}),  (\ref{d265}), (\ref{d267}),  (\ref{d268}), (\ref{d269}), (\ref{d271}) and (\ref{d272}) yields (\ref{57}), (\ref{58}), and (\ref{59}).
\null{\hfill $\Box$}

\end{document}